\documentclass[final]{siamltex}
\usepackage{mathrsfs}
\usepackage{amsfonts}
\usepackage{graphicx}
\usepackage{subfigure}
\usepackage{xcolor}
\usepackage{amssymb}
\usepackage{amsmath}
\DeclareMathOperator{\sech}{sech}
 \usepackage{epstopdf}
 \usepackage{boondox-cal}
\usepackage{mathrsfs,color}
\usepackage{amsfonts}
\usepackage{multirow}
\usepackage{amsmath}
\usepackage{graphicx}
\usepackage{fancybox}
\usepackage{fancyhdr}
\usepackage{fancyvrb}
\usepackage{graphicx}
\usepackage{subfigure}
\usepackage{txfonts}
\usepackage{dsfont}
\usepackage{txfonts}
\usepackage{amssymb}
\usepackage{amsmath}
 \usepackage{epstopdf}
\usepackage{amsmath,amssymb}
\usepackage{mathrsfs,color}
\usepackage[linkcolor=red,pdfborder=001,anchorcolor=gree,citecolor=blue]{hyperref}
\usepackage{pifont}
\usepackage[misc]{ifsym}

\newtheorem{example}{Example}

\begin{document}

%




\title{Mathematical analysis and numerical simulation
of coupled nonlinear space-fractional Ginzburg-Landau 
equations\thanks{This work was partially supported by
 the National Natural Science Foundation of China
 (Grant Nos. 12461069, 11961057), the Science and
technology project of Guangxi (Grant No. GuikeAD21220114).}}

\author{Hengfei Ding\textsuperscript{\Letter}\thanks{1.School of Mathematics and Statistics,
Guangxi Normal University, Guilin 541006, China;
2.The Center for Applied
Mathematics of Guangxi (GXNU), Guilin 541006, China;
3.Guangxi Colleges and Universities Key Laboratory of Mathematical Model and Application (GXNU),
 Guilin 514006, China. (\textsuperscript{\Letter} E-mail:dinghf05@163.com).}
    \and Yuxin Zhang\thanks{School of Mathematics and Statistics,
Guangxi Normal University, Guilin 541006, China. (E-mail: zhangyuxin2006@163.com).}\
    \and Qian Yi\textsuperscript{\Letter}\thanks{School of Mathematics and Statistics,
Guangxi Normal University, Guilin 541006, China. (\textsuperscript{\Letter}
E-mail: yiqianqianyi@163.com).}}

\maketitle
\begin{abstract}
The coupled nonlinear space fractional Ginzburg-Landau (CNLSFGL) equations
with the fractional Laplacian have been widely used to model the dynamical processes in
 a fractal media with fractional dispersion.
Due to the existence of fractional power derivatives and strong nonlinearity, it is extremely
difficult to mathematically analyze the CNLSFGL equations and construct efficient numerical algorithms.
For this reason, this paper aims to investigate the theoretical results about the considered
system and construct a novel high-order numerical scheme for this coupled system. We prove
 rigorously an a priori estimate of the solution to the coupled system and the well-posedness of its weak solution.
Then, to develop the efficient numerical algorithm, we construct a fourth-order numerical
differential formula to approximate the fractional Laplacian. Based on this formula, we
 construct a high-order implicit difference scheme for the coupled system. Furthermore,
 the unique solvability and convergence of the established algorithm are proved in detail.
 To implement the implicit algorithm efficiently, an iterative algorithm is designed in
 the numerical simulation. Extensive numerical examples are reported to further
  demonstrate the correctness of the theoretical analysis and the efficiency of the proposed numerical algorithm.
\end{abstract}



\pagestyle{myheadings} \thispagestyle{plain} \markboth{Hengfei Ding}{The
 coupled nonlinear Ginzburg-Landau equations with the
fractional Laplacian}

\section{Introduction}

Generally, we regard fractional derivatives as a generalization and extension of
classical derivatives, which can compensate the drawbacks of integer
derivatives in characterizing some complex processes or systems.
In particular, the nonlocality of fractional derivatives makes them
perfectly suitable for describing materials or processes with memory
and genetic properties, enabling accurate descriptions and clearer
physical interpretations for complex systems.
Due to these advantages, the fractional differential equations based
on the fractional derivatives have been widely used in physics,
biology, engineering, finance, and various interdisciplinary fields \cite{Podlubny,Metzler,Hilfer,Ortigueira}.

In this paper, we focus on the following coupled  nonlinear space fractional
 Ginzburg-Landau (CNLSFGL) equations with the
fractional Laplacian
\begin{equation}\label{eq.1}
\left\{
\begin{aligned}\displaystyle
\partial_t u+\left(\beta_1+\mathrm{i}\eta_1\right)
\left(-\triangle\right)^\frac{\alpha}{2}_\Omega u
+\left(\mu_1+\mathrm{i}\zeta_1\right)|u|^2u-\gamma_1 u-\mathrm{i}|u|^2v
=0,\;(x,t)\in\Omega\times(0,T],\\
\partial_t v+\left(\beta_2+\mathrm{i}\eta_2\right)
\left(-\triangle\right)^\frac{\alpha}{2}_\Omega v
+\left(\mu_2+\mathrm{i}\zeta_2\right)|v|^2v-\gamma_2 v-\mathrm{i}|v|^2u
=0,\;(x,t)\in\Omega\times(0,T],\;
\end{aligned}
\right.
\end{equation}
with the initial conditions
\begin{equation}\label{eq.2}
u(x,0)=u_0(x),\;v(x,0)=v_0(x),\;x\in\Omega,
\end{equation}
and the periodic boundary conditions
\begin{equation}\label{eq.3}
u(a,t)=u(b,t),\;v(a,t)=v(b,t),\;t\in[0,T],
\end{equation}
where $\mathrm{i}=\sqrt{-1}$,
$u=u(x,t)$ and $v=v(x,t)$ are complex-valued functions of a space variable $x \in \Omega =[a,b]\subset{\mathds{R}}$
 and of a time variable $t \in [0, T ]$, $\beta_s>0,\eta_s>0,\mu_s\geq0$ and $\zeta_s>0$
 are all real parameters for $s=1,2$. $\gamma_1$ and $\gamma_2$ are two real-valued parameters.
  The fractional Laplacian operator $\left(-\triangle\right)^\frac{\alpha}{2}_\Omega$ with $1 < \alpha \leq2$
is defined via the following hypersingular integral \cite{Samko,Landkof}
\begin{equation}\label{eq.4}
\left(-\triangle\right)^\frac{\alpha}{2}_{\Omega}u(x)=K_\alpha\int_\Omega
\frac{u(x)-u(y)}{|x-y|^{1+\alpha}}\mathrm{d}y,
\;\;K_\alpha=\frac{2^\alpha\Gamma\left(\frac{1+\alpha}{2}\right)}
{\pi^{\frac{1}{2}}\left|\Gamma\left(-\frac{\alpha}{2}\right)\right|}.
\end{equation}

The classical nonlinear Ginzburg-Landau (NLGL) equation is one of the most extensively
studied nonlinear evolution equations in physics, which can be used to explain and
 describe various complex physical phenomena, ranging from nonlinear waves to
 second-order phase transitions, and from superconductivity, superfluidity and
  Bose-Einstein condensation to liquid crystal chords in field theory, and so on \cite{Aranson,Du,Ginzburg,Pan}.
The nonlinear space fractional Ginzburg-Landau (NLSFGL) equation can be regarded as a generalization of the
classical NLGL equation, which was first obtained by Tarasova and Zaslavsky in deducing the Euler-Lagrange
equation of fractal substances \cite{Tarasov1,Tarasov2}.
Compared with the classical counterpart, the fractional equation has more extensive applications. For example, it can be
used to describe the nonlinear dynamics of diffuse Hindmarsh-Rose neural network with long-range
 coupling \cite{Mvogo}, the dynamics in media with fractal mass dimension or
 fractal dispersion characteristics \cite{Tarasov},
 and the critical phenomenon transition point affected by phase nonlocal effect \cite{Milovanov}.
 Among these, the coupled NLSFGL equations are often used to demonstrate the phenomenon
of interference between two different beams of light in a nonlinear medium when
the refractive index changes.

So far, there have been a series of works devoted to study and analyze the NLGL equation from
the theoretical point of view, please refer to \cite{Duan,Doering,Gao,Guo,Huo}. On the contrary, the theoretical results
of NLSFGL equations are relatively scarce, and the fundamental reason lies in the fractional
power of the Laplacian, which makes its theoretical analysis more difficult.
In \cite{Guo}, the
global well-posedness and long-time dynamics were studied by
Pu and Guo. Later, Gu et al. also investigated
the well-posedness of the NLSFGL equation in \cite{Gu}. The
asymptotic analysis of the NLSFGL equation in bounded domains were given in \cite{Millot} by Millot and Sire.
Meanwhile, they proposed a general
partial regular result in an arbitrary dimension.
Lu and L\"{u} \cite{Lu} considered the well-posedness and asymptotic behaviors of solutions to the NLSFGL equation with
the periodic boundary conditions.
It is worth mentioning that as far as we know, there is
no study on the well-posedness of the CNLSFGL equations with the fractional
Laplacian.

Generally,
it is
difficult to obtain the analytical solution of the NLSFGL equation due
 to the nonlocal property of the fractional derivatives and the nonlinear nature of the
  equation itself. Therefore, we must resolve it by numerical methods.
In recent years,
there are quite a lot of different numerical methods for the NLSFGL equation.
As far as the finite difference method is concerned, Wang and Huang \cite{Wang}
established an implicit midpoint
difference scheme with second-order convergence in time and space.
Subsequently, with the help of a known fourth-order numerical differential formula that approximates the
Riesz derivatives, they proposed a three-level implicit-explicit difference
scheme to solve the NLSFGL equation.
In \cite{Hao}, based on the compact fractional centered difference formula in
space and the Crank-Nicolson/leap-frog method in time,
a linearized three-level difference scheme has been derived.
Recently, a series of finite difference methods driven by similar core idea
 have also been developed;
see, e.g., \cite{He,Wang1,Ding,Mohebbi,Zhang1,Zhang2} and reference therein.
In addition, a number of research teams have also proposed some other algorithms,
such as the Galerkin finite element \cite{Li,Zhang3},
the Galerkin-Legendre spectral \cite{Fei}, the radial basis functions \cite{Shokri},
the Petrov-Galerkin \cite{Abbaszadeh}
and Fourier
spectral \cite{Zeng,Lu1} techniques to deal with NLSFGL equations.
Here, it should be emphasized that compared with the two-level schemes (Crank-Nicolson method),
the linearized three-level schemes can avoid iterative calculation. However,  its defect is also obvious,
that is, we need to use some other methods
 to calculate the value of the second-level, which will lead to
 a coupled difference scheme and bring some difficulties to the
  theoretical analysis and calculation.

From the above introduction and analysis, it is easy to see that some progress has been made
in the numerical research on the NLSFGL equation. However,
there are limited theoretical analysis and numerical
studies concerning the CNLSFGL equations (\ref{eq.1}). As far as we are aware,
only in \cite{Li&Huang} have researchers proposed
a second-order difference scheme based
on the implicit midpoint method in time and a weighted-shifted Gr\"{u}nwald formula in space.
Given the foregoing, there remains a need to
investigate the well-posedness of the weak solution and
develop higher-order numerical methods
 for the CNLSFGL equations.
 The main contributions in this paper are listed specifically as follows:
\begin{itemize}
\item The well-posedness of the weak solutions to system (\ref{eq.1}) is
studied. Compared with a single NLSFGL equation, the difficulty of the well-posedness analysis lies in how to deal
with the coupling terms $-\mathrm{i}|u|^2v$ and $-\mathrm{i}|v|^2u$ effectively. Here, we successfully solve
these problems by establishing an a priori estimate.

\item A new efficient fourth-order numerical differential formula
approximating the fractional Laplacian is proposed. We obtain an efficient numerical differential
  formula by an indirect method based on the equivalence between the fractional Laplacian and Riesz derivative
   under certain conditions.
In this case, choosing a reasonable generating function is the core of constructing the
higher-order numerical differential formula for the Riesz derivative. Here, we have solved
the most critical problem.

\item It is strictly proved that the generating function $G_2(z)$ proposed in
\cite{Ding&Li} is analyzable on $[-1,1]$. Previously, the generating function $G_2(z)$
 was studied only in the interval $(-1,1)$, and here, we prove that its convergence interval is $[-1,1]$.

\item A two-level efficient numerical algorithm is constructed for system (\ref{eq.1}).
Different from the existing linearized three-level difference schemes, we use
the compact implicit integration factor method and Pad\'{e} approximation strategy
for the
time discretization, while the spatial discretization is based on the fourth-order numerical
differential formula constructed in the paper. The numerical algorithm established in this
way is a two-level scheme with second-order and fourth-order accuracy in time and space, respectively.
The most important advantage of the proposed two-level scheme lies in that we do not
need to use some other numerical methods to compute the value at time level $t_1$,
but achieve 2nd-order in time as well. Moreover, the two-level scheme will not increase the complexity
  of subsequent theoretical analysis and calculations.

  \item
An efficient iterative algorithm is constructed and theoretically analyzed for
the implicit difference scheme established in the paper. The effective iterative algorithm
is designed based on the
characteristics of the established algorithm, which has the advantage
of maintaining the same convergence order as the original difference schemes. Finally,
we prove that this iterative algorithm is convergent.

\end{itemize}

The reminder of this paper is organized as follows.
In Section 2, we analyze the well-posedness of the weak solution to (\ref{eq.1})--(\ref{eq.3}).
Section 3 presents a new effective method for numerical solving
(\ref{eq.1})--(\ref{eq.3}). Moreover, we discuss the numerical method's boundness, unique solvability and convergence.
In Section 4, we validate
the theoretical results through several numerical examples and demonstrate
the effectiveness of the proposed numerical
method. Finally, the paper is concluded with a summary in the last section.

\section{Mathematical analysis}
In this section, our main purpose is to discuss the well-posedness
of system (\ref{eq.1})--(\ref{eq.3}). For the sake of simplicity,
let $\Omega=[0,2\pi]$ denote the periodic interval with periodic $2\pi$ and
$\mathds{Z}$ denote the set of all integers.
Hence, function $u$ and $v$ can be expressed by a Fourier series
$$u=\sum_{s\in\mathds{Z}}u_{s}e^{\mathrm{i}sx},\;
v=\sum_{s\in\mathds{Z}}v_{s}e^{\mathrm{i}sx},$$
where, $u_s$ and $v_s$ are the Fourier coefficients of $u$ and $v$, respectively.
Then we further know that \cite{Lu}
\begin{equation*}
\left(-\triangle\right)^\frac{\alpha}{2}_{\Omega}u(x)=\sum_{s\in\mathds{Z}}
\left|s\right|^\alpha u_{s}e^{\mathrm{i}sx},\;\;
\left(-\triangle\right)^\frac{\alpha}{2}_{\Omega}v(x)=\sum_{s\in\mathds{Z}}
\left|s\right|^\alpha v_{s}e^{\mathrm{i}sx}.
\end{equation*}

Let $H^{\alpha}\left(\Omega\right)$ denote the complete Sobolev space of order
$\alpha$ under the following norm
\begin{equation*}
\begin{aligned}\displaystyle
\left\|u\right\|_{H^\alpha\left(\Omega\right)}=
\left(\sum_{s\in\mathds{Z}}\left|u_{s}\right|^2\left(1+\left|s\right|^{2\alpha}\right)\right)^{\frac{1}{2}}.
\end{aligned}
\end{equation*}
Meanwhile, we denote by $\left\|\cdot\right\|_{L^2(\Omega)}$
 the norm of $L^2(\Omega)$ with usual
inner product $(\cdot,\cdot)$, $\left\|\cdot\right\|_{L^p(\Omega)}$
the norm of the standard Lebesgue space $L^p(\Omega)$ for all $1 \leq p \leq\infty$,
and ${v^*}$  the complex conjugate of $v$. In addition,
$\Re(\cdot)$ and $\Im(\cdot)$ represent
 the real part and the imaginary part, respectively.
 For the convenience of subsequent analysis, we first provide
 several lemmas that play a key role in relevant analysis.

\begin{definition}(See \cite{Guo1})
Let $X$ be a Banach space with norm $\left\|\cdot\right\|_X$. Denote by $L^p\left(0,T;X\right)$ the
space of all measurable function $u:[0,T]\rightarrow X$ with the norm
\begin{equation*}
\begin{aligned}\displaystyle
\left\|u\right\|_{L^p\left(0,T;X\right)}=
\left(\int_{0}^{T}\left\|u\right\|_{X}^p\mathrm{d}t
\right)^{1/p}<\infty,\;\;1\leq p<\infty,
\end{aligned}
\end{equation*}
and when $p=\infty$,
$$
\left\|u\right\|_{L^\infty\left(0,T;X\right)}=
{\sup}_{0\leq t\leq T}\left\|u\right\|_X<\infty.
$$
\end{definition}

\begin{lemma}\label{Le.2.2}
(Differential form of Grownall's inequality \cite{Evans})
Let $u(t)$ be a nonnegative, absolutely continuous
function on $[0,T]$, which satisfies for a.e. $t$, the
differential inequality holds
$$\frac{\mathrm{d}u(t)}{\mathrm{d}t}\leq
c(t)u(t)+f(t),$$
where $c(t)$ and $f(t)$ are nonnegative, summable
 functions on $[0,T]$. Then there holds that
$$u(t)\leq e^{\int_{0}^{t}c(s)\mathrm{d}s}
\left[u(0)+\int_{0}^{t}f(s)\mathrm{d}s\right]$$
for all $0\leq t\leq T$.
\end{lemma}


\begin{lemma}\label{Le.2.3}(See \cite{Macias&Diaz})
The fractional Laplace operator  $\left(-\Delta\right)^{\frac{\alpha}{2}}_\Omega$
 has a unique square-root operator over the space of
sufficiently regular function with compact support. That is to say, there
exists an operator $\left(-\Delta\right)^{\frac{\alpha}{4}}_\Omega $,
such that
\begin{eqnarray*}
\begin{array}{lll}
\displaystyle
\int_{\Omega}\left(\left(-\Delta\right)^{\frac{\alpha}{2}}_\Omega u(x)\right)v(x)\mathrm{d}x
=\int_{\Omega}\left(\left(-\Delta\right)^{\frac{\alpha}{4}}_\Omega u(x)\right)
\left(\left(-\Delta\right)^{\frac{\alpha}{4}}_\Omega v(x)\right)
\mathrm{d}x.
\end{array}
\end{eqnarray*}
\end{lemma}

\begin{lemma}\label{Le.2.4}(See \cite{Ding1})
For any complex function $U, V, u, v $, there holds that
\begin{equation*}
\begin{aligned}
\left|\left|U\right|^2V-\left|u\right|^2v\right| \leq C_{\max}
\left(\left|U-u\right|+\left|V-v\right|\right),
\end{aligned}
\end{equation*}
where $C_{\max}=\max\left\{2\left|U\right|\left|V\right|,2\left|u\right|\left|V\right|,\left|u\right|^2\right\}$.
\end{lemma}

\subsection{A priori estimate}
In this section, we consider the a priori estimate of the system (\ref{eq.1})--(\ref{eq.3}).
\begin{theorem}\label{Th.2.5}
Let $w=(u,v)$ be the solution to the system (\ref{eq.1})
with initial data $w_0=\left(u,v\right)|_{t=0}=\left(u_0,v_0\right)\in L^2(\Omega)$,
then there hold
\begin{equation}\label{eq.2.1}
\begin{aligned}\displaystyle
\left\|u\right\|_{L^2{(\Omega)}}^2+\left\|v\right\|_{L^2{(\Omega)}}^2\leq
\exp\left(2t\max\left\{|\gamma_1|,|\gamma_2|\right\}\right)
\left(\left\|u_0\right\|_{L^2{(\Omega)}}^2+\left\|v_0\right\|_{L^2{(\Omega)}}^2\right),\;\;
0<t\leq T.
\end{aligned}
\end{equation}
and
\begin{equation}\label{eq.2.2}
\begin{aligned}\displaystyle
&\left\|u\right\|_{L^2{(\Omega)}}^2+\left\|v\right\|_{L^2{(\Omega)}}^2
+2\beta_1
\int_{0}^t
\left\|\left(-\triangle\right)^\frac{\alpha}{4}_\Omega u\right\|_{L^2{(\Omega)}}^2\mathrm{d}s
+
2\beta_2
\int_{0}^t
\left\|\left(-\triangle\right)^\frac{\alpha}{4}_\Omega v\right\|_{L^2{(\Omega)}}^2\mathrm{d}s\\
\leq&\exp\left(2t\max\left\{\left|\gamma_1\right|,\left|\gamma_2\right|\right\}\right)
\left(\left\|u_0\right\|_{L^2{(\Omega)}}^2+\left\|v_0\right\|_{L^2{(\Omega)}}^2\right),
\;\;0<t\leq T.
\end{aligned}
\end{equation}
\end{theorem}

\begin{proof}
Taking the inner product of the first equation of system (\ref{eq.1}) with $u$ and using Lemma \ref{Le.2.3},
then taking the real part can result in
 \begin{equation}\label{eq.2.3}
\displaystyle
\Re\left\{\int_{\Omega}\partial_t uu^{*}\mathrm{d}x\right\}+\beta_1
\int_{\Omega}
\left|\left(-\triangle\right)^\frac{\alpha}{4}_\Omega u\right|^2\mathrm{d}x
+\mu_1\int_{\Omega}|u|^4\mathrm{d}x
-\gamma_1\int_{\Omega} \left|u\right|^2\mathrm{d}x=\Im\left\{\int_{\Omega}|u|^2vu^{*}\mathrm{d}x\right\}
.
\end{equation}
For the second equation of system (\ref{eq.1}), using the same method, we can also obtain
 \begin{equation}\label{eq.2.4}
\displaystyle
\Re\left\{\int_{\Omega}\partial_t vv^{*}\mathrm{d}x\right\}+\beta_2
\int_{\Omega}
\left|\left(-\triangle\right)^\frac{\alpha}{4}_\Omega v\right|^2\mathrm{d}x
+\mu_2\int_{\Omega}|v|^4\mathrm{d}x
-\gamma_2\int_{\Omega} \left|v\right|^2\mathrm{d}x=\Im\left\{\int_{\Omega}|v|^2uv^{*}\mathrm{d}x\right\}
.
\end{equation}
Using the Young's inequality, we can know that
\begin{equation*}
\begin{aligned}\displaystyle
\Im\left\{\int_{\Omega}|u|^2vu^{*}\mathrm{d}x\right\}
\leq&\int_{\Omega}|u|^2|v||u^{*}|\mathrm{d}x
\leq\int_{\Omega}|u|^2\left(\varepsilon_1^2|v|^2
+\frac{1}{4\varepsilon_1^2}|u|^2\right)\mathrm{d}x\\
\leq&\int_{\Omega}\left[
\frac{1}{2}\varepsilon_1^2\left(|u|^4+|v|^4\right)
+\frac{1}{4\varepsilon_1^2}|u|^4\right]\mathrm{d}x
\end{aligned}
\end{equation*}
and
\begin{equation*}
\begin{aligned}\displaystyle
\Im\left\{\int_{\Omega}|v|^2uv^{*}\mathrm{d}x\right\}
\leq\int_{\Omega}\left[
\frac{1}{2}\varepsilon_2^2\left(|v|^4+|u|^4\right)
+\frac{1}{4\varepsilon_2^2}|v|^4\right]\mathrm{d}x.
\end{aligned}
\end{equation*}
Taking
\begin{equation*}
\left\{
\begin{aligned}
\frac{1}{2}\left(\varepsilon_1^2+\varepsilon_2^2\right)
+\frac{1}{4\varepsilon_1^2}=\mu_1,\\
\frac{1}{2}\left(\varepsilon_1^2+\varepsilon_2^2\right)
+\frac{1}{4\varepsilon_2^2}=\mu_2,
\end{aligned}
\right.
\end{equation*}
and adding (\ref{eq.2.3}) and (\ref{eq.2.4}) and further simplifying can lead to
\begin{equation}\label{eq.2.5}
\begin{aligned}\displaystyle
&\frac{\mathrm{d}}{\mathrm{d}t}\int_{\Omega}\left(\left| u\right|^2+\left| v\right|^2\right)\mathrm{d}x
+2\beta_1
\int_{\Omega}
\left|\left(-\triangle\right)^\frac{\alpha}{4}_\Omega u\right|^2\mathrm{d}x
+
2\beta_2
\int_{\Omega}
\left|\left(-\triangle\right)^\frac{\alpha}{4}_\Omega v\right|^2\mathrm{d}\\
\leq&
2\gamma_1\int_{\Omega} \left|u\right|^2\mathrm{d}x
+
2\gamma_2\int_{\Omega} \left|v\right|^2\mathrm{d}x
.
\end{aligned}
\end{equation}

By discarding some positive terms and with the help of Lemma \ref{Le.2.2}, we can obtain
the desired result as follows
\begin{equation*}
\begin{aligned}\displaystyle
\left\|u\right\|_{L^2{(\Omega)}}^2+\left\|v\right\|_{L^2{(\Omega)}}^2\leq
\exp\left(t\max\left\{2\gamma_1,2\gamma_2\right\}\right)
\left(\left\|u_0\right\|_{L^2{(\Omega)}}^2+\left\|v_0\right\|_{L^2{(\Omega)}}^2\right),\;\;
0<t\leq T.
\end{aligned}
\end{equation*}

Furthermore, from (\ref{eq.2.5}), we can easily get
\begin{equation}\label{eq.2.6}
\begin{aligned}\displaystyle
&\left\|u\right\|_{L^2{(\Omega)}}^2+\left\|v\right\|_{L^2{(\Omega)}}^2
+2\beta_1
\int_{0}^t
\left\|\left(-\triangle\right)^\frac{\alpha}{4}_\Omega u\right\|_{L^2{(\Omega)}}^2\mathrm{d}s
+
2\beta_2
\int_{0}^t
\left\|\left(-\triangle\right)^\frac{\alpha}{4}_\Omega v\right\|_{L^2{(\Omega)}}^2\mathrm{d}s\\
\leq&\max\left\{2\gamma_1,2\gamma_2\right\}\int_{0}^t
\left(\left\|u\right\|_{L^2{(\Omega)}}^2+\left\|v\right\|_{L^2{(\Omega)}}^2\right)\mathrm{d}s
+\left(\left\|u_0\right\|_{L^2{(\Omega)}}^2+\left\|v_0\right\|_{L^2{(\Omega)}}^2\right)
.
\end{aligned}
\end{equation}
Substituting (\ref{eq.2.1}) into (\ref{eq.2.6}) can ultimately lead to the desired
result (\ref{eq.2.2}). Thus, we complete
the proof.
\end{proof}

\subsection{Well-posedness of the weak solution}
Here, we give the well-posedness result of the weak solution as shown in the following theorem.
\begin{theorem}
Suppose that $w_0=\left(u,v\right)|_{t=0}=\left(u_0,v_0\right)\in L^2(\Omega)$.
Then there is only one pair of solution $w=(u,v)$, such that
\begin{equation*}
\begin{aligned}\displaystyle
w\in L^2\left(0,T;H^{\frac{\alpha}{2}}\left(\Omega\right)\right)
\cap L^4\left(0,T;L^4\left(\Omega\right)\right).
\end{aligned}
\end{equation*}
Moreover, it satisfies the energy estimate
\begin{equation}\label{eq.2.7}
\begin{aligned}\displaystyle
&\left\|u\right\|_{L^2{(\Omega)}}^2+\left\|v\right\|_{L^2{(\Omega)}}^2
+2\beta_1
\int_{0}^t
\left\|\left(-\triangle\right)^\frac{\alpha}{4}_\Omega u\right\|_{L^2{(\Omega)}}^2\mathrm{d}s
+
2\beta_2
\int_{0}^t
\left\|\left(-\triangle\right)^\frac{\alpha}{4}_\Omega v\right\|_{L^2{(\Omega)}}^2\mathrm{d}s\\
\leq&2\gamma_1\int_{0}^t \left\|u\right\|_{L^2{(\Omega)}}^2\mathrm{d}s
+
2\gamma_2\int_{0}^t \left\|v\right\|_{L^2{(\Omega)}}^2\mathrm{d}s+
\left\|u_0\right\|_{L^2{(\Omega)}}^2+\left\|v_0\right\|_{L^2{(\Omega)}}^2.
\end{aligned}
\end{equation}
\end{theorem}

\begin{proof}First, we consider the existence of weak solution.
With the help of the Theorem \ref{Th.2.5}, we can obtain the existence of the weak solution $w$ by using
a similar method in \cite{Pu} with Fourier-Galerkin approximation and compactness argument.

Next, we consider the uniqueness.
%
%
Suppose that there are two pair of solutions
$w=(u,v)$ and $\widetilde{w}=(\widetilde{u},\widetilde{v})$ satisfying the system (\ref{eq.1}) with (\ref{eq.2}).
Denote $u-\widetilde{u}=:\chi$ and $v-\widetilde{v}=:\widetilde{\chi}$, then system (\ref{eq.1}) becomes
\begin{equation}\label{eq.2.8}
\left\{
\begin{aligned}\displaystyle
&\partial_t \chi+\left(\beta_1+\mathrm{i}\eta_1\right)
\left(-\triangle\right)^\frac{\alpha}{2}_\Omega \chi
+\left(\mu_1+\mathrm{i}\zeta_1\right)\left(|u|^2u-|\widetilde{u}|^2\widetilde{u}\right)
-\gamma_1 \chi-\mathrm{i}\left(|u|^2v-|\widetilde{u}|^2v\right)
=0,\\
&\partial_t \widetilde{\chi}+\left(\beta_2+\mathrm{i}\eta_2\right)
\left(-\triangle\right)^\frac{\alpha}{2}_\Omega \widetilde{\chi}
+\left(\mu_2+\mathrm{i}\zeta_2\right)\left(|v|^2v-|\widetilde{v}|^2\widetilde{v}\right)
-\gamma_2 \widetilde{\chi}-\mathrm{i}\left(|v|^2u-|\widetilde{v}|^2u\right)
=0.
\end{aligned}
\right.
\end{equation}

Taking the inner product of the first equation of system (\ref{eq.2.8}) with $\chi$ and
considering the real part, we get
\begin{equation}\label{eq.2.9}
\begin{aligned}\displaystyle
&\frac{1}{2}\frac{\mathrm{d}}{\mathrm{d}t} \left\|\chi\right\|_{L^2{(\Omega)}}^2
+\beta_1
\left\|\left(-\triangle\right)^\frac{\alpha}{4}_\Omega \chi\right\|_{L^2{(\Omega)}}^2
+\Re\left\{\left(\mu_1+\mathrm{i}\zeta_1\right)\left(\left(|u|^2u-|\widetilde{u}|^2
\widetilde{u}\right),\chi\right)\right\}\\
&
-\gamma_1 \left\|\chi\right\|_{L^2{(\Omega)}}^2-
\Im\left\{\left(\left(|u|^2v-|\widetilde{u}|^2v\right),\chi\right)\right\}
=0.
\end{aligned}
\end{equation}
Using the same method to handle the second term of system (\ref{eq.2.8}), we also have
\begin{equation}\label{eq.2.10}
\begin{aligned}\displaystyle
&\frac{1}{2}\frac{\mathrm{d}}{\mathrm{d}t} \left\|\widetilde{\chi}\right\|_{L^2{(\Omega)}}^2
+\beta_2
\left\|\left(-\triangle\right)^\frac{\alpha}{4}_\Omega \widetilde{\chi}\right\|_{L^2{(\Omega)}}^2
+\Re\left\{\left(\mu_2+\mathrm{i}\zeta_2\right)
\left(\left(|v|^2v-|\widetilde{v}|^2\widetilde{v}\right),\widetilde{\chi}\right)\right\}\\
&
-\gamma_2 \left\|\widetilde{\chi}\right\|_{L^2{(\Omega)}}^2-
\Im\left\{\left(\left((|v|^2u-|\widetilde{v}|^2u\right),\widetilde{\chi}\right)\right\}
=0.
\end{aligned}
\end{equation}
Adding (\ref{eq.2.9}) and (\ref{eq.2.10}) can result in
\begin{equation}\label{eq.2.11}
\begin{aligned}\displaystyle
\frac{1}{2}\frac{\mathrm{d}}{\mathrm{d}t}\left( \left\|\chi\right\|_{L^2{(\Omega)}}^2
+\left\|\widetilde{\chi}\right\|_{L^2{(\Omega)}}^2\right)\leq
\gamma_1 \left\|\chi\right\|_{L^2{(\Omega)}}^2
+\gamma_2 \left\|\widetilde{\chi}\right\|_{L^2{(\Omega)}}^2
+\Xi,
\end{aligned}
\end{equation}
where
\begin{equation*}
\begin{aligned}\displaystyle
\Xi=&\Im\left\{\left(\left(|u|^2v-|\widetilde{u}|^2v\right),\chi\right)\right\}
+\Im\left\{\left(\left((|v|^2u-|\widetilde{v}|^2u\right),\widetilde{\chi}\right)\right\}\\
&
-\Re\left\{\left(\mu_1+\mathrm{i}\zeta_1\right)\left(\left(|u|^2u-|\widetilde{u}|^2
\widetilde{u}\right),\chi\right)\right\}-\Re\left\{\left(\mu_2+\mathrm{i}\zeta_2\right)\left(\left(|v|^2v-|\widetilde{v}|^2
\widetilde{v}\right),\widetilde{\chi}\right)\right\}.
\end{aligned}
\end{equation*}
It follows from Lemma \ref{Le.2.4}, we know that
\begin{equation}\label{eq.2.12}
\begin{aligned}\displaystyle
\Xi\leq&\left|\left(\left(|u|^2v-|\widetilde{u}|^2v\right),\chi\right)\right|
+\left|\left(\left((|v|^2u-|\widetilde{v}|^2u\right),\widetilde{\chi}\right)\right|\\
&
+\left|\mu_1+\mathrm{i}\zeta_1\right|\left|\left(\left(|u|^2u-|\widetilde{u}|^2
\widetilde{u}\right),\chi\right)\right|+
\left|\mu_2+\mathrm{i}\zeta_2\right|\left|\left(\left(|v|^2v-|\widetilde{v}|^2
\widetilde{v}\right),\widetilde{\chi}\right)\right|\\
\leq&\max\left\{2|u||v|,2|\widetilde{u}||v|,|\widetilde{u}|^2\right\}\left\|
\chi\right\|_{L^2{(\Omega)}}^2+
\max\left\{2|u||v|,2|{u}||\widetilde{v}|,|\widetilde{v}|^2\right\}
\left\|\widetilde{\chi}\right\|_{L^2{(\Omega)}}^2\\&+
\sqrt{\mu_1^2+\zeta_1^2}\max\left\{2|u|^2,2|\widetilde{u}||u|,|\widetilde{u}|^2\right\}
\left\|\chi\right\|_{L^2{(\Omega)}}^2\\&
+
\sqrt{\mu_2^2+\zeta_2^2}\max\left\{2|v|^2,2|\widetilde{v}||v|,|\widetilde{v}|^2\right\}
\left\|\widetilde{\chi}\right\|_{L^2{(\Omega)}}^2\\
\leq &\widetilde{c}\left(\left\|\chi\right\|_{L^2{(\Omega)}}^2+\left\|\widetilde{\chi}\right\|_{L^2{(\Omega)}}^2\right).
\end{aligned}
\end{equation}
Combining (\ref{eq.2.11}) and (\ref{eq.2.12}) leads to
\begin{equation*}
\begin{aligned}\displaystyle
\frac{\mathrm{d}}{\mathrm{d}t}\left( \left\|\chi\right\|_{L^2{(\Omega)}}^2
+\left\|\widetilde{\chi}\right\|_{L^2{(\Omega)}}^2\right)\leq
c\left( \left\|\chi\right\|_{L^2{(\Omega)}}^2
+\left\|\widetilde{\chi}\right\|_{L^2{(\Omega)}}^2\right),
\end{aligned}
\end{equation*}
where $c=2\max\left\{|\gamma_1|+\widetilde{c},|\gamma_2|+\widetilde{c}\right\}$.

By virtue of  the Lemma \ref{Le.2.2} and note that the initial conditions $\chi(0)=0$ and
 $\widetilde{\chi}(0)=0$, we can conclude that
\begin{equation*}
\begin{aligned}\displaystyle
\left\|\chi\right\|_{L^2{(\Omega)}}^2
+\left\|\widetilde{\chi}\right\|_{L^2{(\Omega)}}^2=0.
\end{aligned}
\end{equation*}
This implies that the weak solution of system (\ref{eq.1}) with (\ref{eq.2}) is unique.

Finally, integrating (\ref{eq.2.5}) with respect to $t$ and using
 the initial condition (\ref{eq.2}) can yield the inequality (\ref{eq.2.7}).
This ends the proof.
\end{proof}

\section{Numerical method}

Here, different from the previous theoretical analysis, we choose $\Omega=[a,b]$
 in order to make our method universal.
In practical computation, it is always necessary to restrict the original problem
on a bounded interval with homogeneous Dirichlet boundary conditions due to
the solutions of the system are fast decaying. i.e.,
\begin{equation}\label{eq.3.15}
\begin{aligned}
u(x,t)=0,\;v(x,t)=0,\;x\in{\mathds{R}}\backslash\Omega,\;t\in(0,T].
\end{aligned}
\end{equation}

In this section, we focus on the construction of high order numerical method to
resolve system  (\ref{eq.1}) with (\ref{eq.2}) and (\ref{eq.3.15})

\subsection{Notations and fractional Sobolev norm}

Let $x_j=a+jh, 0 \leq j\leq N$, where $N$ is a
positive integer and $h= (b-a)/N$ is the space step size,
and $t_k=k\tau, 0 \leq k\leq M$, where $M$ is a
positive integer and $\tau= T/M$ is the time step size.
Then
  the spatial domain $[a,b]$ and the time domain $[0,T]$
  are covered by
 $\Omega_{h}=\left\{x_j\,|\;x_j=a+jh,\;0\leq j\leq M\right\}$
 and  $\Omega_{\tau}=\left\{t_k\,|\;t_k=k\tau,\;0\leq k\leq N\right\}$, respectively. Denote
 by $\left\{u(x_j,t_k),v(x_j,t_k)\right\}$ and $\left\{U_j^k,V_j^k\right\}$ the exact
  solutions and numerical
 solutions of the problem (\ref{eq.1})--(\ref{eq.2}) with (\ref{eq.3.15}), respectively.

Denote
 $h\mathds{Z}=\left\{x_j=jh,\;j\in\mathds{Z}\right\}$ the infinite grid.
For any grid functions $u=\left\{u_j\right\}$,
$v=\left\{v_j\right\}$ on $h\mathds{Z}$,
we define
the discrete inner product and the associated $l_h^2$ norm as
$$\begin{array}{lll}
\displaystyle
\left(u,v\right)_{h}=h\sum_{j\in\mathds{Z}}u_j{v^{*}_j},\;\;
\left\|u\right\|_{h}^2=\left(u,u\right)_{h}.
\end{array}
$$
The discrete $l_h^p$ norm
and $l_h^\infty$ norm are defined by
  $$\begin{array}{lll}
\displaystyle
\left\|u\right\|_{l_h^p}^p=h\sum_{j\in\mathds{Z}}\left|u_j\right|^p,\;
\left\|u\right\|_{l_h^\infty}=\mathrm{sup}_{j\in\mathds{Z}}\left|u_j\right|,\;\;1\leq p<+\infty.
\end{array}
$$

Set $L_h^2=\left\{u\,|\,u=\left\{u_j\right\},\;\left\|u\right\|
_{h}<+\infty,\;j\in\mathds{Z} \right\}$. Then for a given constant $\sigma\in[0,1]$ and
function $u\in L_h^2$,
the fractional Sobolev
semi-norm $|\cdot|_{H_h^\sigma}$ and
norm $||\cdot||_{H_h^\sigma}$
can be defined by
$$\begin{array}{lll}
\displaystyle
|u|_{H_h^\sigma}^2=h\int_{-\pi}^{+\pi}h^{-2\sigma}|s|^{2\sigma}|{\hat{u}}(s)|^2\,\mathrm{d}s,\;\;
\left\|u\right\|_{H_h^\sigma}^2=h\int_{-\pi}^{+\pi}\left(1+h^{-2\sigma}|s|^{2\sigma}\right)|{\hat{u}}(s)|^2\,\mathrm{d}s,
\end{array}
$$
where the discrete
Fourier transform ${\hat{u}}(s)\in L^2\left[-\pi,\pi\right]$ is defined by
$$
{\hat{u}}(s)=\frac{1}{\sqrt{2\pi}}\sum_{j\in\mathds{Z}}
u_je^{-\mathrm{i}js}.
$$
Based on Parseval's identity, we immediately obtain
$\left\|u\right\|_{H_h^\sigma}^2=\left\|u\right\|_{h}^2+\left|u\right|_{H_h^\sigma}^2$.

%

Below, we will provide two lemmas proved in \cite{Ding2} related to the fractional Sobolev norm,
which play a crucial role in subsequent theoretical analysis.

\begin{lemma}\label{Lem.3.1}(Interpolation estimate inequality).
For any $0\leq\sigma_0\leq\sigma\leq1$, there holds that
$$\begin{array}{lll}
\displaystyle
\left\|u\right\|_{H_h^{\sigma_0}}\leq 2^{\frac{\sigma-\sigma_0}{2\sigma}}
\left\|u\right\|_{H_h^{\sigma}}^{\frac{\sigma_0}{\sigma}}\;
\left\|u\right\|_{{h}}^{1-\frac{\sigma_0}{\sigma}}.
\end{array}
$$
\end{lemma}

\begin{lemma}\label{Lem.3.2}
(Discrete Gagliardo-Nirenberg inequality)
For any $\frac{p-2}{2p}<\sigma_0\leq \sigma\leq1$, there holds that
$$\begin{array}{lll}
\displaystyle
\left\|u\right\|_{l_h^{p}}\leq 2^{\frac{\sigma-\sigma_0}{2\sigma}}
{\sigma_0}^{{\frac{2-p}{2p}}}
p^{-\frac{p}{2}}\left(\frac{p}{p-1}\right)^{\frac{p}{2(p-1)}}
\mathrm{B}^{\frac{p-2}{2p}}\left(\frac{1}{2\sigma_0},
\frac{p}{p-2}-\frac{1}{2\sigma_0}\right)
\left\|u\right\|_{H_h^{\sigma}}^{\frac{\sigma_0}{\sigma}}\;
\left\|u\right\|_{{h}}^{1-\frac{\sigma_0}{\sigma}},
\end{array}
$$
where $\mathrm{B}(\cdot,\cdot)$ is the Beta function defined by
$\mathrm{B}\left({\alpha_1},{\alpha_2}\right)=\int_{0}^{1}x^{{\alpha_1}-1}
(1-x)^{{\alpha_2}-1}\mathrm{d}x,\;\alpha_1,\alpha_2>0$.
\end{lemma}

\subsection{Constructing numerical differential formula for the fractional Laplacian}

The fundamental difficulty in designing an efficient numerical method to solve
(\ref{eq.1})--(\ref{eq.2}) with (\ref{eq.3.15}) lies in approximating the fractional Laplacian (\ref{eq.4}) efficiently.
It is worth mentioning that in \cite{Yang},
Yang et al. showed that the fractional Laplacian
is equivalent to the Riesz derivative defined on the infinite domain $-\infty < x < +\infty$, that is
\begin{equation}\label{eq.3.1}
-\left(-\Delta\right)^{\frac{\alpha}{2}}u\left(x\right)
=\partial_{x}^\alpha u(x),
\end{equation}
where the fractional Laplacian and the Riesz derivative
on the infinite domain $-\infty < x < +\infty$ are defined by
\begin{equation}\label{eq.3.2}
\left(-\triangle\right)^\frac{\alpha}{2}u(x)=K_\alpha\int_\mathds{R}
\frac{u(x)-u(y)}{|x-y|^{1+\alpha}}\mathrm{d}y,\;1<\alpha\leq2,
\end{equation}
and
\begin{equation}\label{eq.3.3}
\partial_{x}^\alpha u(x)=
-\frac{1}{2\cos\left(\frac{\pi\alpha}{2}\right)}\left
(\,_{RL}\mathrm{D}_{-\infty,x}^\alpha u(x)+\,_{RL}\mathrm{D}_{x,+\infty}^\alpha u(x)\right),
\;1<\alpha\leq2,
\end{equation}
respectively, in which $\,_{RL}\mathrm{D}_{-\infty,x}^\alpha u(x)$
and $\,_{RL}\mathrm{D}_{x,+\infty}^\alpha u(x)$ are the left-and
right-side Riemann-Liouville derivatives, that are
\begin{equation}\label{eq.3.4}
\,_{RL}{\mathrm{D}}_{-\infty,x}^{\alpha}u(x)=
\frac{1}{\mathrm{\Gamma}(2-\alpha)}\frac{\mathrm{d}^2}{\mathrm{d}
x^2}\int_{-\infty}^{x}\frac{u(y)}{(x-y)^{\alpha-1}}\textmd{d}y,\,
\end{equation}
and
\begin{equation}\label{eq.3.5}
\,_{RL}{\mathrm{D}}_{x,+\infty}^{\alpha}u(x)=
\frac{1}{\mathrm{\Gamma}(2-\alpha)}\frac{\mathrm{d}^2}{\mathrm{d}
x^2}\int_{x}^{+\infty}\frac{u(y)}{(y-x)^{\alpha-1}}\textmd{d}y.
\end{equation}

Inspired by the equivalent relation (\ref{eq.3.1}), we
 indirectly construct a higher-order numerical
 differential formula to approximate the fractional Laplacian (\ref{eq.3.2})
through the establishment
 of numerical differential formulas for left and right
 Riemann-Liouville derivatives (\ref{eq.3.4}) and (\ref{eq.3.5}).

 Let's first review the following generating function
$$
\begin{array}{lll}
\displaystyle
G_2(z)=\left(b_0+b_1z+b_2z^2\right)^{\alpha},\;b_0=\frac{3\alpha-2}{2\alpha},
\;b_1=-\frac{2(\alpha-1)}{\alpha},
\;b_2=\frac{\alpha-2}{2\alpha},
\end{array}
$$
which was established in \cite{Ding&Li}, and we find that it has an important property as follows:

\begin{lemma}\label{Le.3.3}
The generating function $G_2(z)$ is analytical, that is, it holds that
\begin{equation}\label{eq.3.6}
\displaystyle{G}_{2}(z)=
\sum\limits_{m=0}^{\infty}\kappa_{2,m}^{(\alpha)}z^m
\end{equation}
for all $\alpha\geq1$ and $|z|\leq1$.
\end{lemma}

\begin{proof}
Please refer to the appendix for detailed proof.
\end{proof}

In fact, from \cite{Ding&Li}, we find that the numerical differential formula constructed by generating
function ${G}_{2}(z)$ has only second-order accuracy. In order to get a higher order numerical
differential formula, we must find a way to construct a new generation function.
Inspired by the generating function ${G}_{2}(z)$, we establish
a novel generating function ${G}_{4}(z)$, which has the following form
\begin{equation*}
\displaystyle
 {G}_4(z)=\left[1+\eta\left(b_0+b_1z+b_2z^2\right)^{2}\right] {G}_2(z),
\end{equation*}
where $\eta=\frac{\left(3\alpha-2\right)\left(\alpha-2\right)\left(\alpha-1\right)}{24\alpha}$.

Similarly, we can also know that the generating function ${G}_4(z)$ has the same
 properties as ${G}_2(z)$, which can be described by the following theorem.

\begin{theorem}
The generating function $G_4(z)$ is analytical, that is, it holds that
$$
\begin{array}{lll}
\displaystyle{G}_{4}(z)=
\sum\limits_{m=0}^{\infty}\kappa_{4,m}^{(\alpha)}z^m
\end{array}
$$
for all $\alpha\geq1$ and $|z|\leq1$, where $\kappa_{4,m}^{(\alpha)}
=\kappa_{2,m}^{(\alpha)}+
\eta\kappa_{2,m}^{(\alpha+2)}$,
 which can also be calculated using the following recursive formulae
\begin{equation*}
\left\{
\begin{aligned}\kappa_{4,0}^{(\alpha)}=&b_0^{\alpha}
 \left(1+\eta b_0^2\right)
,\\
\kappa_{4,1}^{(\alpha)}=&
\frac{b_1}{b_0}\left[\alpha\kappa_{4,0}^{(\alpha)}+2\eta{\kappa}_{2,0}^{(\alpha+2)}\right],\\
\kappa_{4,m}^{(\alpha)}=&
\frac{1}{mb_0}\left[b_1\left(\alpha-m+1\right)\kappa_{4,m-1}^{(\alpha)}
+b_2\left(2\alpha-m+2\right)\kappa_{4,m-2}^{(\alpha)}
+2\eta b_1 \kappa_{2,m-1}^{(\alpha+2)}
 +
 4\eta b_2{\kappa}_{2,m-2}^{(\alpha+2)}
\right],m\geq2.
\end{aligned}
\right.
\end{equation*}
\end{theorem}
\begin{proof}
The result can be obtained by using the same proof method
as Lemma \ref{Le.3.3}, so we omit the specific details here.
\end{proof}

\subsubsection{Fourth-order numerical differential formulas of Riemann-Liouville derivatives}
Next, we start with the constructing of numerical differential formulas for the left and right
Riemann-Liouville derivatives (\ref{eq.3.4}) and (\ref{eq.3.5}). To this end,
denote
$$
\begin{array}{lll}
\displaystyle
\mathscr{C}^{{n+\alpha}}\left(\mathds{{R}}\right)=\left\{u\,|\,u\in L^1(\mathds{{R}}),
\;{\rm and}\;\int_{\mathds{{R}}}\left(1+\left|\xi\right|\right)^{{n+\alpha}}\left|
\hat{u}(\xi)
\right|\mathrm{d}\xi<\infty\right\},
\end{array}
$$
where $\hat{u}(\xi)$ is the Fourier transform of $u(x)$ defined as
 $\hat{u}(\xi)=\int_{\mathds{{R}}}u(x) e^{-\mathrm{i}\xi x}\mathrm{d}x$.

Meanwhile, the following fractional difference operators
\begin{equation}\label{eq.3.7}
\displaystyle \,^{L}\mathcal{B}_{h}^{\alpha}u(x)
=\frac{1}{h^{\alpha}}\sum\limits_{m=0}^{\infty}
\kappa_{4,m}^{(\alpha)}u\left(x-(m-1)h\right),
\end{equation}
and
\begin{equation}\label{eq.3.8}
 \,^{R}\mathcal{B}_{h}^{\alpha}u(x)
=\frac{1}{h^{\alpha}}\sum\limits_{m=0}^{\infty}
\kappa_{4,m}^{(\alpha)}u\left(x+(m-1)h\right)
\end{equation}
 are introduced.

In \cite{Tadjeran&Meerschaert}, Tadjerran et al. gave an asymptotic expansion of the
shifted Gr\"{u}nwald difference operator,
which plays a vital role in the construction of high-order numerical differential formulas.
Here, we provide a similar asymptotic expansion for operators (\ref{eq.3.7}) and (\ref{eq.3.8}).

\begin{lemma}\label{Th.3}
Let $u(x)\in \mathscr{C}^{{n+\alpha}}\left(\mathds{{R}}\right)$.
Then we can obtain
\begin{equation}\label{eq.3.9}
\begin{aligned}
\,^{L}\mathcal{B}_{h}^{\alpha}u(x)
 =\,_{RL}\mathrm{D}_{-\infty,x}^{\alpha}u(x)+
 \sum\limits_{m=1}^{n-1}\left(\varrho_{_m}^{(\alpha)}
 \,_{RL}\mathrm{D}_{-\infty,x}^{\alpha+m}u(x)\right)h^{m}
+\mathcal{O}\left(h^n\right)
\end{aligned}
\end{equation}
uniformly for $x\in \mathds{R}$, where
 $\varrho_{_m}^{(\alpha)}\;(m=1,2,\ldots)$ are the
 coefficients of the power series expansion of the function
  $\frac{{e}^z}{z^\alpha}{G}_4{({e}^{-z})}$, i.e.,
\begin{equation*}
\begin{aligned}
\frac{{e}^z}{z^\alpha}{G}_4{\left({e}^{-z}\right)}
=1+\sum\limits_{m=1}^{\infty}\varrho_m^{(\alpha)}z^m.
\end{aligned}
\end{equation*}
In particular, the first four coefficients are
\begin{equation*}
\begin{aligned}
&
\varrho_{_1}^{(\alpha)}=0,\;\;
\varrho_{_2}^{(\alpha)}= \frac{3\alpha^3-19\alpha^2+36\alpha-16}{24\alpha},\;\;
\varrho_{_3}^{(\alpha)}=0,
\\ &\varrho_{_4}^{(\alpha)}=- \frac{30\alpha^6-180\alpha^5+459\alpha^4
-835\alpha^3+1210\alpha^2-990\alpha+300}{720\alpha^3}.
\end{aligned}
\end{equation*}
\end{lemma}

\begin{proof}
Since the basic idea of proof is similar to that in \cite{Tadjeran&Meerschaert}, with only slight differences in
 the process, we omit the detailed proof here.
\end{proof}

Now, we derive the following high-order numerical differential
formulas for the Riemann-Liouville derivatives.
Our main results are stated as follows.

\begin{theorem}\label{Th.4} Suppose that
$u(x)\in \mathscr{C}^{{4+\alpha}}\left(\mathds{{R}}\right)$.
Then the following formula
\begin{equation}\label{eq.3.10}
\begin{aligned}
\left[1+\varrho_{_2}^{(\alpha)}h^2
\delta_x^2\right]
\,_{RL}\mathrm{D}_{-\infty,x}^{\alpha}u(x)
=\,^{L}\mathcal{B}_{h}^{\alpha}u(x)
+\mathcal{O}\left(h^4\right)
\end{aligned}
\end{equation}
uniformly holds for $x\in\mathds{{R}}$ as $h\rightarrow0$, where $\delta_x^2$
denote the second-order central difference operator, that is
 $\delta_x^2 u(x)=[u(x-h) -2u(x) +u(x +h)]/h^2$.
\end{theorem}

\begin{proof}
Based on the fact that the operator
$ \,_{RL}\mathrm{D}_{-\infty,x}^{\alpha+2}$
is the composition of operators $ \,_{RL}\mathrm{D}_{-\infty,x}^{\alpha}$
and $\frac{\mathrm{d}^2}{\mathrm{d}x^2}$, i.e.,
$ \,_{RL}\mathrm{D}_{-\infty,x}^{\alpha+2}
=\frac{\mathrm{d}^2}{\mathrm{d}x^2}\left[\,_{RL}\mathrm{D}_{-\infty,x}^{\alpha}\right]$,
then combining it with equation (\ref{eq.3.9}), we can get
\begin{equation}\label{eq.3.11}
\begin{aligned}
\,^{L}\mathcal{B}_{h}^{\alpha}u(x)
 =\left[1+\varrho_{_2}^{(\alpha)}h^2
\frac{\mathrm{d}^2}{\mathrm{d}x^2}\right]
\,_{RL}\mathrm{D}_{-\infty,x}^{\alpha}u(x)
+\mathcal{O}\left(h^4\right).
\end{aligned}
\end{equation}

On the other hand, noting that the approximation formula
$
\frac{\mathrm{d}^2u(x)}{\mathrm{d}x^2}={\delta_x^2u(x)}
+\mathcal{O}\left(h^2\right)
$
 and applying it to equation (\ref{eq.3.11}) will eventually lead to the result (\ref{eq.3.10}).
This completes the proof.
\end{proof}

For the right-side Riemann-Liouville derivative (\ref{eq.3.5}), similarly repeating the above process,
 we can easily obtain the following fourth-order approximation formula
\begin{equation}\label{eq.3.12}
\begin{aligned}
\left[1+\varrho_{_2}^{(\alpha)}h^2
\delta_x^2\right]
\,_{RL}\mathrm{D}_{x,+\infty}^{\alpha}u(x)
=\,^{R}\mathcal{B}_{h}^{\alpha}u(x)
+\mathcal{O}\left(h^4\right).
\end{aligned}
\end{equation}

\subsubsection{Fourth-order numerical differential formula of fractional Laplacian}
Based on the equivalency (\ref{eq.3.2}), we now turn to
construct the higher-order numerical differential
 formula of Riesz derivative. Noting the definition of (\ref{eq.3.3}) and formulas
 (\ref{eq.3.10}) and (\ref{eq.3.12}), we immediately get a fourth-order fractional compact
  numerical differential formula that approximates Riesz derivative as follows
\begin{equation*}
\displaystyle
\left[1+\varrho_{_2}^{(\alpha)}h^2
\delta_x^2\right]\partial_{x}^\alpha u(x)=
-\frac{1}{2\cos\left(\frac{\pi\alpha}{2}\right)}\left(\,^{L}\mathcal{B}_{h}^{\alpha}+
\,^{R}\mathcal{B}_{h}^{\alpha}\right)u(x)
+\mathcal{O}\left(h^4\right).
\end{equation*}

Combined with equation (\ref{eq.3.1}), a fourth-order fractional compact numerical differential formula
that approximates the fractional Laplacian can be obtained naturally and is shown as follows
\begin{equation}\label{eq.3.13}
\displaystyle
\left[1+\varrho_{_2}^{(\alpha)}h^2
\delta_x^2\right]\left(-\Delta\right)^{\frac{\alpha}{2}}u\left(x\right)=
\frac{1}{2\cos\left(\frac{\pi\alpha}{2}\right)}\left(\,^{L}\mathcal{B}_{h}^{\alpha}+
\,^{R}\mathcal{B}_{h}^{\alpha}\right)u(x)
+\mathcal{O}\left(h^4\right).
\end{equation}

For convenience, let $\mathcal{H}_h^{\alpha}=1+\varrho_{_2}^{(\alpha)}h^2
\delta_x^2$ and $\mathcal{B}_h^{\alpha}=-\frac{1}{2\cos\left(\frac{\pi\alpha}{2}\right)}
\left(\,^{L}\mathcal{B}_{h}^{\alpha}+
\,^{R}\mathcal{B}_{h}^{\alpha}\right)$, then (\ref{eq.3.13}) can be further simplified as
\begin{equation*}
\displaystyle
\mathcal{H}_h^{\alpha}\left\{\left(-\Delta\right)^{\frac{\alpha}{2}}u\left(x\right)\right\}=-
\mathcal{B}_h^{\alpha}u(x)
+\mathcal{O}\left(h^4\right),\;-\infty<x<+\infty.
\end{equation*}

Finally, let's consider the approximation formula of the fractional Laplacian
on bounded domain $\Omega=[a,b]$.
For $u\in C[a, b]$,  and $u(a) =u(b) =0$, we can make the zero-extension of $u$, such that $u$ is defined on $\mathds{R}$.
Then for $u(x)\in \mathscr{C}^{{4+\alpha}}\left(\Omega\right)$,
we have
\begin{equation} \label{eq.3.14}
\displaystyle
\mathcal{H}_h^{\alpha}\left\{\left(-\Delta\right)^{\frac{\alpha}{2}}_\Omega u\left(x\right)\right\}=-
\mathcal{A}_h^{\alpha}u(x)
+\mathcal{O}\left(h^4\right),\;a<x<b,
\end{equation}
where
$
\mathcal{A}_h^{\alpha}=-\frac{1}{2\cos\left(\frac{\pi\alpha}{2}\right)}
\left(\,^{L}\mathcal{A}_{h}^{\alpha}+
\,^{R}\mathcal{A}_{h}^{\alpha}\right)$
with
\begin{equation*}
\displaystyle \,^{L}\mathcal{A}_{h}^{\alpha}u(x)
=\frac{1}{h^{\alpha}}\sum\limits_{m=0}^{\left[\frac{x-a}{h}\right]}
\kappa_{4,m}^{(\alpha)}u\left(x-(m-1)h\right),\;\;
 \,^{R}\mathcal{A}_{h}^{\alpha}u(x)
=\frac{1}{h^{\alpha}}\sum\limits_{m=0}^{\left[\frac{b-x}{h}\right]}
\kappa_{4,m}^{(\alpha)}u\left(x+(m-1)h\right).
\end{equation*}

\subsection{Construction of a high order numerical algorithm for system (\ref{eq.1})}

 For any given grid function $U_j^k$ on $\Omega_{\tau h}=\Omega_{\tau}\times\Omega_{h}$,
 we define the
 following difference
  operators
 $$\delta_tU_j^{k-1/2}=\frac{1}{\tau}\left(U_j^k-U_{j}^{k-1}\right)$$
 and
 $$\delta_xU_{j-1/2}^k=\frac{1}{h}\left(U_j^k-U_{j-1}^k\right),
\;\;\delta_x^2U_{j}^k=\frac{1}{h}\left(\delta_xU_{j+1/2}^k-\delta_xU_{j-1/2}^k\right).$$

For the needs of subsequent analysis,
denote the index set $\mathds{Z}_M=\left\{j\;|\;j=1,2,\ldots,M-1\right\}$ and
the grid function space ${\mathcal{U}}_h=\left\{u\;|\;u=(u_1,u_2,\ldots,u_{M-1})\right\}.$
Meanwhile, under the boundary constraint (\ref{eq.3.15}), the inner product and norms previously defined
on an unbounded interval can be carried over to a finite interval $\Omega$. So we just restrict the
index $j $ from $1 $ to $M-1 $ in equations and continue to use them for convenience.
Under these considerations, the inequalities introduced in the previous section
still holds within a finite interval.

Now, considering the equation (\ref{eq.1}) at the grid point $(x_j,t)$, one has
\begin{equation*}
\left\{
\begin{aligned}\displaystyle
\partial_t u(x_j,t)+\left(\beta_1+\mathrm{i}\eta_1\right)
\left(-\triangle\right)^\frac{\alpha}{2}_\Omega u(x_j,t)
+\left(\mu_1+\mathrm{i}\zeta_1\right)|u(x_j,t)|^2u(x_j,t)\\
-\gamma_1 u(x_j,t)-\mathrm{i}|u(x_j,t)|^2v(x_j,t)
=0,\;j\in\mathds{Z}_M,\;0\leq t\leq T,\\
\partial_t v(x_j,t)+\left(\beta_2+\mathrm{i}\eta_2\right)
\left(-\triangle\right)^\frac{\alpha}{2}_\Omega v(x_j,t)
+\left(\mu_2+\mathrm{i}\zeta_2\right)|v(x_j,t)|^2v(x_j,t)\\
-\gamma_2 v(x_j,t)-\mathrm{i}|v(x_j,t)|^2u(x_j,t)
=0,\;j\in\mathds{Z}_M,\;0\leq t\leq T.
\end{aligned}
\right.
\end{equation*}

Performing the compact factor $\mathcal{H}_h^{\alpha}$ on both sides of the above system
and using the formula (\ref{eq.3.14}) can lead to
\begin{equation}\label{eq.3.16}
\left\{
\begin{aligned}\displaystyle
\mathcal{H}_h^{\alpha}\partial_t u(x_j,t)-\left(\beta_1+\mathrm{i}\eta_1\right)
\mathcal{A}_h^{\alpha}u(x_j,t)
+\left(\mu_1+\mathrm{i}\zeta_1\right)\mathcal{H}_h^{\alpha}|u(x_j,t)|^2u(x_j,t)\\
-\gamma_1 \mathcal{H}_h^{\alpha}u(x_j,t)-\mathrm{i}\mathcal{H}_h^{\alpha}|u(x_j,t)|^2v(x_j,t)
=\mathcal{O}\left(h^4\right),\;j\in\mathds{Z}_M,\;0\leq t\leq T,\\
\mathcal{H}_h^{\alpha}\partial_t v(x_j,t)-\left(\beta_2+\mathrm{i}\eta_2\right)
\mathcal{A}_h^{\alpha}v(x_j,t)
+\left(\mu_2+\mathrm{i}\zeta_2\right)\mathcal{H}_h^{\alpha}|v(x_j,t)|^2v(x_j,t)\\
-\gamma_2\mathcal{H}_h^{\alpha} v(x_j,t)-\mathrm{i}\mathcal{H}_h^{\alpha}|v(x_j,t)|^2u(x_j,t)
=\mathcal{O}\left(h^4\right),\;j\in\mathds{Z}_M,\;0\leq t\leq T.
\end{aligned}
\right.
\end{equation}

Omitting the high-order terms in (\ref{eq.3.16}) and denote $u_j(t)$ as
 the approximate solution of $u(x_j,t)$,
nd note that operator $\mathcal{H}_h^{\alpha}$ is
reversible (which will be proved in Lemma \ref{Lem.3.9} later),
we can get a semi-discrete scheme for system (\ref{eq.1}) as follows
\begin{equation}\label{eq.3.17}
\left\{
\begin{aligned}\displaystyle
\frac{\mathrm{d}u_j(t)}{\mathrm{d}t}=\left(\beta_1+\mathrm{i}\eta_1\right)\left(\mathcal{H}_h^{\alpha}\right)^{-1}
\mathcal{A}_h^{\alpha}u_j(t)
-\left(\mu_1+\mathrm{i}\zeta_1\right)|u_j(t)|^2u_j(t)\\
+\gamma_1 u_j(t)+\mathrm{i}|u_j(t)|^2v_j(t)
,\;j\in\mathds{Z}_M,\\
\frac{\mathrm{d}v_j(t)}{\mathrm{d}t}=\left(\beta_2+\mathrm{i}\eta_2\right)\left(\mathcal{H}_h^{\alpha}\right)^{-1}
\mathcal{A}_h^{\alpha}v_j(t)
-\left(\mu_2+\mathrm{i}\zeta_2\right)|v_j(t)|^2v_j(t)\\
+\gamma_2 v_j(t)+\mathrm{i}|v_j(t)|^2u_j(t),\;j\in\mathds{Z}_M.
\end{aligned}
\right.
\end{equation}

Denote ${\mathbf{u}}(t)=\left[u_1(t),u_2(t),\ldots,u_{M-1}(t)\right]^{T}$
with initial value ${\mathbf{u}}_0=\left[u_1(t_0),u_2(t_0),\ldots,u_{M-1}(t_{0})\right]^{T}$,
and ${\mathbf{v}}(t)=\left[v_1(t),v_2(t),\ldots,v_{M-1}(t)\right]^{T}$
with initial value ${\mathbf{v}}_0=\left[v_1(t_0),v_2(t_0),\ldots,v_{M-1}(t_{0})\right]^{T}$. Then
 the semi-discrete scheme (\ref{eq.3.17}) can be written in the matrix form:
\begin{equation}\label{eq.3.18}
\left\{
\begin{aligned}\displaystyle
&\frac{\mathrm{d}\mathbf{u}(t)}{\mathrm{d}t}=\left(\beta_1+\mathrm{i}\eta_1\right)\mathbf{Q}\mathbf{u}(t)
-\left(\mu_1+\mathrm{i}\zeta_1\right)|\mathbf{u}(t)|^2\mathbf{u}(t)
+\gamma_1 \mathbf{u}(t)+\mathrm{i}|\mathbf{u}(t)|^2\mathbf{v}(t)
,\;\mathbf{u}(0)=\mathbf{u}_0,\\
&\frac{\mathrm{d}\mathbf{v}(t)}{\mathrm{d}t}=\left(\beta_2+\mathrm{i}\eta_2\right)\mathbf{Qv}(t)
-\left(\mu_2+\mathrm{i}\zeta_2\right)|\mathbf{v}(t)|^2\mathbf{v}(t)
+\gamma_2 \mathbf{v}(t)+\mathrm{i}|\mathbf{v}(t)|^2\mathbf{u}(t),\;\mathbf{v}(0)=\mathbf{v}_0,
\end{aligned}
\right.
\end{equation}
where
$\mathbf{Q}={\mathbf{A}}^{-1}{\mathbf{B}}$ with
$${\mathbf{A}}=
\left(
  \begin{array}{ccccc}
1-2\varrho_{_2}^{(\alpha)} & \varrho_{_2}^{(\alpha)}&0
    & 0& \vspace{0.2 cm}\\
 \varrho_{_2}^{(\alpha)} &1-2\varrho_{_2}^{(\alpha)}
    & \varrho_{_2}^{(\alpha)}& \vspace{0.2 cm}\\
    \vdots &\vdots & \ddots & \ddots & \ddots \vspace{0.2 cm}\\
0 & \ldots &\varrho_{_2}^{(\alpha)}&1-2\varrho_{_2}^{(\alpha)}&
 \varrho_{_2}^{(\alpha)} \vspace{0.2 cm}\\
 0& 0& \ldots &\varrho_{_2}^{(\alpha)}&
  1-2\varrho_{_2}^{(\alpha)}
  \end{array}
\right)\in {\mathds{R}}^{(M-1)\times(M-1)},
$$
$${\mathbf{B}}=-\frac{1}{2h^{\alpha}\cos\left(\frac{\pi}{2}\alpha\right)}
 \left({\mathbf{H}}+{\mathbf{H}}^T\right),\;\;{\mathbf{H}}=
\left(
  \begin{array}{ccccc}
    \kappa_{4,1}^{(\alpha)} & \kappa_{4,0}^{(\alpha)}&
    & & \vspace{0.2 cm}\\
  \kappa_{4,2}^{(\alpha)} & \kappa_{4,1}^{(\alpha)} & \kappa_{4,0}^{(\alpha)}
    & & \vspace{0.2 cm}\\
    \vdots &\vdots & \ddots & \ddots & \ddots \vspace{0.2 cm}\\
 \kappa_{4,M-2}^{(\alpha)} & \ldots & \kappa_{4,2}^{(\alpha)}&\kappa_{4,1}^{(\alpha)}&
 \kappa_{4,0}^{(\alpha)} \vspace{0.2 cm}\\
 \kappa_{4,M-1}^{(\alpha)}& \kappa_{4,M-2}^{(\alpha)}&
 \ldots & \kappa_{4,2}^{(\alpha)}&
  \kappa_{4,1}^{(\alpha)}
  \end{array}
\right)\in {\mathds{R}}^{(M-1)\times(M-1)}.
$$

By using the compact implicit integration factor method, the solution of the
system (\ref{eq.3.18}) at time $t_{k+1}$
can be represented as
\begin{equation}\label{eq.3.19}
\left\{
\begin{aligned}\displaystyle
{{\mathbf{u}}}\left(t_{k+1}\right)
=&e^{\tau\left[\left(\beta_1+\mathrm{i}\eta_1\right)\mathbf{Q}+\gamma_1\mathbf{I}\right]}{\mathbf{u}}\left(t_{k}\right)\\
&+e^{\tau\left[\left(\beta_1+\mathrm{i}\eta_1\right)\mathbf{Q}+\gamma_1\mathbf{I}\right]}\int_{0}^{\tau}
e^{-s\left[\left(\beta_1+\mathrm{i}\eta_1\right)\mathbf{Q}+\gamma_1\mathbf{I}\right]}
F_1\left({\mathbf{u}}\left(t_{k}+s\right),{\mathbf{v}}\left(t_{k}+s\right)\right)\mathrm{d}{s},\\
{{\mathbf{v}}}\left(t_{k+1}\right)
=&e^{\tau\left[\left(\beta_2+\mathrm{i}\eta_2\right)\mathbf{Q}+\gamma_2\mathbf{I}\right]}{\mathbf{v}}\left(t_{k}\right)\\
&+e^{\tau\left[\left(\beta_2+\mathrm{i}\eta_2\right)\mathbf{Q}+\gamma_2\mathbf{I}\right]}\int_{0}^{\tau}
e^{-s\left[\left(\beta_2+\mathrm{i}\eta_2\right)\mathbf{Q}+\gamma_2\mathbf{I}\right]}
F_2\left({\mathbf{u}}\left(t_{k}+s\right),{\mathbf{v}}\left(t_{k}+s\right)\right)\mathrm{d}{s}.
\end{aligned}
\right.
\end{equation}
where $\mathbf{I}$ denotes the $(M-1)\times(M-1)$ identity matrix, and the functions
 $$F_1\left({\mathbf{u}}\left(t_{k}+s\right),{\mathbf{v}}\left(t_{k}+s\right)\right)=
-\left(\mu_1+\mathrm{i}\zeta_1\right)|\mathbf{u}\left(t_{k}+s\right)|^2\mathbf{u}\left(t_{k}+s\right)
+\mathrm{i}|\mathbf{u}\left(t_{k}+s\right)|^2\mathbf{v}\left(t_{k}+s\right),
$$
 and
 $$F_2\left({\mathbf{u}}\left(t_{k}+s\right),{\mathbf{v}}\left(t_{k}+s\right)\right)=
-\left(\mu_2+\mathrm{i}\zeta_2\right)|\mathbf{v}\left(t_{k}+s\right)|^2\mathbf{v}\left(t_{k}+s\right)
+\mathrm{i}|\mathbf{v}\left(t_{k}+s\right)|^2\mathbf{u}\left(t_{k}+s\right).$$

For the integrands in (\ref{eq.3.19}), we use
the $(r-1)$-th order Lagrange interpolation polynomial with interpolation points
$t_{k+1},t_{k},\ldots,t_{k-r+2}$ to approximate them and
 obtain the following
$r$-th order scheme
\begin{equation}\label{eq.3.20}
\left\{
\begin{aligned}\displaystyle
{{\mathbf{u}}}\left(t_{k+1}\right)
=&e^{\tau\left[\left(\beta_1+\mathrm{i}\eta_1\right)\mathbf{Q}+\gamma_1\mathbf{I}\right]}{\mathbf{u}}\left(t_{k}\right)
+\tau\left[\rho_{_1}F_1\left({\mathbf{u}}\left(t_{k+1}\right),{\mathbf{v}}\left(t_{k+1}\right)\right)\right.\\&\left.+
\sum_{j=0}^{r-2}\rho_{_{-j}}e^{\tau(j+1)\left[\left(\beta_1+\mathrm{i}\eta_1\right)\mathbf{Q}+\gamma_1\mathbf{I}\right]}
F_1\left({\mathbf{u}}\left(t_{k-j}\right),{\mathbf{v}}\left(t_{k-j}\right)\right)\right]+\mathcal{O}\left(\tau^r\right),\;
r\geq2,
\\
{{\mathbf{v}}}\left(t_{k+1}\right)
=&e^{\tau\left[\left(\beta_2+\mathrm{i}\eta_2\right)\mathbf{Q}+\gamma_2\mathbf{I}\right]}{\mathbf{v}}\left(t_{k}\right)
+\tau\left[\rho_{_1}F_2\left({\mathbf{u}}\left(t_{k+1}\right),{\mathbf{v}}\left(t_{k+1}\right)\right)\right.\\&\left.+
\sum_{j=0}^{r-2}\rho_{_{-j}}e^{\tau(j+1)\left[\left(\beta_2+\mathrm{i}\eta_2\right)\mathbf{Q}+\gamma_2\mathbf{I}\right]}
F_2\left({\mathbf{u}}\left(t_{k-j}\right),{\mathbf{v}}\left(t_{k-j}\right)\right)\right]+\mathcal{O}\left(\tau^r\right),\;
r\geq2,
\end{aligned}
\right.
\end{equation}
where the values of coefficients $\rho_{_{1}},\rho_{_{0}},\rho_{_{-1}},\ldots,\rho_{_{2-r}}$ are
obtained by \cite{Jian}
\begin{equation*}
\begin{aligned}
\displaystyle
\rho_{_{-j}}=\frac{1}{\tau}\int_{0}^{\tau}\mathop{\Pi}\limits^{r-2}_{n=-1,n\neq j}
\frac{s+n\tau}{(n-j)\tau}\mathrm{d}s,\;\;-1\leq j\leq r-2.
\end{aligned}
\end{equation*}

Below, in particular, we choose $r=2$ in equations (\ref{eq.3.20}) to obtain a temporal second-order scheme, that is
\begin{equation}\label{eq.3.21}
\left\{
\begin{aligned}\displaystyle
{{\mathbf{u}}}\left(t_{k+1}\right)
=&e^{\tau\left[\left(\beta_1+\mathrm{i}\eta_1\right)\mathbf{Q}+\gamma_1\mathbf{I}\right]}{\mathbf{u}}\left(t_{k}\right)
+\frac{\tau}{2}\left[F_1\left({\mathbf{u}}\left(t_{k+1}\right),{\mathbf{v}}\left(t_{k+1}\right)\right)\right.\\&\left.+
e^{\tau \left[\left(\beta_1+\mathrm{i}\eta_1\right)\mathbf{Q}+\gamma_1\mathbf{I}\right]}
F_1\left({\mathbf{u}}\left(t_{k}\right),{\mathbf{v}}\left(t_{k}\right)\right)\right]+\mathcal{O}\left(\tau^2\right),
\\
{{\mathbf{v}}}\left(t_{k+1}\right)
=&e^{\tau\left[\left(\beta_2+\mathrm{i}\eta_2\right)\mathbf{Q}+\gamma_2\mathbf{I}\right]}{\mathbf{v}}\left(t_{k}\right)
+\frac{\tau}{2}\left[F_2\left({\mathbf{u}}\left(t_{k+1}\right),{\mathbf{v}}\left(t_{k+1}\right)\right)\right.\\&\left.+
e^{\tau\left[\left(\beta_2+\mathrm{i}\eta_2\right)\mathbf{Q}+\gamma_2\mathbf{I}\right]}
F_2\left({\mathbf{u}}\left(t_{k}\right),{\mathbf{v}}\left(t_{k}\right)\right)\right]+\mathcal{O}\left(\tau^2\right).
\end{aligned}
\right.
\end{equation}
Dealing with the exponential matrix function by means of the $[1,1]$ pad{\'{e}} approximation formula \cite{Moler},
the scheme (\ref{eq.3.21}) further degenerates into
\begin{equation}\label{eq.3.22}
\left\{
\begin{aligned}\displaystyle
&\delta_t{{\mathbf{u}}}\left(t_{k+1/2}\right)
=\left[\left(\beta_1+\mathrm{i}\eta_1\right)\mathbf{Q}+\gamma_1\mathbf{I}\right]{\mathbf{u}}\left(t_{k+1/2}\right)
+F_1\left({\mathbf{u}}\left(t_{k+1/2}\right),{\mathbf{v}}\left(t_{k+1/2}\right)\right)+\mathcal{O}\left(\tau^2\right),
\\
&\delta_t{{\mathbf{v}}}\left(t_{k+1/2}\right)
=\left[\left(\beta_2+\mathrm{i}\eta_2\right)\mathbf{Q}+\gamma_2\mathbf{I}\right]{\mathbf{v}}\left(t_{k+1/2}\right)+
F_2\left({\mathbf{u}}\left(t_{k+1/2}\right),{\mathbf{v}}\left(t_{k+1/2}\right)\right)+\mathcal{O}\left(\tau^2\right).
\end{aligned}
\right.
\end{equation}

Once again, the higher-order term in equation (\ref{eq.3.22}) is discarded, and finally the
 following numerical scheme for solving system (\ref{eq.1}) with (\ref{eq.2}) and (\ref{eq.3.15}) can be obtained
\begin{equation}\label{eq.3.23}
\left\{
\begin{aligned}\displaystyle
\mathcal{H}_h^{\alpha}\delta_tU_j^{k+1/2}
=&\left[\left(\beta_1+\mathrm{i}\eta_1\right)
\mathcal{A}_h^{\alpha}+\gamma_1\mathcal{H}_h^{\alpha}\right]U_j^{k+1/2}
-\left(\mu_1+\mathrm{i}\zeta_1\right)\mathcal{H}_h^{\alpha}\left|U_j^{k+1/2}\right|^2U_j^{k+1/2}\\
&+\mathrm{i}\mathcal{H}_h^{\alpha}\left|U_j^{k+1/2}\right|^2V_j^{k+1/2},\;j\in\mathds{Z}_M,\;0\leq k\leq N-1,
\\
\mathcal{H}_h^{\alpha}\delta_tV_j^{k+1/2}
=&\left[\left(\beta_2+\mathrm{i}\eta_2\right)
\mathcal{A}_h^{\alpha}+\gamma_2\mathcal{H}_h^{\alpha}\right]V_j^{k+1/2}
-\left(\mu_2+\mathrm{i}\zeta_2\right)\mathcal{H}_h^{\alpha}\left|V_j^{k+1/2}\right|^2V_j^{k+1/2}\\
&+\mathrm{i}\mathcal{H}_h^{\alpha}\left|V_j^{k+1/2}\right|^2U_j^{k+1/2},\;j\in\mathds{Z}_M,\;0\leq k\leq N-1,
\end{aligned}
\right.
\end{equation}
\vspace{-0.5cm}
\begin{equation}\label{eq.3.24}
\begin{aligned}
U_j^0=u_0(x_j),\;V_j^0=v_0(x_j),\;j\in\mathds{Z},
\end{aligned}
\end{equation}
\vspace{-0.5cm}
\begin{equation}\label{eq.3.25}
\begin{aligned}
U_j^k=0,\;V_j^k=0,\;j\in{\mathds{Z}}\backslash\mathds{Z}_M,\;0\leq k\leq N.
\end{aligned}
\end{equation}

\subsection{Theory analysis of  the numerical algorithm}
In this section, we study the boundedness, unique solvability and
convergence of the difference scheme (\ref{eq.3.23})
 with (\ref{eq.3.24}) and (\ref{eq.3.25}).
For this purpose, we first introduce several auxiliary lemmas.

\subsubsection{Auxiliary lemmas}

\begin{lemma}(See \cite{Ding&Li})\label{Lem.3.7}
Let
\begin{equation*}
\begin{aligned}
Z_1(\alpha,s)=
\cos\left[\left(\frac{s-\pi}{2}+\theta\right)\alpha-s\right],\;
\end{aligned}
\end{equation*}
with $\theta=\arctan\frac{d_2}{d_1},\;d_1=b_0-b_2\cos s,\;d_2=-b_2\sin s$, $b_0$
 and $b_2$ as defined above. Then there holds that
\begin{equation*}
\begin{aligned}
-1\leq Z_1(\alpha,s)\leq\cos\left(\frac{\pi}{2}\alpha\right),
\end{aligned}
\end{equation*}
for $1<\alpha\leq2$ and $0\leq s\leq\pi$.
\end{lemma}

\begin{lemma}\label{Lem.3.8}
Denote
\begin{equation*}
\begin{aligned}
Z(\alpha,s)=Z_1(\alpha,s)+4\eta\left(d_1^2+d_2^2\right)\sin^2\left(\frac{s}{2}\right)Z_1(\alpha+2,s),
\end{aligned}
\end{equation*}
then we have
\begin{equation}\label{eq.3.26}
\begin{aligned}
-1+\frac{16\eta(\alpha-1)^2}{\alpha^2}\leq
 Z(\alpha,s)\leq \cos\left(\frac{\pi}{2}\alpha\right)-4\eta,\
\end{aligned}
\end{equation}
for $1<\alpha\leq2$ and $0\leq s\leq\pi$, where
$\eta=\frac{(3\alpha-2)(\alpha-2)(\alpha-1)}{24\alpha}$. Furthermore, a more direct result is $Z(\alpha,s)\leq0$.
\end{lemma}

\begin{proof}
On one hand, it follows from the Lemma \ref{Lem.3.7} that
\begin{equation}\label{eq.3.27}
\begin{aligned}
Z(\alpha,s)\leq Z_1(\alpha,s)-4\eta\left(d_1^2+d_2^2\right)\leq
Z_1(\alpha,s)-4\eta\leq
\cos\left(\frac{\pi}{2}\alpha\right)-4\eta.
\end{aligned}
\end{equation}
On the contrary, we can also know that
\begin{equation}\label{eq.3.28}
\begin{aligned}
Z(\alpha,s)\geq Z_1(\alpha,s)+4\eta\left(d_1^2+d_2^2\right)\geq
-1+4\eta\left(b_0+b_2\right)^2.
\end{aligned}
\end{equation}
Combining (\ref{eq.3.27}) and (\ref{eq.3.28}), the result (\ref{eq.3.26}) can be obtained naturally.

Furthermore, from result (\ref{eq.3.27}), we can conclude that
\begin{equation*}
\begin{aligned}
Z(\alpha,s)\leq
\cos\left(\frac{\pi}{2}\alpha\right)+\frac{1}{3}(2-\alpha)(\alpha-1)=:p(\alpha).
\end{aligned}
\end{equation*}
Using Kober inequality \cite{Bhayo}, i.e., $\cos x\geq 1-\frac{2}{\pi}x$, $x\in[0,\pi/2]$, it can lead to
\begin{equation*}
\begin{aligned}
p'(\alpha)=
-\frac{\pi}{2}\cos\left(\frac{\pi}{2}\alpha-\frac{\pi}{2}\right)+\frac{1}{3}(3-2\alpha)
\leq-\pi+1+\frac{\alpha}{6}(3\pi-4)\leq-\frac{1}{3}.
\end{aligned}
\end{equation*}
This implies that $p(\alpha)$ is a decreasing function, so
$p_{\max}(\alpha)\leq p(1)=0$. This ends the proof.
\end{proof}

Denote $U=[U_1,U_2,\ldots,U_{M-1}]^T$, note that
the matrices $A$ and
 ${B}$
are the associate matrices of fractional difference quotient operators $\mathcal{H}_h^{\alpha}$
and $\mathcal{A}_h^{\alpha}$, respectively, i.e., $\mathcal{H}_h^{\alpha}U=AU$ and
$\mathcal{A}_h^{\alpha}U=BU$. Then we have the following results:

\begin{lemma}\label{Lem.3.9}The matrix ${A}^{-1}$ is a
real-valued symmetry
  positive definite for $1<\alpha\leq2$.
\end{lemma}
\begin{proof}
It is not difficult to find that the eigenvalues of tridiagonal matrix $A$ are
\begin{equation*}
\begin{aligned}
\lambda_j(A)=1-\frac{3\alpha^3-19\alpha^2+36\alpha-16}{6\alpha}\sin^2\left(\frac{j\pi}{2M}\right),
\;j=1,2,\ldots,M-1,
\end{aligned}
\end{equation*}
which implies
\begin{equation*}
\begin{aligned}
1\geq\lambda_j(A)\geq \frac{-3\alpha^3+19\alpha^2-30\alpha+16}{6\alpha}=:M_1(\alpha)>0,
\;1<\alpha\leq2.
\end{aligned}
\end{equation*}
This indicates that the matrix $A$ is positive definite. Hence, from the knowledge of
matrices, we know that associate matrix $A^{-1}$ of fractional difference quotient
operator $\left(\mathcal{H}_h^{\alpha}\right)^{-1}$
is also real, symmetric and positive definite.
\end{proof}

\begin{lemma}\label{Lem.3.10}The matrix ${B}$ is a
real-valued symmetry
 negative semi-definite for $1<\alpha\leq2$.
\end{lemma}

\begin{proof}
It is obvious that ${B}$ is a real-valued symmetry matrix.
Next, we only consider it to be negative semi-definite.
Let's consider its generating function
$
F(\alpha,s)=
-\frac{G(\alpha,s)}{2h^\alpha\cos\left(\frac{\pi}{2}\alpha\right)}
$
. Note that
\begin{equation*}
\begin{aligned}
G(\alpha,s)=&e^{-\mathrm{i}s}G_4\left(e^{\mathrm{i}s}\right)
+e^{\mathrm{i}s}G_4\left(e^{-\mathrm{i}s}\right)\\
=&b_0^{\alpha}
\left[e^{-\mathrm{i}s}
\left(1-e^{\mathrm{i}s}\right)^{\alpha}
\left(1+\frac{b_0+b_1}{b_0}e^{\mathrm{i}s}
+\frac{b_0+b_1+b_2}{b_0}e^{2\mathrm{i}s}\right)^{\alpha}\right.\\&
+\eta
e^{-\mathrm{i}s}
\left(1-e^{\mathrm{i}s}\right)^{\alpha+2}
\left(1+\frac{b_0+b_1}{b_0}e^{\mathrm{i}s}
+\frac{b_0+b_1+b_2}{b_0}e^{2\mathrm{i}s}\right)^{\alpha+2}\\
&+e^{\mathrm{i}s}
\left(1-e^{-\mathrm{i}s}\right)^{\alpha}
\left(1+\frac{b_0+b_1}{b_0}e^{-\mathrm{i}s}
+\frac{b_0+b_1+b_2}{b_0}e^{-2\mathrm{i}s}\right)^{\alpha}\\
&\left.+\eta
e^{\mathrm{i}s}
\left(1-e^{-\mathrm{i}s}\right)^{\alpha+2}
\left(1+\frac{b_0+b_1}{b_0}e^{-\mathrm{i}s}
+\frac{b_0+b_1+b_2}{b_0}e^{-2\mathrm{i}s}\right)^{\alpha+2}
\right].
\end{aligned}
\end{equation*}

Obviously, it is not difficult to see that $G(\alpha,s)$ is a
real-valued even function with respect to $s$. Thus, we can only consider its principle
value for $s\in[0,\pi]$. Hence, we further know that
\begin{equation*}
\begin{aligned}
G(\alpha,s)=&\left(2\sin\frac{s}{2}\right)^\alpha
\left(d_1^2+d_2^2\right)^{\frac{\alpha}{2}}Z_1(\alpha,s)
+\eta\left(2\sin\frac{s}{2}\right)^{\alpha+2}
\left(d_1^2+d_2^2\right)^{\frac{\alpha+2}{2}}Z_1(\alpha+2,s)\\
=&\left(2\sin\frac{s}{2}\right)^\alpha
\left(d_1^2+d_2^2\right)^{\frac{\alpha}{2}}
Z(\alpha,s)\leq0
\end{aligned}
\end{equation*}
by Lemma \ref{Lem.3.8}.
Combined with Grenander-Szeg\"{o} Theorem
\cite{Chan},
 we can know that
the matrix ${B}$ is
 negative semi-definite.
This ends the proof.
\end{proof}

\begin{lemma}\label{Lem.5}
For any $1<\alpha\leq2$, we have
\begin{equation}\label{eq.3.29}
\begin{aligned}
C_2(\alpha)\left|u\right|_{H_h^\frac{\alpha}{2}}^2
\leq\left(\left(\mathcal{H}_h^{\alpha}\right)^{-1}
\mathcal{A}_h^{\alpha}u,u\right)_{h}\leq C_1(\alpha)
\left|u\right|_{H_h^\frac{\alpha}{2}}^2,
\end{aligned}
\end{equation}
where
\begin{equation*}
\begin{aligned}
C_1(\alpha)=\frac{(b_0+b_2)^\alpha\left[4\eta-\cos\left(\frac{\pi}{2}\alpha\right)\right]
}{2\cos\left(\frac{\pi}{2}\alpha\right)}
\left(\frac{2}{\pi}\right)^{\alpha},\;\;
C_2(\alpha)=\frac{(b_0-b_2)^\alpha\left[\alpha^2-16\eta(\alpha-1)^2\right]
}{2\alpha^2\cos\left(\frac{\pi}{2}\alpha\right)\left[1-4\varrho_{_2}^{(\alpha)}\right]}
.
\end{aligned}
\end{equation*}
\end{lemma}

\begin{proof}
It follows from the former Parseval's identity that
\begin{equation}\label{eq.3.30}
\begin{aligned}
\left(\left(\mathcal{H}_h^{\alpha}\right)^{-1}
\mathcal{A}_h^{\alpha}u,u\right)_{h}&
=\displaystyle h\int_{-\pi}^{+\pi}f(\alpha,s)|{\hat{u}}(s)|^2\; \mathrm{d}s,
\end{aligned}
\end{equation}
where
\begin{equation*}
\begin{aligned}
f(\alpha,s)=F(\alpha,s)\left[1-4\varrho_{_2}^{(\alpha)}
\sin^2\left(\frac{s}{2}\right)\right]^{-1},
\end{aligned}
\end{equation*}
that is
\begin{equation*}
\begin{aligned}
f(\alpha,s)=-
\frac{Z(\alpha,s)}{2h^\alpha\cos\left(\frac{\pi}{2}\alpha\right)}
\left(2\sin\frac{s}{2}\right)^\alpha
\left(d_1^2+d_2^2\right)^{\frac{\alpha}{2}}
\left[1-4\varrho_{_2}^{(\alpha)}
\sin^2\left(\frac{s}{2}\right)\right]^{-1}.
\end{aligned}
\end{equation*}
Noting the fact that $\frac{2s}{\pi}\leq2\sin\frac{s}{2}\leq s$ for
$s\in[0,\pi]$, $(b_0+b_2)^2\leq d_1^2+d_2^2\leq(b_0-b_2)^2$,
$1-4\varrho_{_2}^{(\alpha)}\leq1-4\varrho_{_2}^{(\alpha)}
\sin^2\left(\frac{s}{2}\right)\leq1$ and using Lemma \ref{Lem.3.8}, we have
\begin{equation*}
\begin{aligned}
f(\alpha,s)\leq\frac{(b_0+b_2)^\alpha\left[4\eta-\cos\left(\frac{\pi}{2}\alpha\right)\right]
}{2\cos\left(\frac{\pi}{2}\alpha\right)}
\left(\frac{2}{\pi}\right)^{\alpha}|s|^{\alpha}
h^{-\alpha}=:C_1(\alpha)|s|^{\alpha}
h^{-\alpha},
\end{aligned}
\end{equation*}
and
\begin{equation*}
\begin{aligned}
f(\alpha,s)\geq\frac{(b_0-b_2)^\alpha\left[\alpha^2-16\eta(\alpha-1)^2\right]
}{2\alpha^2\cos\left(\frac{\pi}{2}\alpha\right)\left[1-4\varrho_{_2}^{(\alpha)}\right]}
|s|^{\alpha}
h^{-\alpha}=:C_2(\alpha)|s|^{\alpha}
h^{-\alpha}.
\end{aligned}
\end{equation*}
These together with (\ref{eq.3.30}) implies (\ref{eq.3.29}) holds and thus completes the proof.
\end{proof}

Next, based on the Lemmas \ref{Lem.3.9} and \ref{Lem.3.10}, we give two lemmas without
proof. For the detailed proof process,
please refer to \cite{Ding2}.

\begin{lemma}\label{Lem.3.12}For any grid function
${U},V\in U_h$, there exists an
asymmetric positive difference quotient operator denoted
by $\nabla_h^{\frac{\alpha}{2}}$, such that
$$
\begin{array}{lll}
\displaystyle
\left(\left(\mathcal{H}_h^{\alpha}\right)^{-1}U,{V}\right)_{h}
=\left(\nabla_h^{\frac{\alpha}{2}}U,\nabla_h^{\frac{\alpha}
{2}}V\right)_{h}.
\end{array}
$$
\end{lemma}


\begin{lemma}\label{Lem.3.13}For any grid function
${U},V\in U_h$, there exists an
asymmetric positive difference quotient operator denoted
by $\delta_h^{\frac{\alpha}{2}}$, such that
$$
\begin{array}{lll}
\displaystyle
-\left(\mathcal{A}_h^{\alpha}U,{V}\right)_{h}
=\left(\delta_h^{\frac{\alpha}{2}}U,\delta_h^{\frac{\alpha}
{2}}V\right)_{h}.
\end{array}
$$
\end{lemma}

Finally, we introduce a discrete norm as follows
\begin{equation*}
\begin{aligned}
\left\||U|\right\|=\sqrt{\left(\left(\mathcal{H}_h^{\alpha}\right)^{-1}U,U\right)_h}\;.
\end{aligned}
\end{equation*}
Norms $\left\||\cdot|\right\|$ and $\left\|\cdot\right\|_h$ have the following relationship:

\begin{lemma}
The discrete norms $\left\||\cdot|\right\|$ and $\left\|\cdot\right\|$ are equivalent, in the sense that,
\begin{equation*}
\begin{aligned}
\left\|U\right\|_h\leq\left\||U|\right\|\leq\frac{1}{\sqrt{M_1(\alpha)}}\left\|U\right\|_h.
\end{aligned}
\end{equation*}
\end{lemma}

\begin{proof}
It follows from the Lemma \ref{Lem.3.9} that the eigenvalues of $A^{-1}$ satisfy
\begin{equation*}
\begin{aligned}
1\leq\lambda_j\left(A^{-1}\right)\leq \frac{1}{M_1(\alpha)},\;j=1,2,\ldots,M-1,
\end{aligned}
\end{equation*}
which means that the spectral radius $\rho(A^{-1})\leq\frac{1}{M_1(\alpha)}$, and
consequently
\begin{equation*}
\begin{aligned}
\left\|A^{-1}\right\|_h=\rho(A^{-1})\leq\frac{1}{M_1(\alpha)}.
\end{aligned}
\end{equation*}
Thus, it can be concluded that
\begin{equation}\label{eq.3.31}
\begin{aligned}
\left\||U|\right\|^2=\left(\left(\mathcal{H}_h^{\alpha}\right)^{-1}U,U\right)_h
\leq\left\|A^{-1}U\right\|_h\left\|U\right\|_h\leq\frac{1}{M_1(\alpha)}\left\|U\right\|_h^2.
\end{aligned}
\end{equation}

On the other hand, we also have
\begin{equation}\label{eq.3.32}
\begin{aligned}
\left\|U\right\|^2_h=&\left(\mathcal{H}_h^{\alpha}\left(\mathcal{H}_h^{\alpha}\right)^{-1}U,U\right)_h
=\left(\left(\mathcal{H}_h^{\alpha}\right)^{-1}U,U\right)_h
+\varrho_{_2}^{(\alpha)}h^2\left(\delta_x^2\left(\mathcal{H}_h^{\alpha}\right)^{-1}U,U\right)_h\\
=&\left\||U|\right\|^2-\varrho_{_2}^{(\alpha)}h^2\left\|\delta_x\nabla_h^{\frac{\alpha}{2}}U\right\|^2_h
\leq\left\||U|\right\|^2.
\end{aligned}
\end{equation}
by follows from the Lemma \ref{Lem.3.12}. Combining (\ref{eq.3.31}) with (\ref{eq.3.32}) gives the desired result.
\end{proof}

\subsubsection{Boundedness of solution}
\begin{theorem}\label{Th.3.15}
The difference solution of scheme (\ref{eq.3.23}) -- (\ref{eq.3.25}) is
bounded. That is, there exists a positive
constant $C_M$, such that
\begin{equation*}
\displaystyle
\left\|U^k\right\|_{h}+\left\|V^k\right\|_{h}\leq C_M,\;\;k=0,1,\cdots, N.
\end{equation*}
\end{theorem}

\begin{proof}
Computing the discrete inner product of the first equation of (\ref{eq.3.23}) with $U^{k+\frac{1}{2}}$
and using the Lemma \ref{Lem.3.13},
again taking the real part of the resulting equation yield
\begin{equation}\label{eq.3.33}
\begin{aligned}\displaystyle
\frac{1}{2\tau}\left(\left\|U^{k+1}\right\|_h^2-\left\|U^{k}\right\|_h^2\right)
=&-\beta_1\left\|\left|\delta_h^{\frac{\alpha}{2}}U\right|\right\|^2
+\gamma_1\left\|U^{k+1/2}\right\|_h^2
-\mu_1\left\|U^{k+1/2}\right\|_{l_h^4}^4\\
&+\Im\left(\left|U^{k+1/2}\right|^2V^{k+1/2},U^{k+1/2}\right)_h.
\end{aligned}
\end{equation}

For the second equation in (\ref{eq.3.23}), by performing a similar operation, we can also obtain
\begin{equation}\label{eq.3.34}
\begin{aligned}\displaystyle
\frac{1}{2\tau}\left(\left\|V^{k+1}\right\|_h^2-\left\|V^{k}\right\|_h^2\right)
=&-\beta_2\left\|\left|\delta_h^{\frac{\alpha}{2}}V\right|\right\|^2
+\gamma_2\left\|V^{k+1/2}\right\|_h^2
-\mu_2\left\|V^{k+1/2}\right\|_{l_h^4}^4\\
&+\Im\left(\left|V^{k+1/2}\right|^2U^{k+1/2},V^{k+1/2}\right)_h.
\end{aligned}
\end{equation}

For the last terms of equation (\ref{eq.3.33}) and (\ref{eq.3.34}), it follows from the H\"{o}lder inequality and Young's
inequality that there are positive constants $\epsilon_1$,
$\epsilon_2$, $\epsilon_3$ and $\epsilon_4$ that are independent of $\tau$ and $h$, such that
\begin{equation}\label{eq.3.35}
\begin{aligned}\displaystyle
\Im\left(\left|U^{k+1/2}\right|^2V^{k+1/2},U^{k+1/2}\right)_h\leq&
\left\|\left|U^{k+1/2}\right|^2\right\|_h
\left\|V^{k+1/2}\right\|_h\left\|U^{k+1/2}\right\|_h\\
\leq&\epsilon_1^2\left\|U^{k+1/2}\right\|_{l_h^4}^4+
\frac{1}{4\epsilon_1^2}\left(\left\|V^{k+1/2}\right\|_h^2\left\|U^{k+1/2}\right\|_h^2\right)
\\
\leq&\epsilon_1^2\left\|U^{k+1/2}\right\|_{l_h^4}^4+
\frac{1}{4\epsilon_1^2}\left(\epsilon_2^2\left\|V^{k+1/2}\right\|_{l_h^4}^4+
\frac{1}{4\epsilon_2^2}\left\|U^{k+1/2}\right\|_{l_h^4}^4\right)
\\
=&\left(\epsilon_1^2+\frac{1}{16\epsilon_1^2\epsilon_2^2}\right)\left\|U^{k+1/2}\right\|_{l_h^4}^4
+\frac{\epsilon_2^2}{4\epsilon_1^2}\left\|V^{k+1/2}\right\|_{l_h^4}^4,
\end{aligned}
\end{equation}
and
\begin{equation}\label{eq.3.36}
\begin{aligned}\displaystyle
\Im\left(\left|V^{k+1/2}\right|^2U^{k+1/2},V^{k+1/2}\right)_h\leq&
\left(\epsilon_3^2+\frac{1}{16\epsilon_3^2\epsilon_4^2}\right)\left\|V^{k+1/2}\right\|_{l_h^4}^4
+\frac{\epsilon_4^2}{4\epsilon_3^2}\left\|U^{k+1/2}\right\|_{l_h^4}^4.
\end{aligned}
\end{equation}

Denoting $W^k=:\left\|U^{k}\right\|_h^2+\left\|V^{k}\right\|_h^2$, then
combining (\ref{eq.3.33})--(\ref{eq.3.36}) can lead to
\begin{equation}\label{eq.3.37}
\begin{aligned}\displaystyle
W^{k+1}
\leq&W^{k}
+2\tau\gamma_1\left\|U^{k+1/2}\right\|_h^2+
2\tau\left(\epsilon_1^2-\mu_1+\frac{1}{16\epsilon_1^2\epsilon_2^2}+\frac{\epsilon_4^2}{4\epsilon_3^2}\right)
\left\|U^{k+1/2}\right\|_{l_h^4}^4\\
&+2\tau\gamma_2\left\|V^{k+1/2}\right\|_h^2
+2\tau\left(\epsilon_3^2-\mu_2+\frac{1}{16\epsilon_3^2\epsilon_4^2}+\frac{\epsilon_2^2}{4\epsilon_1^2}\right)
\left\|V^{k+1/2}\right\|_{l_h^4}^4\\
\leq&W^{k}
+2\tau\max\left\{|\gamma_1|,|\gamma_2|\right\}\left(\left\|U^{k+1/2}
\right\|_h^2+\left\|V^{k+1/2}\right\|_h^2\right).
\end{aligned}
\end{equation}
Here, the following conditions
\begin{equation}\label{eq.3.38}
\begin{aligned}\displaystyle
\mu_1\geq\epsilon_1^2+\frac{1}{16\epsilon_1^2\epsilon_2^2}+\frac{\epsilon_4^2}{4\epsilon_3^2}
,\;\;
\mu_2\geq\epsilon_3^2+\frac{1}{16\epsilon_3^2\epsilon_4^2}+\frac{\epsilon_2^2}{4\epsilon_1^2},
\end{aligned}
\end{equation}
are used.

Hence, denoting $\gamma=\max\left\{|\gamma_1|,|\gamma_2|\right\}$,
then (\ref{eq.3.37}) can be further rewritten as
\begin{equation}\label{eq.3.39}
\begin{aligned}\displaystyle
W^{k+1}
-W^{k}\leq
2\tau\gamma\left(\left\|U^{k+1/2}\right\|_h^2+\left\|V^{k+1/2}\right\|_h^2\right)
\leq \tau\gamma\left(W^{k+1}+W^{k}\right),\;k=0,1,\cdots,N.
\end{aligned}
\end{equation}

If we choose $\tau\gamma\leq1/2$, we can immediately infer from (\ref{eq.3.39}) that
\begin{equation*}
\begin{aligned}\displaystyle
W^{k}
\leq\exp\left(\frac{2k\tau\gamma}{1-\tau\gamma}\right)W^{0}
\leq\exp\left({4\gamma}T\right)W^{0},
\;k=0,1,\cdots,N.
\end{aligned}\vspace{-0.4cm}
\end{equation*}
This implies that
\begin{equation*}
\displaystyle
\left\|U^k\right\|_{h}+\left\|V^k\right\|_{h}\leq\sqrt{2}\exp\left({2\gamma}T\right)
\left(\left\|U^0\right\|_{h}+\left\|V^0\right\|_{h}\right)
=: C_M,\;k=0,1,\cdots,N.
\end{equation*}
This ends the proof.
\end{proof}

\subsubsection{Existence and uniqueness of solution}
In this section, we will show that the solution of difference scheme (\ref{eq.3.23})--(\ref{eq.3.25})
 exists and is unique.
To this end, we first give several lemmas that play an important role in the follow-up analysis.

\begin{lemma}\label{Lem.3.16}(Brouwder fixed point
theorem \cite{Akrivis}) Let $\left(\mathcal{H},\langle\cdot,\cdot\rangle\right)$ be
a finite-dimensional inner product space,
$\left\|\cdot\right\|_{\mathcal{H}}$ be the associated norm,
and $\mathcal{G}:\mathcal{H}\rightarrow \mathcal{H}$ be continuous.
 Suppose
 \begin{equation*}
\begin{aligned}
\exists \;\zeta>0, \;\forall\; z\in \mathcal{H}, \;\left\|z\right\|_{\mathcal{H}}=\zeta,\;
\Re\left\langle \mathcal{G}(z),z\right\rangle\geq0.
\end{aligned}
\end{equation*}
Then, there exists a $z^*\in \mathcal{H}$ such that $\mathcal{G}(z^*)=0$ and
 $\left\|z^*\right\|_{\mathcal{H}}\leq\zeta$.
\end{lemma}

\begin{lemma}(\cite{Sun})\label{Lem.3.17}
For any complex functions $U, V, u$ and $v$, we have
$$\left|\left|U\right|^2V-\left|u\right|^2v\right|\leq
\left(\max\left\{\left|U\right|,\left|V\right|,
\left|u\right|,\left|v\right|\right\}\right)^2
\left(2\left|U-u\right|+\left|V-v\right|\right).
$$
\end{lemma}
\vspace{-0.5cm}
\begin{theorem}
The solution of system (\ref{eq.3.23})--(\ref{eq.3.25}) exists.
\end{theorem}
\begin{proof}
Let $\mathcal{Z}=\left\{z\,|\,z=[z_1,z_2],\;z_1,
z_2\in\mathcal{\mathring V}_ h\right\}$, and
$z=[z_1,z_2],\;z^{\prime}=[z_1^{\prime},\;z_2^{\prime}]\in\mathcal{Z}$,
define
$$\left<z,z^{\prime}\right>=\left(z_1,z_1^{\prime}\right)_h+\left(z_2,z_2^{\prime}\right)_h,\;
\left\|z\right\|_\mathcal{H}^2=\left\|z_1\right\|^2_{h}+\left\|z_2\right\|^2_{h}.
$$

The following proof process can be completed by mathematical induction method.
Next, we use mathematical induction to prove it.
First of all, it is not difficult to find that $[U^0,V^0]\in\mathcal{Z}$
 satisfies system (\ref{eq.3.23})--(\ref{eq.3.25}). Now suppose that $\left[U^k,V^k\right]\in\mathcal{Z}$
 satisfies system (\ref{eq.3.23})--(\ref{eq.3.25}), then our main task is to prove that there is
 $\left[U^{k+1},V^{k+1}\right]\in\mathcal{Z}$, which
 also satisfies system (\ref{eq.3.23}) -- (\ref{eq.3.25}). For analysis, we will rewrite system (\ref{eq.3.23}) as
 \begin{equation}\label{eq.3.40}
\left\{
\begin{aligned}\displaystyle
U_j^{k+1/2}=&U_j^k+\frac{\tau}{2}\left\{
\left[\left(\beta_1+\mathrm{i}\eta_1\right)\left(\mathcal{H}_h^{\alpha}\right)^{-1}
\mathcal{A}_h^{\alpha}+\gamma_1\right]U_j^{k+1/2}
-\left(\mu_1+\mathrm{i}\zeta_1\right)\left|U_j^{k+1/2}\right|^2U_j^{k+1/2}\right.\\
&\left.+\mathrm{i}\left|U_j^{k+1/2}\right|^2V_j^{k+1/2}\right\},\;j\in\mathds{Z}_M,\;0\leq k\leq N-1,
\\
V_j^{k+1/2}=&V_j^k+\frac{\tau}{2}\left\{
\left[\left(\beta_2+\mathrm{i}\eta_2\right)\left(\mathcal{H}_h^{\alpha}\right)^{-1}
\mathcal{A}_h^{\alpha}+\gamma_2\right]V_j^{k+1/2}
-\left(\mu_2+\mathrm{i}\zeta_2\right)\left|V_j^{k+1/2}\right|^2V_j^{k+1/2}\right.\\
&\left.+\mathrm{i}\left|V_j^{k+1/2}\right|^2U_j^{k+1/2}\right\},\;j\in\mathds{Z}_M,\;0\leq k\leq N-1.
\end{aligned}
\right.
\end{equation}

Denote $z_1=U^{k+1/2},\;z_2=V^{k+1/2}$. Based on the system (\ref{eq.3.40}), the mapping
$\mathcal{G}(z)=[\mathcal{G}_1(z),\mathcal{G}_2(z)]$ on $\mathcal{Z}$
can be defined by
\begin{equation}\label{eq.3.41}
\left\{
\begin{aligned}\displaystyle
\left(\mathcal{G}_1(z)\right)_j=&
(z_1)_j-U_j^k-\frac{\tau}{2}\left\{
\left[\left(\beta_1+\mathrm{i}\eta_1\right)\left(\mathcal{H}_h^{\alpha}\right)^{-1}
\mathcal{A}_h^{\alpha}+\gamma_1\right]U_j^{k+1/2}\right.\\
&\left.
-\left(\mu_1+\mathrm{i}\zeta_1\right)\left|U_j^{k+1/2}\right|^2U_j^{k+1/2}
+\mathrm{i}\left|U_j^{k+1/2}\right|^2V_j^{k+1/2}\right\},
\\
\left(\mathcal{G}_2(z)\right)_j=&
(z_2)_j-V_j^k-\frac{\tau}{2}\left\{
\left[\left(\beta_2+\mathrm{i}\eta_2\right)\left(\mathcal{H}_h^{\alpha}\right)^{-1}
\mathcal{A}_h^{\alpha}+\gamma_2\right]V_j^{k+1/2}\right.\\
&\left.
-\left(\mu_2+\mathrm{i}\zeta_2\right)\left|V_j^{k+1/2}\right|^2V_j^{k+1/2}
+\mathrm{i}\left|V_j^{k+1/2}\right|^2U_j^{k+1/2}\right\}.
\end{aligned}
\right.
\end{equation}

It is not difficult to find that the mapping $\mathcal{G}(z)$ is continuous.
Computing the inner product
of the system (\ref{eq.3.41}) with $z=[z_1,z_2]$, and with help of the Lemma \ref{Lem.3.13}, we can get
\begin{equation}\label{eq.3.42}
\begin{aligned}\displaystyle
\left<\mathcal{G}(z),z\right>=&\left(\mathcal{G}_1(z),z_1\right)_h
+\left(\mathcal{G}_2(z),z_2\right)_h\vspace{0.2cm}\\
=&\displaystyle\left\|z\right\|_\mathcal{H}-\left(U^{k},z_1\right)_h
-\left(V^{k},z_2\right)_h+\frac{\tau}{2}\left[
\left(\beta_1+\mathrm{i}\eta_1\right)\left\|\left|\delta_h^{\frac{\alpha}{2}}z_1\right|\right\|^2
\right.\\
&+\left(\beta_2+\mathrm{i}\eta_2\right)\left\|\left|\delta_h^{\frac{\alpha}{2}}z_2\right|\right\|^2
-\gamma_1\left\|z_1\right\|_h^2
-\gamma_2\left\|z_2\right\|_h^2
+\left(\mu_1+\mathrm{i}\zeta_1\right)\left\|z_1\right\|_{l_h^4}^4
\\
&+\left(\mu_2+\mathrm{i}\zeta_2\right)\left\|z_2\right\|_{l_h^4}^4
-\left.\mathrm{i}\left(\left|z_1\right|^2z_2,z_1\right)-
\mathrm{i}\left(\left|z_2\right|^2z_1,z_2\right)\right].
\end{aligned}
\end{equation}

Further taking the real part of the (\ref{eq.3.42}) and using the Cauchy-Schwarz inequality, it can lead to
\begin{equation}\label{eq.3.43}
\begin{aligned}\displaystyle
\Re\left<\mathcal{G}(z),z\right>
\geq&\left(1-\frac{\tau\gamma_1}{2}\right)\left\|z_1\right\|_h^2+
\left(1-\frac{\tau\gamma_2}{2}\right)\left\|z_2\right\|_h^2
-\Re\left(U^{k},z_1\right)_h-\Re\left(V^{k},z_2\right)_h\\
&+\left(\mu_1-\epsilon_1^2-\frac{1}{16\epsilon_1^2\epsilon_2^2}-\frac{\epsilon_4^2}{4\epsilon_3^2}\right)
\left\|z_1\right\|_{l_h^4}^4+
\left(\mu_2-\epsilon_3^2-\frac{1}{16\epsilon_3^2\epsilon_4^2}-\frac{\epsilon_2^2}{4\epsilon_1^2}\right)
\left\|z_2\right\|_{l_h^4}^4
\\
\geq&\frac{1}{2}\left(1-\tau\gamma_1\right)\left\|z_1\right\|_h^2+
\frac{1}{2}\left(1-\tau\gamma_2\right)\left\|z_2\right\|_h^2
-\frac{1}{2}\left(\left\|U^k\right\|_h^2+\left\|V^k\right\|_h^2\right)
\\
\geq&\frac{1}{2}\left(1-\tau|\gamma_1|\right)\left\|z_1\right\|_h^2+
\frac{1}{2}\left(1-\tau|\gamma_2|\right)\left\|z_2\right\|_h^2
-\frac{1}{2}\left(\left\|U^k\right\|_h^2+\left\|V^k\right\|_h^2\right)
\\
\geq&\frac{1}{2}\left(1-\tau\gamma\right)\left\|z\right\|_\mathcal{H}^2
-\frac{1}{2}\left(\left\|U^k\right\|_h^2+\left\|V^k\right\|_h^2\right)
\\
\geq&\frac{1}{4}\left\|z\right\|_\mathcal{H}^2
-\frac{1}{2}\left(\left\|U^k\right\|_h^2+\left\|V^k\right\|_h^2\right)
\end{aligned}
\end{equation}
by using condition (\ref{eq.3.38}) and choosing a sufficiently small $\tau$ to satisfy
$\tau\gamma=\tau\left\{|\gamma_1|,|\gamma_2|\right\}\leq1/2$.

Finally, take $\zeta^2=2\left(\left\|U^{k}\right\|^2_{h}+
\left\|V^{k}\right\|^2_{h}+1\right)$, then inequality (\ref{eq.3.43}) implies that
$\Re\left<\mathcal{G}(z),z\right>\geq0$
as for $\left\|z\right\|_\mathcal{H}=\zeta$.
It follows from the
Lemma \ref{Lem.3.16} that there is an element $z^{*}=[z_1^{*},z_2^{*}]\in\mathcal{Z}$,
which satisfies $\mathcal{G}(z^*)=0$ and
$\left\|z^{*}\right\|_\mathcal{H}=\zeta$.
So we let  $\left[U^{k+1},V^{k+1}\right]=2\left[z_1^{*},z_2^{*}\right]-\left[U^{k},V^{k}\right]$,
and it can be proved that it satisfies system (\ref{eq.3.23})--(\ref{eq.3.25}). The proof is completed.
\end{proof}

Next, we give its uniqueness result, and it is described in the following theorem.
\begin{theorem}
The difference solution of system (\ref{eq.3.23})--(\ref{eq.3.25}) has uniqueness.
\end{theorem}

\begin{proof}
Suppose $\left(U,V\right)$ and
$\left({\widetilde{U}},{\widetilde{V}}\right)$ are two sets
of solutions of system (\ref{eq.3.23})-- (\ref{eq.3.25}).

Denote
$$\xi_j=U_j^{k+1/2}-{\widetilde{U}}_j^{k+1/2},\;
\widetilde{\xi}_j=V_j^{k+1/2}-\widetilde{V}_j^{k+1/2},
\;0\leq j\leq M-1,\;0\leq k\leq N-1,$$
then from (\ref{eq.3.23}), we can see that $\xi_j$ and $\widetilde{\xi}_j$ satisfy
 \begin{equation}\label{eq.3.44}
\left\{
\begin{aligned}\displaystyle
\frac{2}{\tau}\xi_j=&
\left[\left(\beta_1+\mathrm{i}\eta_1\right)\left(\mathcal{H}_h^{\alpha}\right)^{-1}
\mathcal{A}_h^{\alpha}+\gamma_1\right]\xi_j
-\left(\mu_1+\mathrm{i}\zeta_1\right)S_j^{(1)}+
\mathrm{i}S_j^{(2)},
\\
\frac{2}{\tau}\widetilde{\xi}_j=&
\left[\left(\beta_2+\mathrm{i}\eta_2\right)\left(\mathcal{H}_h^{\alpha}\right)^{-1}
\mathcal{A}_h^{\alpha}+\gamma_2\right]\widetilde{\xi}_j
-\left(\mu_2+\mathrm{i}\zeta_2\right)\widetilde{S}_j^{(1)}+
\mathrm{i}\widetilde{S}_j^{(2)},
\end{aligned}
\right.
\end{equation}
where
\begin{equation*}
\begin{aligned}\displaystyle
S_j^{(1)}=\left|U_j^{k+1/2}\right|^2U_j^{k+1/2}-
\left|\widetilde{U}_j^{k+1/2}\right|^2\widetilde{U}_j^{k+1/2},\;\;
S_j^{(2)}=\left|U_j^{k+1/2}\right|^2V_j^{k+1/2}-
\left|\widetilde{U}_j^{k+1/2}\right|^2\widetilde{V}_j^{k+1/2},\\
\widetilde{S}_j^{(1)}=\left|V_j^{k+1/2}\right|^2V_j^{k+1/2}-
\left|\widetilde{V}_j^{k+1/2}\right|^2\widetilde{V}_j^{k+1/2},\;\;
\widetilde{S}_j^{(2)}=\left|V_j^{k+1/2}\right|^2U_j^{k+1/2}-
\left|\widetilde{V}_j^{k+1/2}\right|^2\widetilde{U}_j^{k+1/2}.
\end{aligned}
\end{equation*}

Based on Lemma \ref{Le.2.4}, we can know that there are positive constants $C_1$,
$C_2$, $\widetilde{C}_1$
 and $\widetilde{C}_2$, such that
\begin{equation}\label{eq.3.45}
\begin{aligned}\displaystyle
\left|S_j^{(1)}\right|\leq C_1\left|\xi_j\right|,\;
\left|S_j^{(2)}\right|\leq C_2\left(\left|\xi_j\right|+\left|\widetilde{\xi}_j\right|\right),\;
\left|\widetilde{S}_j^{(1)}\right|\leq \widetilde{C}_1\left|\widetilde{\xi}_j\right|,\;
\left|\widetilde{S}_j^{(2)}\right|\leq \widetilde{C}_2\left(\left|\xi_j\right|+\left|\widetilde{\xi}_j\right|\right).
\end{aligned}
\end{equation}

Computing the inner product for the first and second equations of (\ref{eq.3.44}) with $\xi$ and $\chi$,
respectively,  and then adding the resulting two equations and taking the real part
yields
 \begin{equation}\label{eq.3.46}
\begin{aligned}\displaystyle
\frac{2}{\tau}\left(\left\|\xi\right\|_h^2+\left\|\widetilde{\xi}\right\|_h^2\right)
\leq&\gamma_1\left\|\xi\right\|_h^2+\gamma_2\left\|\widetilde{\xi}\right\|_h^2
-\Re\left\{
\left(\mu_1+\mathrm{i}\zeta_1\right)\left(S^{(1)},\xi\right)\right\}\\
&-\Re\left\{
\left(\mu_2+\mathrm{i}\zeta_2\right)\left(\widetilde{S}^{(1)},\widetilde{\xi}\right)\right\}
+\Im\left(S^{(2)},{\xi}\right)+\Im\left(\widetilde{S}^{(2)},\widetilde{\xi}\right).
\end{aligned}
\end{equation}

Applying the Cauchy-Schwarz inequality to (\ref{eq.3.46}) and combining with (\ref{eq.3.45}) can finally lead to
 \begin{equation}\label{eq.3.47}
\begin{aligned}\displaystyle
\left\|\xi\right\|_h^2+\left\|\widetilde{\xi}\right\|_h^2
\leq \tau C_{\max}\left(\left\|\xi\right\|_h^2+\left\|\widetilde{\xi}\right\|_h^2\right),
\end{aligned}
\end{equation}
where
 \begin{equation*}
\begin{aligned}\displaystyle
C_{\max}=\max\left\{2\gamma_1+3C_2+\widetilde{C}_2
+2C_1\sqrt{\mu_1^2+\zeta_1^2},
2\gamma_2+3\widetilde{C}_2+{C}_2
+2\widetilde{C}_1\sqrt{\mu_2^2+\zeta_2^2}
\right\}.
\end{aligned}
\end{equation*}

For a sufficiently small time step $\tau$, when $\tau C_{\max}<1$, we can infer from (\ref{eq.3.47}) that
$
\left\|\xi\right\|_h^2+\left\|\widetilde{\xi}\right\|_h^2=0,
$
this further implies that
$U_j^{k+1/2}={\widetilde{U}}_j^{k+1/2}$ and $V_j^{k+1/2}={\widetilde{V}}_j^{k+1/2}$,
which completes the proof.
\end{proof}

\subsubsection{Convergence}
Firstly, we denote the local truncation errors of numerical
 scheme (\ref{eq.3.23})--(\ref{eq.3.25}) as $R_j^{k+1/2}$ and $\widetilde{R}_j^{k+1/2}$.
 Based on the previous process of establishing
 the difference scheme, we can easily obtain the following result:

\begin{lemma}
Assuming that the coupled system (\ref{eq.3.23})--(\ref{eq.3.25}) has solutions $U(x,t)$ and $V(x,t)$ that meet certain requirements,
then we know that there are positive constants $C_R$ and $\widetilde{C}_{\widetilde{R}}$,
such that
\begin{equation*}
\begin{aligned}\displaystyle
\left|R_j^{k+1/2}\right|\leq C_R\left(\tau^2+h^4\right),\;\;
\left|\widetilde{R}_j^{k+1/2}\right|\leq \widetilde{C}_{\widetilde{R}}\left(\tau^2+h^4\right),\;0\leq k\leq N-1.
\end{aligned}
\end{equation*}
\end{lemma}
and further results
\begin{equation}\label{eq.3.48}
\begin{aligned}\displaystyle
\left\|R^{k+1/2}\right\|_{h}^2\leq (b-a)\;C_R^2\left(\tau^2+h^4\right)^2,\;\;
\left\|\widetilde{R}^{k+1/2}\right\|_{h}^2\leq (b-a)\;
\widetilde{C}_{\widetilde{R}}^2\left(\tau^2+h^4\right)^2.
\end{aligned}
\end{equation}

Below, we present the convergence result.

\begin{theorem}\label{Th.3.21}Suppose that problem (\ref{eq.1}) with (\ref{eq.2}) and (\ref{eq.3.15}) has a group
of unique smooth solution $\left(u,v\right)=:w\in
C^{4}\left(0,T; \mathscr{C}^{{4+\alpha}}\left(\Omega\right)\right)$ and $\left\{\left(U_j^k,V_j^k\right)\,|\,0\leq j\leq M,
0\leq k\leq N\right\}$ is the solution of system (\ref{eq.3.23})--(\ref{eq.3.25}).
Let's denote the error functions $E_j^k$ and $\widetilde{E}_j^k$ as
\begin{equation*}
\begin{aligned}\displaystyle
E_j^k=U_j^k-u(x_j,t_k),\;\;
\widetilde{E}_j^k=V_j^k-v(x_j,t_k),\;0\leq j\leq M,\;0\leq k\leq N.
\end{aligned}
\end{equation*}
Then, there exists a positive constant $C$ such that
\begin{eqnarray*}
\begin{array}{lll}
\displaystyle \left\|E^k\right\|_{h}+\left\|\widetilde{E}^k\right\|_{h}
\leq C\left(\tau^2+h^4\right),\;\;0\leq k\leq N.
\end{array}
\end{eqnarray*}
\end{theorem}

\begin{proof}
Subtracting (\ref{eq.3.23}) from (\ref{eq.1}), we get the error
system as follows
\begin{equation}\label{eq.3.49}
\left\{
\begin{aligned}\displaystyle
&\delta_tE_j^{k+1/2}
-\left[\left(\beta_1+\mathrm{i}\eta_1\right)\left(\mathcal{H}_h^{\alpha}\right)^{-1}
\mathcal{A}_h^{\alpha}+\gamma_1\right]E_j^{k+1/2}
+\left(\mu_1+\mathrm{i}\zeta_1\right)\left|U_j^{k+1/2}\right|^2E_j^{k+1/2}
\\&+\left(\mu_1+\mathrm{i}\zeta_1\right)X_j^{k+1/2}
-\mathrm{i}Y_j^{k+1/2}=R_j^{k+1/2},\;j\in\mathds{Z}_M,\;0\leq k\leq N-1,\vspace{0.2cm}
\\
&\delta_t\widetilde{E}_j^{k+1/2}
-\left[\left(\beta_2+\mathrm{i}\eta_2\right)\left(\mathcal{H}_h^{\alpha}\right)^{-1}
\mathcal{A}_h^{\alpha}+\gamma_2\right]\widetilde{E}_j^{k+1/2}
+\left(\mu_2+\mathrm{i}\zeta_2\right)\left|V_j^{k+1/2}\right|^2\widetilde{E}_j^{k+1/2}\\
&
+\left(\mu_2+\mathrm{i}\zeta_2\right)\widetilde{X}_j^{k+1/2}
-\mathrm{i}\widetilde{Y}_j^{k+1/2}=\widetilde{R}_j^{k+1/2},\;j\in\mathds{Z}_M,\;0\leq k\leq N-1,
\end{aligned}
\right.
\end{equation}
where\vspace{-0.2cm}
\begin{equation*}
\begin{aligned}\displaystyle
&X_j^{k+1/2}=\left|U_j^{k+1/2}\right|^2u\left(x_j,t_{k+1/2}\right)-
\left|u\left(x_j,t_{k+1/2}\right)\right|^2u\left(x_j,t_{k+1/2}\right),\\
&Y_j^{k+1/2}=\left|U_j^{k+1/2}\right|^2V_j^{k+1/2}-
\left|u\left(x_j,t_{k+1/2}\right)\right|^2v\left(x_j,t_{k+1/2}\right),\\
&\widetilde{X}_j^{k+1/2}=\left|V_j^{k+1/2}\right|^2v\left(x_j,t_{k+1/2}\right)-
\left|v\left(x_j,t_{k+1/2}\right)\right|^2v\left(x_j,t_{k+1/2}\right),\\
&\widetilde{Y}_j^{k+1/2}=\left|V_j^{k+1/2}\right|^2U_j^{k+1/2}-
\left|v\left(x_j,t_{k+1/2}\right)\right|^2u\left(x_j,t_{k+1/2}\right).
\end{aligned}
\end{equation*}

Denoting $c_u=\sup_{0\leq t\leq T,a\leq x\leq b}|u(x,t)|$ and
$c_v=\sup_{0\leq t\leq T,a\leq x\leq b}|v(x,t)|$,
 it follows from the Lemma \ref{Lem.3.17} that
\begin{equation}\label{eq.3.50}
\begin{aligned}\displaystyle
\left|X_j^{k+1/2}\right|=&\left|\left|U_j^{k+1/2}\right|^2u\left(x_j,t_{k+1/2}\right)-
\left|u\left(x_j,t_{k+1/2}\right)\right|^2u\left(x_j,t_{k+1/2}\right)\right|\\
=&\left|\left(\left|U_j^{k+1/2}\right|-\left|u\left(x_j,t_{k+1/2}\right)\right|\right)
\left(\left|U_j^{k+1/2}\right|+\left|u\left(x_j,t_{k+1/2}\right)\right|\right)
\left|u\left(x_j,t_{k+1/2}\right)\right|\right|\\
\leq&
\left[\left|E_j^{k+1/2}\right|+
2\left|u\left(x_j,t_{k+1/2}\right)\right|\right]\left|E_j^{k+1/2}\right|
\left|u\left(x_j,t_{k+1/2}\right)\right|
\\
\leq&
c_u
\left|E_j^{k+1/2}\right|
\left[2c_u+\left|E_j^{k+1/2}\right|\right]
\\
=&
2c_u^2
\left|E_j^{k+1/2}\right|+
c_u\left|E_j^{k+1/2}\right|^2,
\end{aligned}
\end{equation}
\begin{equation}\label{eq.3.51}
\begin{aligned}\displaystyle
\left|Y_j^{k+1/2}\right|
=&\left|\left|U_j^{k+1/2}\right|^2V_j^{k+1/2}-
\left|u\left(x_j,t_{k+1/2}\right)\right|^2v\left(x_j,t_{k+1/2}\right)\right|\\
\leq&
\left(\max\left\{\left|U_j^{k+1/2}\right|,\left|V_j^{k+1/2}\right|,
\left|u\left(x_j,t_{k+1/2}\right)\right|,\left|v\left(x_j,t_{k+1/2}\right)\right|\right\}\right)^2
\left(2\left|{E}_j^{k+1/2}\right|+\left|\widetilde{E}_j^{k+1/2}\right|\right)
\\
\leq&
\left(\left|U_j^{k+1/2}\right|+\left|V_j^{k+1/2}\right|+
\left|u\left(x_j,t_{k+1/2}\right)\right|+\left|v\left(x_j,t_{k+1/2}\right)\right|\right)^2
\left(2\left|{E}_j^{k+1/2}\right|+\left|\widetilde{E}_j^{k+1/2}\right|\right)
\\
\leq&
\left(\left|E_j^{k+1/2}\right|+\left|\widetilde{E}_j^{k+1/2}\right|+
2\left|u\left(x_j,t_{k+1/2}\right)\right|+2\left|v\left(x_j,t_{k+1/2}\right)\right|\right)^2
\left(2\left|{E}_j^{k+1/2}\right|+\left|\widetilde{E}_j^{k+1/2}\right|\right)
\\
\leq&
2\left[\left(\left|E_j^{k+1/2}\right|+\left|\widetilde{E}_j^{k+1/2}\right|\right)^2+
4\left(\left|u\left(x_j,t_{k+1/2}\right)\right|+\left|v\left(x_j,t_{k+1/2}\right)\right|\right)^2\right]
\left(2\left|{E}_j^{k+1/2}\right|+\left|\widetilde{E}_j^{k+1/2}\right|\right)
\\
\leq&
4\left[\left|E_j^{k+1/2}\right|^2+\left|\widetilde{E}_j^{k+1/2}\right|^2+
2\left(c_u+c_v\right)^2\right]
\left(2\left|{E}_j^{k+1/2}\right|+\left|\widetilde{E}_j^{k+1/2}\right|\right),
\end{aligned}
\end{equation}
\begin{equation}\label{eq.3.52}
\begin{aligned}\displaystyle
\left|\widetilde{X}_j^{k+1/2}\right|
\leq2c_v^2
\left|\widetilde{E}_j^{k+1/2}\right|+
c_v\left|\widetilde{E}_j^{k+1/2}\right|^2,
\end{aligned}
\end{equation}
and
\begin{equation}\label{eq.3.53}
\begin{aligned}\displaystyle
\left|\widetilde{Y}_j^{k+1/2}\right|
\leq&
4\left[\left|E_j^{k+1/2}\right|^2+\left|\widetilde{E}_j^{k+1/2}\right|^2+
2\left(c_u+c_v\right)^2\right]
\left(2\left|\widetilde{E}_j^{k+1/2}\right|+\left|{E}_j^{k+1/2}\right|\right).
\end{aligned}
\end{equation}

In addition, it follows from Theorem \ref{Th.3.15} that
\begin{equation}\label{eq.3.54}
\begin{aligned}\displaystyle
\left\|E^{k+1/2}\right\|_h\leq \left\|U^{k+1/2}\right\|_h+\left\|u\left(x_j,t_{k+1/2}\right)\right\|_h
\leq C_M+(b-a)c_u=:C_E,
\end{aligned}
\end{equation}
and
\begin{equation}\label{eq.3.55}
\begin{aligned}\displaystyle
\left\|\widetilde{E}^{k+1/2}\right\|_h\leq \left\|V^{k+1/2}\right\|_h+\left\|v\left(x_j,t_{k+1/2}\right)\right\|_h
\leq C_M+(b-a)c_v=:C_{\widetilde{E}}.
\end{aligned}
\end{equation}

Next, for the sake of clarity, we will proceed in three steps in detail.

{\bf Step 1: Deal with the first equation in (\ref{eq.3.49})}.
For the first equation of (\ref{eq.3.49}), calculating the discrete inner product with respect
to $E^{k+1/2}$ and taking the real part of the obtained result can ultimately lead to
\begin{equation}\label{eq.3.56}
\begin{aligned}
&\frac{1}{2\tau}\left(\left\|E^{k+1}\right\|_h^2-\left\|E^{k}\right\|_h^2\right)
+\beta_1\left\|\left|\delta_h^{\frac{\alpha}{2}}E\right|\right\|^2
-\gamma_1\left\|E^{k+1/2}\right\|_h^2\\
&=-
\Re\left\{\left(\mu_1+\mathrm{i}\zeta_1\right)X^{k+1/2},E^{k+1/2}\right\}
+\Im\left(Y^{k+1/2},E^{k+1/2}\right)+\Re\left(R^{k+1/2},E^{k+1/2}\right)
\\
&=:\mathbb{I}_1+\mathbb{I}_2+\mathbb{I}_3.
\end{aligned}\vspace{-0.2cm}
\end{equation}

Now, let's begin to analyze each term on the right side of equation (\ref{eq.3.56}). First of all, for the
first item on the right, using (\ref{eq.3.50}), we have
\begin{equation}\label{eq.3.57}
\begin{aligned}\displaystyle
\left|\mathbb{I}_1\right|=&\left|-\Re\left\{\left(\mu_1+
\mathrm{i}\zeta_1\right)X^{k+1/2},E^{k+1/2}\right\}\right|
=\left|-\Re\left\{\left(\mu_1+\mathrm{i}\zeta_1\right)h\sum_{j=1}^{M-1}X_j^{k+1/2}
\left(E_j^{k+1/2}\right)^{*}\right\}\right|
\\
\leq&\sqrt{\mu_1^2+\zeta_1^2}\,h\sum_{j=1}^{M-1}\left|X_j^{k+1/2}\right|
\left|E_j^{k+1/2}\right|
\leq\sqrt{\mu_1^2+\zeta_1^2}\,h\sum_{j=1}^{M-1}
\left(
2c_u^2
\left|E_j^{k+1/2}\right|^2+
c_u\left|E_j^{k+1/2}\right|^3\right)\\
=&\sqrt{\mu_1^2+\zeta_1^2}\left(
2c_u^2
\left\|E^{k+1/2}\right\|_h^2+
c_u\left\|E^{k+1/2}\right\|_{l_h^3}^3\right).
\end{aligned}
\end{equation}

Based on Lemma \ref{Lem.3.1} and \ref{Lem.3.2}, and the Young's inequality, we get
\begin{equation}\label{eq.3.58}
\begin{aligned}\displaystyle
\left\|E^{k+1/2}\right\|_{l_h^3}^3\leq& C_3\left\|E^{k+1/2}\right\|_{H_h^{\frac{\alpha}{2}}}
\left\|E^{k+1/2}\right\|_{{h}}^{2}
\leq C_3\left(\epsilon_5^2\left\|E^{k+1/2}\right\|_{H_h^{\frac{\alpha}{2}}}^2+\frac{1}{4\epsilon_5^2}
\left\|E^{k+1/2}\right\|_{{h}}^{4}\right)\\
=&C_3\left(
\epsilon_5^2\left\|E^{k+1/2}\right\|_{h}^2+\epsilon_5^2\left|E^{k+1/2}\right|_{H_h^{\frac{\alpha}{2}}}^2+\frac{1}{4\epsilon_5^2}
\left\|E^{k+1/2}\right\|_{{h}}^{4}\right)\\\leq&
C_3\left[-\frac{\epsilon_5^2}{C_1(\alpha)}\left\|\left|\delta_h^{\frac{\alpha}{2}}E^{k+1/2}\right|\right\|^2+
\left(\epsilon_5^2+
\frac{C_E^2}{4\epsilon_5^2}\right)\left\|E^{k+1/2}\right\|_{h}^2\right],
\end{aligned}
\end{equation}
where
$
C_3=\frac{\sqrt[4]{24}}{4\alpha}\sqrt{6\alpha\mathrm{B}\left(\frac{3}{\alpha},3-\frac{3}{\alpha}\right)}.
$
Substituting (\ref{eq.3.58}) into (\ref{eq.3.57}) can ultimately lead to
\begin{equation}\label{eq.3.59}
\begin{aligned}\displaystyle
|\mathbb{I}_1|
\leq&-\frac{c_uC_3\epsilon_5^2\sqrt{\mu_1^2+\zeta_1^2}}{C_1(\alpha)}
\left\|\left|\delta_h^{\frac{\alpha}{2}}E^{k+1/2}\right|\right\|^2\\
&+c_u\sqrt{\mu_1^2+\zeta_1^2}\left[2c_u+C_3\left(\epsilon_5^2+
\frac{C_E^2}{4\epsilon_5^2}\right)\left\|E^{k+1/2}\right\|_{h}^2\right]\\
=:&C_4\left\|\left|\delta_h^{\frac{\alpha}{2}}E^{k+1/2}\right|\right\|^2
+C_5\left\|E^{k+1/2}\right\|_{h}^2.
\end{aligned}
\end{equation}

Secondly, let's analyze the second term at the right end of equation (\ref{eq.3.56}).
With the help of (\ref{eq.3.51}), we can see that
\begin{equation}\label{eq.3.60}
\begin{aligned}\displaystyle
|\mathbb{I}_2|=&\left|\Im\left(Y^{k+1/2},E^{k+1/2}\right)\right|\leq
h\sum_{j=1}^{M-1}\left|Y_j^{k+1/2}\right|
\left|E_j^{k+1/2}\right|\\\leq&
4h\sum_{j=1}^{M-1}\left[\left|E_j^{k+1/2}\right|^2+\left|\widetilde{E}_j^{k+1/2}\right|^2+
2\left(c_u+c_v\right)^2\right]
\left(2\left|{E}_j^{k+1/2}\right|+\left|\widetilde{E}_j^{k+1/2}\right|\right)
\left|E_j^{k+1/2}\right|\\=&
8\left(c_u+c_v\right)^2h\sum_{j=1}^{M-1}
\left(2\left|{E}_j^{k+1/2}\right|+\left|\widetilde{E}_j^{k+1/2}\right|\right)
\left|E_j^{k+1/2}\right|+\\
&
4h\sum_{j=1}^{M-1}\left|E_j^{k+1/2}\right|^2
\left(2\left|{E}_j^{k+1/2}\right|+\left|\widetilde{E}_j^{k+1/2}\right|\right)
\left|E_j^{k+1/2}\right|\\
&+
4h\sum_{j=1}^{M-1}\left|\widetilde{E}_j^{k+1/2}\right|^2
\left(2\left|{E}_j^{k+1/2}\right|+\left|\widetilde{E}_j^{k+1/2}\right|\right)
\left|E_j^{k+1/2}\right|\\
=&
2\left(c_u+c_v\right)^2\left(8+\frac{1}{\epsilon_6^2}\right)\left\|{E}^{k+1/2}\right\|_h^2+
8\left(c_u+c_v\right)^2\epsilon_6^2\left\|{\widetilde{E}}^{k+1/2}\right\|_h^2+
8\left\|E^{k+1/2}\right\|_{l_h^4}^4\\
&+4h\sum_{j=1}^{M-1}\left|E_j^{k+1/2}\right|^2
\left|\widetilde{E}_j^{k+1/2}\right|
\left|E_j^{k+1/2}\right|+8h\sum_{j=1}^{M-1}\left|\widetilde{E}_j^{k+1/2}\right|^2
\left|{E}_j^{k+1/2}\right|^2\\
&
+4h\sum_{j=1}^{M-1}\left|\widetilde{E}_j^{k+1/2}\right|^2
\left|\widetilde{E}_j^{k+1/2}\right|
\left|E_j^{k+1/2}\right|\\
=:&2\left(c_u+c_v\right)^2\left(8+\frac{1}{\epsilon_6^2}\right)\left\|{E}^{k+1/2}\right\|_h^2+
8\left(c_u+c_v\right)^2\epsilon_6^2\left\|{\widetilde{E}}^{k+1/2}\right\|_h^2+
8\left\|E^{k+1/2}\right\|_{l_h^4}^4\\
&+\mathbb{II}_1+\mathbb{II}_2+\mathbb{II}_3.
\end{aligned}
\end{equation}
Due to
\begin{equation}\label{eq.3.61}
\begin{aligned}\displaystyle
\mathbb{II}_1=&4h\sum_{j=1}^{M-1}\left|E_j^{k+1/2}\right|^2
\left|\widetilde{E}_j^{k+1/2}\right|
\left|E_j^{k+1/2}\right|\\
\leq&4h\sum_{j=1}^{M-1}
\left(\epsilon_7^2\left|E_j^{k+1/2}\right|^4
+\frac{1}{4\epsilon_7^2}\left|\widetilde{E}_j^{k+1/2}\right|^2
\left|E_j^{k+1/2}\right|^2\right)\\
\leq&4h\sum_{j=1}^{M-1}
\left[\epsilon_7^2\left|E_j^{k+1/2}\right|^4
+\frac{1}{4\epsilon_7^2}\left(\epsilon_8^2\left|\widetilde{E}_j^{k+1/2}\right|^4+\frac{1}{4\epsilon_8^2}
\left|E_j^{k+1/2}\right|^4\right)\right]\\
=&\frac{\epsilon_8^2}{\epsilon_7^2}\left\|\widetilde{E}^{k+1/2}\right\|^4_{l_h^4}+
4\left(\epsilon_7^2+\frac{1}{16\epsilon_7^2\epsilon_8^2}\right)\left\|{E}^{k+1/2}\right\|^4_{l_h^4},
\end{aligned}
\end{equation}
\begin{equation}\label{eq.3.62}
\begin{aligned}\displaystyle
\mathbb{II}_2=&8h\sum_{j=1}^{M-1}
\left|\widetilde{E}_j^{k+1/2}\right|^2
\left|E_j^{k+1/2}\right|^2\leq8h\sum_{j=1}^{M-1}
\left(\epsilon_9^2\left|\widetilde{E}_j^{k+1/2}\right|^4
+\frac{1}{4\epsilon_9^2}\left|{E}_j^{k+1/2}\right|^4\right)\\
=&8\epsilon_9^2\left\|\widetilde{E}^{k+1/2}\right\|^4_{l_h^4}+
\frac{2}{\epsilon_9^2}\left\|{E}^{k+1/2}\right\|^4_{l_h^4},
\end{aligned}
\end{equation}
and
\begin{equation}\label{eq.3.63}
\begin{aligned}\displaystyle
\mathbb{II}_3=&4h\sum_{j=1}^{M-1}\left|\widetilde{E}_j^{k+1/2}\right|^2
\left|\widetilde{E}_j^{k+1/2}\right|
\left|E_j^{k+1/2}\right|\leq
4h\sum_{j=1}^{M-1}
\left(\epsilon_{10}^2\left|\widetilde{E}_j^{k+1/2}\right|^4
+\frac{1}{4\epsilon_{10}^2}\left|\widetilde{E}_j^{k+1/2}\right|^2
\left|E_j^{k+1/2}\right|^2\right)\\
\leq&4h\sum_{j=1}^{M-1}
\left[\epsilon_{10}^2\left|\widetilde{E}_j^{k+1/2}\right|^4
+\frac{1}{4\epsilon_{10}^2}\left(\epsilon_{11}^2\left|{E}_j^{k+1/2}\right|^4+\frac{1}{4\epsilon_{11}^2}
\left|\widetilde{E}_j^{k+1/2}\right|^4\right)\right]\\
=&\frac{\epsilon_{11}^2}{\epsilon_{10}^2}\left\|{E}^{k+1/2}\right\|^4_{l_h^4}+
4\left(\epsilon_{10}^2+\frac{1}{16\epsilon_{10}^2\epsilon_{11}^2}\right)
\left\|{\widetilde{E}}^{k+1/2}\right\|^4_{l_h^4}.
\end{aligned}
\end{equation}
Substituting (\ref{eq.3.61})--(\ref{eq.3.63}) into (\ref{eq.3.60}) further leads to
\begin{equation}\label{eq.3.64}
\begin{aligned}\displaystyle
|\mathbb{I}_2|\leq&
2\left(c_u+c_v\right)^2\left(8+\frac{1}{\epsilon_6^2}\right)\left\|{E}^{k+1/2}\right\|_h^2+
8\left(c_u+c_v\right)^2\epsilon_6^2\left\|{\widetilde{E}}^{k+1/2}\right\|_h^2+
8\left\|E^{k+1/2}\right\|_{l_h^4}^4\\
&+\mathbb{II}_1+\mathbb{II}_2+\mathbb{II}_3\\
\leq&2\left(c_u+c_v\right)^2\left(8+\frac{1}{\epsilon_6^2}\right)\left\|{E}^{k+1/2}\right\|_h^2+
8\left(c_u+c_v\right)^2\epsilon_6^2\left\|{\widetilde{E}}^{k+1/2}\right\|_h^2+
8\left\|E^{k+1/2}\right\|_{l_h^4}^4\\
&+\frac{\epsilon_8^2}{\epsilon_7^2}\left\|\widetilde{E}^{k+1/2}\right\|^4_{l_h^4}+
4\left(\epsilon_7^2+\frac{1}{16\epsilon_7^2\epsilon_8^2}\right)\left\|{E}^{k+1/2}\right\|^4_{l_h^4}
+8\epsilon_9^2\left\|\widetilde{E}^{k+1/2}\right\|^4_{l_h^4}+
\frac{2}{\epsilon_9^2}\left\|{E}^{k+1/2}\right\|^4_{l_h^4}\\&
+\frac{\epsilon_{11}^2}{\epsilon_{10}^2}\left\|{E}^{k+1/2}\right\|^4_{l_h^4}+
4\left(\epsilon_{10}^2+\frac{1}{16\epsilon_{10}^2\epsilon_{11}^2}\right)
\left\|{\widetilde{E}}^{k+1/2}\right\|^4_{l_h^4}\\
\leq&
2\left(c_u+c_v\right)^2\left(8+\frac{1}{\epsilon_6^2}\right)\left\|{E}^{k+1/2}\right\|_h^2+
8\left(c_u+c_v\right)^2\epsilon_6^2\left\|{\widetilde{E}}^{k+1/2}\right\|_h^2\\
&+\left(8+4\epsilon_7^2+\frac{2}{\epsilon_9^2}+\frac{\epsilon_{11}^2}{\epsilon_{10}^2}
+\frac{1}{4\epsilon_7^2\epsilon_8^2}\right)\left\|{E}^{k+1/2}\right\|^4_{l_h^4}\\&
+\left(8\epsilon_9^2+4\epsilon_{10}^2+\frac{\epsilon_8^2}{\epsilon_7^2}
+\frac{1}{4\epsilon_{10}^2\epsilon_{11}^2}\right)\left\|{\widetilde{E}}^{k+1/2}\right\|^4_{l_h^4}\\
=:&C_6\left\|{E}^{k+1/2}\right\|_h^2+C_7\left\|{\widetilde{E}}^{k+1/2}\right\|_h^2
+C_8\left\|{E}^{k+1/2}\right\|^4_{l_h^4}+C_9\left\|{\widetilde{E}}^{k+1/2}\right\|^4_{l_h^4}.
\end{aligned}
\end{equation}
In view of Lemma \ref{Lem.3.2} with $\frac{1}{4}<\sigma_0\leq\frac{\alpha}{2}$ and
using (\ref{eq.3.54}) and Young's inequality, one has
\begin{equation}\label{eq.3.65}
\begin{aligned}\displaystyle
\left\|{E}^{k+1/2}\right\|^4_{l_h^4}\leq &C_6
\left\|{E}^{k+1/2}\right\|_{H_h^{\frac{\alpha}{2}}}^{\frac{8\sigma_0}{\alpha}}\;
\left\|{E}^{k+1/2}\right\|_{{h}}^{4\left(1-\frac{2\sigma_0}{\alpha}\right)}\\
\leq& C_6\left[\frac{4\sigma_0}{\alpha}\left\|{E}^{k+1/2}\right\|_{H_h^{\frac{\alpha}{2}}}^{2}
+\frac{\alpha-4\sigma_0}{\alpha}
\left\|{E}^{k+1/2}\right\|_{{h}}^{\frac{4\left(\alpha-2\sigma_0\right)}{\alpha-4\sigma_0}}
\right]\\
=&
C_6\left[\frac{4\sigma_0}{\alpha}\left(\left|{E}^{k+1/2}\right|_{H_h^{\frac{\alpha}{2}}}^{2}
+\left\|{E}^{k+1/2}\right\|_{h}^{2}
\right)
+\frac{\alpha-4\sigma_0}{\alpha}
\left\|{E}^{k+1/2}\right\|_{{h}}^{2}
\left\|{E}^{k+1/2}\right\|_{{h}}^{\frac{2\alpha}{\alpha-4\sigma_0}}
\right]\\
\leq&
C_6\left[-\frac{4\sigma_0}{\alpha C_1(\alpha)}\left\|\left|
\delta_h^{\frac{\alpha}{2}}E^{k+1/2}\right|\right\|^2
+\frac{4\sigma_0}{\alpha}\left\|{E}^{k+1/2}\right\|_{h}^{2}\right.\\&\left.
+\frac{\alpha-4\sigma_0}{\alpha}\left(C_E\right)^{\frac{2\alpha}{\alpha-4\sigma_0}}
\left\|{E}^{k+1/2}\right\|_{{h}}^{2}
\right]\\
=:&C_{10}\left\|\left|\delta_h^{\frac{\alpha}{2}}E^{k+1/2}\right|\right\|^2
+C_{11}\left\|{E}^{k+1/2}\right\|_{h}^{2},
\end{aligned}
\end{equation}
and
\begin{equation}\label{eq.3.66}
\begin{aligned}\displaystyle
\left\|{\widetilde{E}}^{k+1/2}\right\|^4_{l_h^4}\leq C_{10}
\left\|\left|\delta_h^{\frac{\alpha}{2}}\widetilde{E}^{k+1/2}\right|\right\|^2
+C_{11}\left\|{\widetilde{E}}^{k+1/2}\right\|_{h}^{2},
\end{aligned}
\end{equation}
where
$
C_6=\frac{\sqrt[3]{6}}{27\sigma_0}
4^{\frac{3\alpha-2\sigma_0}{\alpha}}
\mathrm{B}\left(\frac{1}{2\sigma_0},2-\frac{1}{2\sigma_0}\right).
$

Combining (\ref{eq.3.64}), (\ref{eq.3.65}) and (\ref{eq.3.66}), we have
\begin{equation}\label{eq.3.67}
\begin{aligned}\displaystyle
|\mathbb{I}_2|\leq&C_8C_{10}\left\|\left|\delta_h^{\frac{\alpha}{2}}E^{k+1/2}\right|\right\|^2+
C_9C_{10}\left\|\left|\delta_h^{\frac{\alpha}{2}}\widetilde{E}^{k+1/2}\right|\right\|^2\\
&+
\left(C_6+C_8C_{11}\right)\left\|{E}^{k+1/2}\right\|_h^2
+\left(C_7+C_9C_{11}\right)\left\|{\widetilde{E}}^{k+1/2}\right\|_h^2\\
=:&C_{12}\left\|\left|\delta_h^{\frac{\alpha}{2}}E^{k+1/2}\right|\right\|^2+
C_{13}\left\|\left|\delta_h^{\frac{\alpha}{2}}\widetilde{E}^{k+1/2}\right|\right\|^2
+
C_{14}\left\|{E}^{k+1/2}\right\|_h^2
+C_{15}\left\|{\widetilde{E}}^{k+1/2}\right\|_h^2.
\end{aligned}
\end{equation}

Once again, let's deal with the last term on the right side of equation (\ref{eq.3.56}).
Obviously, we know that
\begin{equation}\label{eq.3.68}
\begin{aligned}\displaystyle
\left|\mathbb{I}_3\right|\leq\left\|R^{k+1/2}\right\|_h\left\|E^{k+1/2}\right\|_h
\leq\frac{1}{2}\left(\left\|R^{k+1/2}\right\|_h^2+\left\|E^{k+1/2}\right\|_h^2\right).
\end{aligned}
\end{equation}

Finally, substituting (\ref{eq.3.59}), (\ref{eq.3.67}) and (\ref{eq.3.68}) into (\ref{eq.3.56}) and simplifying
them to obtain the final result
\begin{equation}\label{eq.3.69}
\begin{aligned}\displaystyle
&\frac{1}{2\tau}\left(\left\|E^{k+1}\right\|_h^2-\left\|E^{k}\right\|_h^2\right)
+\beta_1\left\|\left|\delta_h^{\frac{\alpha}{2}}E^{k+1/2}\right|\right\|^2
-\gamma_1\left\|E^{k+1/2}\right\|_h^2
\leq\left(C_4+C_{12}\right)\left\|\left|\delta_h^{\frac{\alpha}{2}}E^{k+1/2}\right|\right\|^2\\
&+C_{13}\left\|\left|\delta_h^{\frac{\alpha}{2}}\widetilde{E}^{k+1/2}\right|\right\|^2
+\left(\frac{1}{2}+C_5+C_{14}\right)\left\|{E}^{k+1/2}\right\|_h^2
+C_{15}\left\|{\widetilde{E}}^{k+1/2}\right\|_h^2
+\frac{1}{2}\left\|R^{k+1/2}\right\|_h^2.
\end{aligned}
\end{equation}

{\bf Step 2: Deal with the second equation in (\ref{eq.3.49})}.
Using the same method as in the step 1, we have
\begin{equation}\label{eq.3.70}
\begin{aligned}\displaystyle
&\frac{1}{2\tau}\left(\left\|\widetilde{E}^{k+1}\right\|_h^2-\left\|\widetilde{E}^{k}\right\|_h^2\right)
+\beta_2\left\|\left|\delta_h^{\frac{\alpha}{2}}\widetilde{E}^{k+1/2}\right|\right\|^2
-\gamma_2\left\|\widetilde{E}^{k+1/2}\right\|_h^2\\
&=-
\Re\left\{\left(\mu_2+\mathrm{i}\zeta_2\right)\widetilde{X}^{k+1/2},\widetilde{E}^{k+1/2}\right\}
+\Im\left(\widetilde{Y}^{k+1/2},\widetilde{E}^{k+1/2}\right)+\Re\left(\widetilde{R}^{k+1/2},\widetilde{E}^{k+1/2}\right)
\\
&=:{\mathbb{\widetilde{I}}_1}+{\mathbb{{\widetilde{I}}}_2}+{\mathbb{{\widetilde{I}}}_3}.
\end{aligned}
\end{equation}

Similarly, using the (\ref{eq.3.52}) and (\ref{eq.3.53}), ignoring the detailed process, we can directly obtain the following results
\begin{equation}\label{eq.3.71}
\begin{aligned}\displaystyle
|{\widetilde{\mathbb{I}}_1}|
\leq C_{16}\left\|\left|\delta_h^{\frac{\alpha}{2}}\widetilde{E}^{k+1/2}\right|\right\|^2
+C_{17}\left\|\widetilde{E}^{k+1/2}\right\|_{h}^2,
\end{aligned}
\end{equation}
\begin{equation}\label{eq.3.72}
\begin{aligned}\displaystyle
|\widetilde{\mathbb{I}}_2|\leq C_{18}\left\|\left|\delta_h^{\frac{\alpha}{2}}\widetilde{E}^{k+1/2}\right|\right\|^2+
C_{19}\left\|\left|\delta_h^{\frac{\alpha}{2}}{E}^{k+1/2}\right|\right\|^2
+
C_{20}\left\|{\widetilde{E}}^{k+1/2}\right\|_h^2
+C_{21}\left\|{{E}}^{k+1/2}\right\|_h^2,
\end{aligned}
\end{equation}
and
\begin{equation}\label{eq.3.73}
\begin{aligned}\displaystyle
\left|{\widetilde{\mathbb{I}}_3}\right|\leq\left\|\widetilde{R}^{k+1/2}\right\|_h\left\|\widetilde{E}^{k+1/2}\right\|_h
\leq\frac{1}{2}\left(\left\|\widetilde{R}^{k+1/2}\right\|_h^2+\left\|\widetilde{E}^{k+1/2}\right\|_h^2\right).
\end{aligned}
\end{equation}

Substituting (\ref{eq.3.71})--(\ref{eq.3.73}) into (\ref{eq.3.70}) can eventually lead to
\begin{equation}\label{eq.3.74}
\begin{aligned}\displaystyle
&\frac{1}{2\tau}\left(\left\|\widetilde{E}^{k+1}\right\|_h^2-\left\|\widetilde{E}^{k}\right\|_h^2\right)
+\beta_2\left\|\left|\delta_h^{\frac{\alpha}{2}}\widetilde{E}^{k+1/2}\right|\right\|^2
-\gamma_2\left\|\widetilde{E}^{k+1/2}\right\|_h^2\leq
\left(C_{16}+C_{18}\right)\left\|\left|\delta_h^{\frac{\alpha}{2}}\widetilde{E}^{k+1/2}\right|\right\|^2\\&
+
C_{19}\left\|\left|\delta_h^{\frac{\alpha}{2}}{E}^{k+1/2}\right|\right\|^2
+\left(\frac{1}{2}+C_{17}+C_{20}\right)\left\|\widetilde{E}^{k+1/2}\right\|_{h}^2
+C_{21}\left\|{{E}}^{k+1/2}\right\|_h^2+\frac{1}{2}\left\|\widetilde{R}^{k+1/2}\right\|_h^2.
\end{aligned}
\end{equation}

{\bf Step 3: Deal with the results obtained in steps 1 and 2}.
Adding equations (\ref{eq.3.69}) and (\ref{eq.3.74}), we have
\begin{equation}\label{eq.3.75}
\begin{aligned}\displaystyle
&\frac{1}{2\tau}\left[\left(\left\|E^{k+1}\right\|_h^2+\left\|\widetilde{E}^{k+1}\right\|_h^2\right)
-\left(\left\|E^{k}\right\|_h^2+\left\|\widetilde{E}^{k}\right\|_h^2\right)\right]
+
\left(\beta_1-C_4-C_{12}-C_{19}\right)\left\|\left|\delta_h^{\frac{\alpha}{2}}{E}^{k+1/2}\right|\right\|^2
\\
&+\left(\beta_2-C_{13}-C_{16}-C_{18}\right)\left\|\left|\delta_h^{\frac{\alpha}{2}}\widetilde{E}\right|\right\|^2
\leq
\left(\frac{1}{2}+\gamma_1+C_5+C_{14}+C_{21}\right)\left\|{E}^{k+1/2}\right\|_h^2\\
&+\left(\frac{1}{2}+\gamma_2+C_{15}+C_{17}+C_{20}\right)\left\|\widetilde{E}^{k+1/2}\right\|_{h}^2
+\frac{1}{2}\left(\left\|R^{k+1/2}\right\|_h^2+\left\|\widetilde{R}^{k+1/2}\right\|_h^2\right).
\end{aligned}
\end{equation}

Selecting appropriate parameters, such that
\begin{equation*}
\begin{aligned}\displaystyle
C_4+C_{12}+C_{19}\leq\beta_1,\;
C_{13}+C_{16}+C_{18}\leq\beta_2,
\end{aligned}
\end{equation*}
then (\ref{eq.3.75}) is further degenerated into
\begin{equation}\label{eq.3.76}
\begin{aligned}\displaystyle
&\frac{1}{2\tau}\left[\left(\left\|E^{k+1}\right\|_h^2+\left\|\widetilde{E}^{k+1}\right\|_h^2\right)
-\left(\left\|E^{k}\right\|_h^2+\left\|\widetilde{E}^{k}\right\|_h^2\right)\right]\\
&\leq
C\left(\gamma_1,\gamma_2\right)\left(\left\|{E}^{k+1/2}\right\|_h^2+
\left\|\widetilde{E}^{k+1/2}\right\|_{h}^2\right)
+\frac{1}{2}\left(\left\|R^{k+1/2}\right\|_h^2+\left\|\widetilde{R}^{k+1/2}\right\|_h^2\right)\\
&\leq\frac{1}{2}C\left(\gamma_1,\gamma_2\right)\left(\left\|{E}^{k+1}\right\|_h^2+
\left\|\widetilde{E}^{k+1}\right\|_{h}^2
+\left\|{E}^{k}\right\|_h^2+
\left\|\widetilde{E}^{k}\right\|_{h}^2
\right)
+\frac{1}{2}\left(\left\|R^{k+1/2}\right\|_h^2+\left\|\widetilde{R}^{k+1/2}\right\|_h^2\right),
\end{aligned}
\end{equation}
where $C\left(\gamma_1,\gamma_2\right)=\max\left\{\frac{1}{2}+\gamma_1+C_5+C_{14}+C_{21},
\frac{1}{2}+\gamma_2+C_{15}+C_{17}+C_{20}\right\}$.

If $\tau\leq\frac{1}{2C\left(\gamma_1,\gamma_2\right)}$, then from (\ref{eq.3.76}) can gives
\begin{equation*}
\begin{aligned}\displaystyle
\left\|E^{k+1}\right\|_h^2+\left\|\widetilde{E}^{k+1}\right\|_h^2
\leq&\left[1+4\tau C\left(\gamma_1,\gamma_2\right)\right]
\left(\left\|E^{k}\right\|_h^2+\left\|\widetilde{E}^{k}\right\|_h^2\right)\\&
+2\tau\left(\left\|R^{k+1/2}\right\|_h^2+\left\|\widetilde{R}^{k+1/2}\right\|_h^2\right).
\end{aligned}
\end{equation*}

Combining with (\ref{eq.3.48}) and using discrete Gronwall's inequality, we finally get
\begin{equation*}
\begin{aligned}\displaystyle
\left\|E^{k+1}\right\|_h^2+\left\|\widetilde{E}^{k+1}\right\|_h^2
\leq2T(b-a)\left(C_R^2+\widetilde{C}_{\widetilde{R}}^2\right)
\exp\left(8C\left(\gamma_1,\gamma_2\right)T\right)\left(\tau^2+h^4\right)^2,
\end{aligned}
\end{equation*}
i.e.,
\begin{equation*}
\begin{aligned}\displaystyle
\left\|E^{k+1}\right\|_h+\left\|\widetilde{E}^{k+1}\right\|_h&\leq
\exp\left(4C\left(\gamma_1,\gamma_2\right)T\right)\sqrt{ T(b-a)\left(C_R^2+\widetilde{C}_{\widetilde{R}}^2\right)}
\left(\tau^2+h^4\right)\\&=:C\left(\tau^2+h^4\right),\;\;0\leq k\leq N-1.
\end{aligned}
\end{equation*}
Hence, we finish the proof.
\end{proof}

\section{Numerical experiment}

\subsection{Design of iterative algorithm}
It is not difficult to find that the difference scheme (\ref{eq.3.23})
 with (\ref{eq.3.24}) and (\ref{eq.3.25}) is an implicit system. To implement the numerical
experiment, we usually need to construct an iterative
 algorithm to realize it. In this section, we will design an iterative algorithm for system (\ref{eq.3.23})
 with (\ref{eq.3.24}) and (\ref{eq.3.25}),
 and also discuss the convergence of the designed iterative algorithm in detail.
In fact, assuming that  $\left\{\left(U_j^k,V_j^k\right)|0\leq k\leq N,0\leq j\leq M\right\}$
 has been determined, if
  $$\left\{\left(U_j^{k+1/2},V_j^{k+1/2}\right)\,\Big|\,0\leq k\leq N,0\leq j\leq M\right\}$$
 can be obtained from system (\ref{eq.3.23})
 with (\ref{eq.3.24}) and (\ref{eq.3.25}), then we naturally get the values
 $$\left\{\left(U_j^{k+1},V_j^{k+1}\right)\,\Big|\,0\leq k\leq N,0\leq j\leq M\right\}$$  basing
 on the following formulas
 \begin{equation*}
\begin{aligned}\displaystyle
U_j^{k+1}=2U_j^{k+1/2}-U_j^{k},\;V_j^{k+1}=2V_j^{k+1/2}-V_j^{k}.
\end{aligned}
\end{equation*}
 In order to achieve this goal smoothly, we first rewrite system (\ref{eq.3.23}) as
\begin{equation}\label{eq.4.1}
\left\{
\begin{aligned}\displaystyle
\mathcal{H}_h^{\alpha}U_j^{k+1/2}
=&\mathcal{H}_h^{\alpha}U_j^{k}+\frac{\tau}{2}\left\{\left[\left(\beta_1+\mathrm{i}\eta_1\right)
\mathcal{A}_h^{\alpha}+\gamma_1\mathcal{H}_h^{\alpha}\right]U_j^{k+1/2}\right.\\
&\left.-\left(\mu_1+\mathrm{i}\zeta_1\right)\mathcal{H}_h^{\alpha}\left|U_j^{k+1/2}\right|^2U_j^{k+1/2}
+\mathrm{i}\mathcal{H}_h^{\alpha}\left|U_j^{k+1/2}\right|^2V_j^{k+1/2}\right\},
\\
\mathcal{H}_h^{\alpha}V_j^{k+1/2}
=&\mathcal{H}_h^{\alpha}V_j^{k}+\frac{\tau}{2}\left\{\left[\left(\beta_2+\mathrm{i}\eta_2\right)
\mathcal{A}_h^{\alpha}+\gamma_2\mathcal{H}_h^{\alpha}\right]V_j^{k+1/2}\right.\\
&\left.
-\left(\mu_2+\mathrm{i}\zeta_2\right)\mathcal{H}_h^{\alpha}\left|V_j^{k+1/2}\right|^2V_j^{k+1/2}
+\mathrm{i}\mathcal{H}_h^{\alpha}\left|V_j^{k+1/2}\right|^2U_j^{k+1/2}\right\}.
\end{aligned}
\right.
\end{equation}

Based on this, we construct an iterative algorithm to compute the
 implicit difference system (\ref{eq.3.23})
 with (\ref{eq.3.24}) and (\ref{eq.3.25}), which is similar to that in \cite{Wang}.
 The specific form is as follows
\begin{equation}\label{eq.4.2}
\left\{
\begin{aligned}\displaystyle
U_j^{k+1/2,\left(n+1\right)}
=&U_j^{k}+\frac{\tau}{2}\left\{\left[\left(\beta_1+\mathrm{i}\eta_1\right)\left(\mathcal{H}_h^{\alpha}\right)^{-1}
\mathcal{A}_h^{\alpha}+
\gamma_1\right]U_j^{k+1/2,\left(n+1\right)}\right.\\
&\left.\left.-\left(\mu_1+\mathrm{i}\zeta_1\right)
\left|U_j^{k+1/2,\left(n\right)}\right|^2U_j^{k+1/2,\left(n\right)}
+\mathrm{i}\left|U_j^{k+1/2,\left(n\right)}\right|^2V_j^{k+1/2,\left(n\right)}\right]\right\},
\\
V_j^{k+1/2,\left(n+1\right)}
=&V_j^{k}+\frac{\tau}{2}\left\{\left[\left(\beta_2+\mathrm{i}\eta_2\right)\left(\mathcal{H}_h^{\alpha}\right)^{-1}
\mathcal{A}_h^{\alpha}+
\gamma_2\right]V_j^{k+1/2,\left(n+1\right)}\right.\\
&\left.\left.
-\left(\mu_2+\mathrm{i}\zeta_2\right)
\left|V_j^{k+1/2,\left(n\right)}\right|^2V_j^{k+1/2,\left(n\right)}
+\mathrm{i}\left|V_j^{k+1/2,\left(n\right)}\right|^2U_j^{k+1/2,\left(n\right)}\right]\right\},\\
&\hspace{3.5cm}\;j\in\mathds{Z}_M,\;0\leq k\leq N-1,\;n=0,1,2\ldots,\\
&U_j^{k+1/2,\left(n\right)}=V_j^{k+1/2,\left(n\right)}=0,\;
j\in{\mathds{Z}}\backslash\mathds{Z}_M,\;0\leq k\leq N,\;n=0,1,2\ldots,
\end{aligned}
\right.
\end{equation}
with the initial iteration values
\begin{equation}\label{eq.4.3}
U_j^{k+1/2,\left(0\right)}=
\left\{
\begin{aligned}\displaystyle
&U_j^{0}+\frac{\tau}{2}\left[\left(\beta_1+\mathrm{i}\eta_1\right)\left(\mathcal{H}_h^{\alpha}\right)^{-1}
\mathcal{A}_h^{\alpha}U_j^{0}+
\gamma_1U_j^{0}-\left(\mu_1+\mathrm{i}\zeta_1\right)
\left|U_j^{0}\right|^2U_j^{0}\right.\\
&\left.\hspace{0.5cm}+\mathrm{i}\left|U_j^{0}\right|^2V_j^{0}\right],\;k=0,
\\
&\frac{3}{2}U_j^k-\frac{1}{2}U_j^{k-1},\;k\geq1,
\end{aligned}
\right.
\end{equation}
and
\begin{equation}\label{eq.4.4}
V_j^{k+1/2,\left(0\right)}=
\left\{
\begin{aligned}\displaystyle
&V_j^{0}+\frac{\tau}{2}\left[\left(\beta_2+\mathrm{i}\eta_2\right)\left(\mathcal{H}_h^{\alpha}\right)^{-1}
\mathcal{A}_h^{\alpha}V_j^{0}+
\gamma_2V_j^{0}-\left(\mu_2+\mathrm{i}\zeta_2\right)
\left|V_j^{0}\right|^2V_j^{0}\right.\\
&\left.\hspace{0.5cm}+\mathrm{i}\left|V_j^{0}\right|^2U_j^{0}\right],\;k=0,
\\
&\frac{3}{2}V_j^k-\frac{1}{2}V_j^{k-1},\;k\geq1.
\end{aligned}
\right.
\end{equation}

Next, we present the convergence results of iterative algorithm (\ref{eq.4.2})
 with (\ref{eq.4.3}) and (\ref{eq.4.4}).

\begin{theorem}
Let
\begin{equation*}
\begin{aligned}\displaystyle
Y_j^{(n)}=U_j^{k+1/2}-U_j^{k+1/2,\left(n\right)},\;
\widetilde{Y}_j^{(n)}=V_j^{k+1/2}-V_j^{k+1/2,\left(n\right)},\;
j\in\mathds{Z}_M,\;n=0,1,2\ldots,
\end{aligned}
\end{equation*}
 then for sufficiently small $\tau$ and $h$, iterative algorithm (\ref{eq.4.2})
 with (\ref{eq.4.3}) and (\ref{eq.4.4}) is convergent, i.e.
 \begin{equation*}
\begin{aligned}\displaystyle
\left\|Y_j^{(n)}\right\|_h+\left\|\widetilde{Y}_j^{(n)}\right\|_h
\rightarrow0,\;\mathrm{when}\; n\rightarrow\infty.
\end{aligned}
\end{equation*}
\end{theorem}

\begin{proof}Here, we will take two steps to prove the conclusion.

{\bf Step 1: First, we consider the case of $n=0$}. When $k=0$, combine (\ref{eq.4.2})
 with (\ref{eq.4.3}) and
 notice that equation (\ref{eq.1}), then we have
\begin{equation}\label{eq.4.5}
\begin{aligned}\displaystyle
Y_j^{(0)}=&U_j^{1/2}-U_j^{1/2,\left(0\right)}\\
=&\frac{1}{2}\left(U_j^{1}+U_j^{0}\right)
-\left\{U_j^{0}+\frac{\tau}{2}\left[\left(\beta_1+\mathrm{i}\eta_1\right)\left(\mathcal{H}_h^{\alpha}\right)^{-1}
\mathcal{A}_h^{\alpha}U_j^{0}+
\gamma_1U_j^{0}\right.\right.
\\&\left.\left.
-\left(\mu_1+\mathrm{i}\zeta_1\right)
\left|U_j^{0}\right|^2U_j^{0}+\mathrm{i}\left|U_j^{0}\right|^2V_j^{0}\right]
\right\}\\
=&\frac{1}{2}\left[U_j^{1}-u\left(x_j,t_{1}\right)\right]+
\frac{1}{2}\left[u\left(x_j,t_{1}\right)-u\left(x_j,t_{0}\right)\right]
\\&-\frac{\tau}{2}\left[\left(\beta_1+\mathrm{i}\eta_1\right)\left(\mathcal{H}_h^{\alpha}\right)^{-1}
\mathcal{A}_h^{\alpha}U_j^{0}+
\gamma_1U_j^{0}
-\left(\mu_1+\mathrm{i}\zeta_1\right)
\left|U_j^{0}\right|^2U_j^{0}+\mathrm{i}\left|U_j^{0}\right|^2V_j^{0}\right]\\
=&\mathcal{O}\left(\tau^2+h^4\right)+\mathcal{O}\left(\tau^2\right)+
\frac{\tau}{2}\partial_t u\left(x_j,t_0\right)-
\frac{\tau}{2}\left[\left(\beta_1+\mathrm{i}\eta_1\right)\left(\mathcal{H}_h^{\alpha}\right)^{-1}
\mathcal{A}_h^{\alpha}U_j^{0}\right.
\\&\left.
+
\gamma_1U_j^{0}
-\left(\mu_1+\mathrm{i}\zeta_1\right)
\left|U_j^{0}\right|^2U_j^{0}+\mathrm{i}\left|U_j^{0}\right|^2V_j^{0}\right]\\
=&\mathcal{O}\left(\tau^2+h^4\right)+\mathcal{O}\left(\tau^2\right)-\frac{\tau}{2}
\left(\beta_1+\mathrm{i}\eta_1\right)
\left[\left(-\triangle\right)^\frac{\alpha}{2}u\left(x_j,t_{0}\right)+
\left(\mathcal{H}_h^{\alpha}\right)^{-1}
\mathcal{A}_h^{\alpha}u\left(x_j,t_{0}\right)\right]\\
=&\mathcal{O}\left(\tau^2+h^4\right)+\mathcal{O}\left(\tau^2\right)
+\mathcal{O}\left(\tau h^4\right)\\
=&\mathcal{O}\left(\tau^2+h^4\right),
\end{aligned}
\end{equation}
 If $k\geq1$, then we also have
\begin{equation}\label{eq.4.6}
\begin{aligned}\displaystyle
Y_j^{(0)}=&\frac{1}{2}\left(U_j^{k+1}+U_j^{k}\right)
-\left(\frac{3}{2}U_j^k-\frac{1}{2}U_j^{k-1}\right)
=\frac{1}{2}\left(U_j^{k+1}-2U_j^{k}+U_j^{k-1}\right)
\\
=&\mathcal{O}\left(\tau^2+h^4\right).
\end{aligned}
\end{equation}
Using the same method, we can easily get
\begin{equation}\label{eq.4.7}
\widetilde{Y}_j^{(0)}=
\left\{
\begin{aligned}\displaystyle
&\mathcal{O}\left(\tau^2+h^4\right),\;k=0,\\
&\mathcal{O}\left(\tau^2+h^4\right),\;k\geq1.
\end{aligned}
\right.
\end{equation}

Combining (\ref{eq.4.5}), (\ref{eq.4.6}) and (\ref{eq.4.7}), we know that there exists a positive constant $ C_{22}$,
 such that
 \begin{equation}\label{eq.4.8}
\begin{aligned}\displaystyle
\left\|Y_j^{(0)}\right\|_h^2+\left\|\widetilde{Y}_j^{(0)}\right\|_h^2
\leq C_{22}\left(\tau^2+h^4\right)^2.
\end{aligned}
\end{equation}

{\bf Step 2: Next, we consider the case of $n\geq1$}.
Subtracting (\ref{eq.4.2}) from (\ref{eq.4.1}), we can obtain
\begin{equation}\label{eq.4.9}
\left\{
\begin{aligned}\displaystyle
Y_j^{\left(n+1\right)}
=&\frac{\tau}{2}\left[\left(\beta_1+\mathrm{i}\eta_1\right)\left(\mathcal{H}_h^{\alpha}\right)^{-1}
\mathcal{A}_h^{\alpha}+
\gamma_1\right]Y_j^{\left(n+1\right)}
+J_j^{\left(n\right)}+P_j^{\left(n\right)},
\\
\widetilde{Y}_j^{\left(n+1\right)}
=&\frac{\tau}{2}\left[\left(\beta_2+\mathrm{i}\eta_2\right)\left(\mathcal{H}_h^{\alpha}\right)^{-1}
\mathcal{A}_h^{\alpha}+
\gamma_2\right]\widetilde{Y}_j^{\left(n+1\right)}
+\widetilde{J}_j^{\left(n\right)}+\widetilde{P}_j^{\left(n\right)},\\
&\hspace{6cm}\;j\in\mathds{Z}_M,\;n=1,2\ldots.
\end{aligned}
\right.
\end{equation}
where
\begin{equation*}
\begin{aligned}\displaystyle
&J_j^{\left(n\right)}=\frac{\tau}{2}\left(\mu_1+\mathrm{i}\zeta_1\right)
\left[\left|U_j^{k+1/2,\left(n\right)}\right|^2U_j^{k+1/2,\left(n\right)}
-\left|U_j^{k+1/2}\right|^2U_j^{k+1/2}\right],
\\
&P_j^{\left(n\right)}=\frac{\tau}{2}\mathrm{i}\left[\left|U_j^{k+1/2}\right|^2V_j^{k+1/2}
-\left|U_j^{k+1/2,\left(n\right)}\right|^2V_j^{k+1/2,\left(n\right)}\right],\\
&\widetilde{J}_j^{\left(n\right)}=\frac{\tau}{2}\left(\mu_2+\mathrm{i}\zeta_2\right)
\left[\left|V_j^{k+1/2,\left(n\right)}\right|^2V_j^{k+1/2,\left(n\right)}
-\left|V_j^{k+1/2}\right|^2V_j^{k+1/2}\right],
\\
&\widetilde{P}_j^{\left(n\right)}=\frac{\tau}{2}\mathrm{i}\left[\left|V_j^{k+1/2}\right|^2U_j^{k+1/2}
-\left|V_j^{k+1/2,\left(n\right)}\right|^2U_j^{k+1/2,\left(n\right)}\right],
\end{aligned}
\end{equation*}

According to Lemma \ref{Le.2.4}, we get
\begin{equation}\label{eq.4.10}
\begin{aligned}\displaystyle
\left|J_j^{\left(n\right)}\right|&\leq
\tau\sqrt{\mu_1^2+\zeta_1^2}
\max\left\{2\left|U_j^{k+1/2,\left(n\right)}\right|^2,
\left|U_j^{k+1/2,\left(n\right)}\right|\left|U_j^{k+1/2}\right|,
\left|U_j^{k+1/2}\right|^2\right\}\left|Y_j^{(n)}\right|\\
&=:\tau C_{23}\left|Y_j^{(n)}\right|,
\end{aligned}
\end{equation}
\begin{equation}\label{eq.4.11}
\begin{aligned}\displaystyle
\left|P_j^{\left(n\right)}\right|\leq&
\frac{\tau}{2}\max\left\{2\left|U_j^{k+1/2}\right|\left|V_j^{k+1/2}\right|,
2\left|U_j^{k+1/2,\left(n\right)}\right|\left|V_j^{k+1/2}\right|,
\left|U_j^{k+1/2,\left(n\right)}\right|^2
\right\}\left(\left|Y_j^{(n)}\right|+\left|\widetilde{Y}_j^{(n)}\right|\right)
\\=:&\tau C_{24}\left(\left|Y_j^{(n)}\right|+\left|\widetilde{Y}_j^{(n)}\right|\right),
\end{aligned}
\end{equation}
\begin{equation}\label{eq.4.12}
\begin{aligned}\displaystyle
\left|\widetilde{J}_j^{\left(n\right)}\right|
&\leq \tau \sqrt{\mu_2^2+\zeta_2^2}
\max\left\{2\left|V_j^{k+1/2,\left(n\right)}\right|^2,
\left|V_j^{k+1/2,\left(n\right)}\right|\left|V_j^{k+1/2}\right|,
\left|V_j^{k+1/2}\right|^2\right\}\left|\widetilde{Y}_j^{(n)}\right|\\
&=:\tau C_{25}\left|\widetilde{Y}_j^{(n)}\right|,
\end{aligned}
\end{equation}
\begin{equation}\label{eq.4.13}
\begin{aligned}\displaystyle
\left|\widetilde{P}_j^{\left(n\right)}\right|\leq&
\frac{\tau}{2}\max\left\{2\left|V_j^{k+1/2}\right|\left|U_j^{k+1/2}\right|,
2\left|V_j^{k+1/2,\left(n\right)}\right|\left|U_j^{k+1/2}\right|,
\left|V_j^{k+1/2,\left(n\right)}\right|^2
\right\}\left(\left|Y_j^{(n)}\right|+\left|\widetilde{Y}_j^{(n)}\right|\right)
\\=:&\tau C_{26}\left(\left|Y_j^{(n)}\right|+\left|\widetilde{Y}_j^{(n)}\right|\right).
\end{aligned}
\end{equation}

Multiplying the first equation of (\ref{eq.4.9}) by $h\left(Y_j^{\left(n+1\right)}\right)^{*}$,
summing $j$ from $1$ to $M-1$ and taking the real part lead to
\begin{equation}\label{eq.4.14}
\begin{aligned}\displaystyle
\left\|Y_j^{\left(n+1\right)}\right\|_h^2
=&-\frac{\tau}{2}\beta_1
\left\|\left|\delta_h^{\frac{\alpha}{2}}Y_j^{\left(n+1\right)}\right|\right\|^2
+\frac{\tau}{2}\gamma_1\left\|Y_j^{\left(n+1\right)}\right\|_h^2
+S_1+S_2,
\end{aligned}
\end{equation}
where
\begin{equation*}
\begin{aligned}\displaystyle
S_1=\Re\left\{h\sum_{j=1}^{M-1}J_j^{\left(n\right)}
\left(Y_j^{\left(n+1\right)}\right)^{*}
\right\},\;S_2=
\Re\left\{h\sum_{j=1}^{M-1}P_j^{\left(n\right)}
\left(Y_j^{\left(n+1\right)}\right)^{*}
\right\}.
\end{aligned}
\end{equation*}
In the same way, we also have
\begin{equation}\label{eq.4.15}
\begin{aligned}\displaystyle
\left\|\widetilde{Y}_j^{\left(n+1\right)}\right\|_h^2
=&-\frac{\tau}{2}\beta_2
\left\|\left|\delta_h^{\frac{\alpha}{2}}\widetilde{Y}_j^{\left(n+1\right)}\right|\right\|^2
+\frac{\tau}{2}\gamma_2\left\|\widetilde{Y}_j^{\left(n+1\right)}\right\|_h^2
+\widetilde{S}_1+\widetilde{S}_2,
\end{aligned}
\end{equation}
where
\begin{equation*}
\begin{aligned}\displaystyle
\widetilde{S}_1=\Re\left\{h\sum_{j=1}^{M-1}\widetilde{J}_j^{\left(n\right)}
\left(\widetilde{Y}_j^{\left(n+1\right)}\right)^{*}
\right\},\;
\widetilde{S}_2=
\Re\left\{h\sum_{j=1}^{M-1}\widetilde{P}_j^{\left(n\right)}
\left(\widetilde{Y}_j^{\left(n+1\right)}\right)^{*}
\right\}.
\end{aligned}
\end{equation*}
By using the Young's inequality and combining with (\ref{eq.4.10})--(\ref{eq.4.13}), we can conclude that
\begin{equation}\label{eq.4.16}
\begin{aligned}\displaystyle
S_1+S_2\leq& h\sum_{j=1}^{M-1}\left(\left|J_j^{\left(n\right)}\right|
+\left|P_j^{\left(n\right)}\right|\right)
\left|\left(Y_j^{\left(n+1\right)}\right)^{*}\right|\\
\leq &h\sum_{j=1}^{M-1}\left[\epsilon_{12}^2\left(\left|J_j^{\left(n\right)}\right|
+\left|P_j^{\left(n\right)}\right|\right)^2+
\frac{1}{4\epsilon_{12}^2}\left|\left(Y_j^{\left(n+1\right)}\right)^{*}\right|\right]\\
\leq &h\sum_{j=1}^{M-1}\left[2\tau^2\epsilon_{12}^2
\max\left\{C_{23}^2,C_{24}^2\right\}\left(\left|Y_j^{(n)}\right|+\left|\widetilde{Y}_j^{(n)}\right|\right)^2+
\frac{1}{4\epsilon_{12}^2}\left|\left(Y_j^{\left(n+1\right)}\right)^{*}\right|\right]\\
\leq&4\tau^2\epsilon_{12}^2
\max\left\{C_{23}^2,C_{24}^2\right\}\left(\left\|Y^{(n)}\right\|_h^2
+\left\|\widetilde{Y}^{(n)}\right\|_h^2\right)+\frac{1}{4\epsilon_{12}^2}
\left\|Y^{\left(n+1\right)}\right\|_h^2,
\end{aligned}
\end{equation}
and
\begin{equation}\label{eq.4.17}
\begin{aligned}\displaystyle
\widetilde{S}_1+\widetilde{S}_2\leq4\tau^2\epsilon_{13}^2
\max\left\{C_{25}^2,C_{26}^2\right\}\left(\left\|Y^{(n)}\right\|_h^2
+\left\|\widetilde{Y}^{(n)}\right\|_h^2\right)+\frac{1}{4\epsilon_{12}^2}
\left\|\widetilde{Y}^{\left(n+1\right)}\right\|_h^2.
\end{aligned}
\end{equation}

Taking $\tau\gamma_1+\frac{1}{2\epsilon_{12}^2}=1$ and
$\tau\gamma_2+\frac{1}{2\epsilon_{13}^2}=1$.
Adding (\ref{eq.4.14}) and (\ref{eq.4.15}) and using (\ref{eq.4.16}) and (\ref{eq.4.17}) can lead to
\begin{equation}\label{eq.4.18}
\begin{aligned}\displaystyle
\left\|Y_j^{\left(n+1\right)}\right\|_h^2
+
\left\|\widetilde{Y}_j^{\left(n+1\right)}\right\|_h^2
\leq&8\tau^2C_{27}^2
\left(\left\|Y^{(n)}\right\|_h^2
+\left\|\widetilde{Y}^{(n)}\right\|_h^2\right),\;m=0,1,\ldots.
\end{aligned}
\end{equation}
where $C_{27}^2
=\max\left\{\epsilon_{12}^2
\max\left\{C_{23}^2,C_{24}^2\right\},
\epsilon_{13}^2
\max\left\{C_{25}^2,C_{26}^2\right\}\right\}$.

Consequently, according to (\ref{eq.4.18}) and noticing (\ref{eq.4.8}), when $\tau\leq2/C_{27}$ and
$n\rightarrow\infty$, we get
\begin{equation*}
\begin{aligned}\displaystyle
\left\|Y_j^{\left(n\right)}\right\|_h^2
+
\left\|\widetilde{Y}_j^{\left(n\right)}\right\|_h^2
\leq\left(\frac{1}{2}\right)^{n-1}\left(\left\|Y^{(0)}\right\|_h^2
+\left\|\widetilde{Y}^{(0)}\right\|_h^2\right)\rightarrow0,\;n=1,2\ldots.
\end{aligned}
\end{equation*}
This implies that the iterative algorithm (\ref{eq.4.2})
 with (\ref{eq.4.3}) and (\ref{eq.4.4}) that we constructed is convergent when $\tau$ is small
enough. This ends the proof.
\end{proof}

\subsection{Numerical results}
In this section,  some numerical results are presented to demonstrate the
effectiveness of the constructed numerical differential formula for the fractional Laplacian and
the numerical scheme for the coupled space fractional Ginzburg-Landau equations.

\begin{example}\label{Ex:1}
(Test the accuracy of numerical differential equation (\ref{eq.3.14})).
We choose the function $u(x) =x^4(1-x)^4,\;\;0\leq x\leq1$ as the test function.
\end{example}

 Here, we calculate the exact value of the fractional Laplacian of $u$ at $x=0.5$ and
 the approximate value by formula (\ref{eq.3.14}).
 The computation errors and accuracy for different step sizes
  $h$ and different $\alpha\in(1,2]$ are shown in Table \ref{Tab.1}.
From the numerical results in this table (the last column), it is not difficult to verify that
 formula (\ref{eq.3.14}) is of fourth-order accuracy.

\begin{table}[htbp]\renewcommand\arraystretch{1.30}
 \begin{center}
 \caption{ The absolute errors and corresponding accuracy of Example \ref{Ex:1}
 with different $h$ and $\alpha$ by the numerical differential
 formula (\ref{eq.3.14}) .}\label{Tab.1}
 \begin{footnotesize}
\begin{tabular}{c c c c c c }\hline
  $\alpha$ &\; $h$\;&  \textrm{The absolute errors}&\;\;\;\;$\textrm {The corresponding accuracy}$\vspace{0.1cm}\\
  \hline
 $1.2$& $1/200$ &1.326188e-009\;\; &  ---\\
  $$  & ${1}/{220}$& 9.145181e-010	 \;\;&  3.8995\\
  $$& ${1}/{240}$ & 6.508301e-010 \;\; &   		3.9092 \\
  $$& ${1}/{260}$ & 4.756608e-010	\;\; &  	3.9172\\
   $$& ${1}/{280}$ &3.556359e-010\;\; &  	3.9240\\ \hline
  $1.4 $& ${1}/{200}$ &2.482153e-009 \;\; &  ---\\
  $$  & ${1}/{220}$&1.702516e-009		 \;\;&  3.9557\\
  $$& ${1}/{240}$ & 1.206293e-009 \;\; &  	3.9599 \\
  $$& ${1}/{260}$ &	8.783753e-010\;\; &  	3.9633\\
   $$& ${1}/{280}$ &6.546772e-010 \;\; &  3.9663\\ \hline
 $1.6 $& ${1}/{200}$ &3.208787e-009\;\; &  ---\\
  $$  & ${1}/{220}$&2.195965e-009 \;\;&   3.9793\\
  $$& ${1}/{240}$ & 1.553023e-009 \;\; & 3.9813 \\
  $$& ${1}/{260}$ & 1.129083e-009\;\; &   3.9828\\
   $$& ${1}/{280}$ &8.403958e-010	 \;\; &  3.9846\\ \hline
   $1.8 $& ${1}/{200}$ &3.083128e-009 \;\; &  ---\\
  $$  & ${1}/{220}$&2.107216e-009 \;\;&   3.9930\\
  $$& ${1}/{240}$ &	1.488610e-009 \;\; &  	3.9940\\
  $$& ${1}/{260}$ &1.081272e-009		\;\; &   	3.9942\\
   $$& ${1}/{280}$ &8.041375e-010	 \;\; & 3.9958\\ \hline
   $2$& ${1}/{200}$ &1.874849e-009		 \;\; &  ---\\
  $$  & ${1}/{220}$& 1.280677e-009	 \;\;& 3.9989\\
  $$& ${1}/{240}$ &9.041987e-010 \;\; &   4.0006 \\
  $$& ${1}/{260}$ &6.564869e-010	\;\; &   3.9997\\
   $$& ${1}/{280}$ &4.879529e-010 \;\; &  	4.0034\\ \hline
\end{tabular}
 \end{footnotesize}
 \end{center}
 \end{table}

\begin{example}\label{Ex:2}
(Check the convergence order of the constructed difference scheme (\ref{eq.3.23})
 with (\ref{eq.3.24}) and (\ref{eq.3.25})). We consider the system
 (\ref{eq.1}) with the following initial conditions
\begin{eqnarray*}
\begin{array}{lll}
\displaystyle u\left(x,0\right)={e}^{-8}\left(1-x^2\right)^2
,\;\;
v\left(x,0\right)={e}^{-8}\left(1-x^2\right)^2,\;x\in[-1,1].
\end{array}
\end{eqnarray*}\vspace{-0.4cm}
\end{example}

Here, without loss of generality, we set
the iterative
tolerance in iterative algorithm (\ref{eq.4.2})
 with (\ref{eq.4.3}) and (\ref{eq.4.4}) as $10^{-14}$, and the final time $T=1$.
As obtaining the exact solution of the system is not feasible, we usually resort to
a ``numerical exact'' solution by using a very small steps $\tau$ and $h$.
The related calculation results of our method are listed in Tables
\ref{Tab.2} and \ref{Tab.3},  respectively.
These numerical results show that the difference scheme (\ref{eq.3.23})
 with (\ref{eq.3.24}) and (\ref{eq.3.25}) has convergence
 order $\mathcal{O}\left(\tau^2+h^4\right)$,
which further confirms the theoretical results in Theorem \ref{Th.3.21}.

\begin{table}[htbp]\renewcommand\arraystretch{1.3}
 \begin{center}
 \caption{
 The $l_h^2$-norm errors and related convergence orders of the Example \ref{Ex:2}
 for $\left(\beta_1,\beta_2\right)=\left(10^{-2},10^{-2}\right)$, $\left(\eta_1,\eta_2\right)=\left(10^{-2},10^{-2}\right)$,
 $\left(\mu_1,\mu_2\right)=\left(1,1\right)$, $\left(\zeta_1,\zeta_2\right)=\left(10^{-2},10^{-2}\right)$ and
 $\left(\gamma_1,\gamma_2\right)=\left(2/23,2/23\right)$.}\label{Tab.2}\vspace{-0.4cm}
 \begin{footnotesize}
\begin{tabular}{c c c c c c }\hline\vspace{-0.4cm}
   & \;\;\;\; \;\; &\;\;\;\;\;\;&\;\;\;\;\;\;&\;\;\;\;\;\;&\\\vspace{0.1cm}
  $\alpha$ &\;$\tau,\;h$\;&  \;The $l_h^2$-norm errors
  \; &\; Temporal convergence order\; & Spatial convergence order\\\hline
  $1.1$& $\tau=1/{5},h={1}/{5}$ &      4.157790e-002   	&  ---&  ---\\
  $$ & $\tau={1}/{20},h={1}/{10}$ &      2.463945e-003    		&       2.0384&  	     4.0768\\ \hline
$1.2$& $\tau=1/{5},h={1}/{5}$ &       4.175023e-002  	&  ---&  ---\\
  $$ & $\tau={1}/{20},h={1}/{10}$ &    2.473321e-003     		&      2.0386&  	       4.0773\\ \hline
$1.3$& $\tau=1/{5},h={1}/{5}$ &       4.188363e-002    	&  ---&  ---\\
  $$ & $\tau={1}/{20},h={1}/{10}$ &     2.481693e-003    		&     2.0385&  	      4.0770\\ \hline
   $1.4$& $\tau=1/{5},h={1}/{5}$ &      4.198741e-002   	&  ---&  ---\\
  $$ & $\tau={1}/{20},h={1}/{10}$ &     2.489342e-003    		&      2.0381&  	      4.0761\\ \hline
 $1.5$& $\tau=1/{5},h={1}/{5}$ &       4.206859e-002  	&  ---&  ---\\
  $$ & $\tau={1}/{20},h={1}/{10}$ &    2.496501e-003    		&      2.0374&  	      4.0748\\ \hline
  $1.6$& $\tau=1/{5},h={1}/{5}$ &      4.213249e-002   	&  ---&  ---\\
  $$ & $\tau={1}/{20},h={1}/{10}$ &    2.503368e-003     		&      2.0365&  	       4.0730\\ \hline
    $1.7$& $\tau=1/{5},h={1}/{5}$ &       4.218318e-002    	&  ---&  ---\\
  $$ & $\tau={1}/{20},h={1}/{10}$ &      2.510110e-003  		&     2.0354&  	   4.0708\\ \hline
  $1.8$& $\tau=1/{5},h={1}/{5}$ &      4.222375e-002    	&  ---&  ---\\
  $$ & $\tau={1}/{20},h={1}/{10}$ &      2.516875e-003     		&     2.0342&  	      4.0683\\ \hline
  $1.9$& $\tau=1/{5},h={1}/{5}$ &        4.225658e-002 	&  ---&  ---\\
  $$ & $\tau={1}/{20},h={1}/{10}$ &      2.523786e-003  		&    2.0503&  	     4.0655\\ \hline
 $2$& $\tau=1/{5},h={1}/{5}$ &        4.228351e-002    	&  ---&  ---\\
  $$ & $\tau={1}/{20},h={1}/{10}$ &    2.530956e-003     		&    2.0312&  	     4.0623\\ \hline
\end{tabular}
 \end{footnotesize}
 \end{center}
 \end{table}

 \begin{table}[htbp]\renewcommand\arraystretch{1.32}
 \begin{center}
 \caption{
 The $l_h^2$-norm errors and related convergence orders of the Example \ref{Ex:2}
 for $\left(\beta_1,\beta_2\right)=\left(10^{-3},10^{-3}\right)$, $\left(\eta_1,\eta_2\right)=\left(10^{-3},10^{-3}\right)$,
 $\left(\mu_1,\mu_2\right)=\left(10^{-1},10^{-1}\right)$, $\left(\zeta_1,\zeta_2\right)=\left(10^{-1},10^{-1}\right)$ and
 $\left(\gamma_1,\gamma_2\right)=\left(-480,-20\right)$.}\label{Tab.3}\vspace{-0.4cm}
 \begin{footnotesize}
\begin{tabular}{c c c c c c }\hline\vspace{-0.4cm}
   & \;\;\;\; \;\; &\;\;\;\;\;\;&\;\;\;\;\;\;&\;\;\;\;\;\;&\\\vspace{0.10cm}
  $\alpha$ &\;$\tau,\;h$\;&  \;The $l_h^2$-norm errors
  \; &\; Temporal convergence order\; & Spatial convergence order\\\hline
  $1.1$& $\tau=1/{5},h={1}/{5}$ &     5.519908e-004 	&  ---&  ---\\
  $$ & $\tau={1}/{20},h={1}/{10}$ &      3.385323e-005     		&         2.0136&  	      4.0273\\ \hline
$1.2$& $\tau=1/{5},h={1}/{5}$ &        5.519909e-004 	&  ---&  ---\\
  $$ & $\tau={1}/{20},h={1}/{10}$ &     3.385325e-005    		&      2.0136&  	         4.0273\\ \hline
$1.3$& $\tau=1/{5},h={1}/{5}$ &         5.519910e-004     	&  ---&  ---\\
  $$ & $\tau={1}/{20},h={1}/{10}$ &     3.385328e-005       		&      2.0136&  	       4.0273\\ \hline
   $1.4$& $\tau=1/{5},h={1}/{5}$ &      5.519912e-004   	&  ---&  ---\\
  $$ & $\tau={1}/{20},h={1}/{10}$ &     3.385331e-005    		&     2.0136&  	       4.0273\\ \hline
 $1.5$& $\tau=1/{5},h={1}/{5}$ &      5.519914e-004  	&  ---&  ---\\
  $$ & $\tau={1}/{20},h={1}/{10}$ &    3.385334e-005    		&        2.0136&  	       4.0273\\ \hline
  $1.6$& $\tau=1/{5},h={1}/{5}$ &      5.519915e-004     	&  ---&  ---\\
  $$ & $\tau={1}/{20},h={1}/{10}$ &   3.385338e-005       		&        2.0136&  	        4.0273\\ \hline
    $1.7$& $\tau=1/{5},h={1}/{5}$ &        5.519917e-004      	&  ---&  ---\\
  $$ & $\tau={1}/{20},h={1}/{10}$ &      3.385342e-005  		&      2.0136&  	     4.0273\\ \hline
  $1.8$& $\tau=1/{5},h={1}/{5}$ &      5.519920e-004     	&  ---&  ---\\
  $$ & $\tau={1}/{20},h={1}/{10}$ &   3.385347e-005      		&     2.0136&  	        4.0273\\ \hline
  $1.9$& $\tau=1/{5},h={1}/{5}$ &          5.519922e-004 	&  ---&  ---\\
  $$ & $\tau={1}/{20},h={1}/{10}$ &      3.385353e-005   		&     2.0136&  	       4.0273\\ \hline
 $2$& $\tau=1/{5},h={1}/{5}$ &           5.519925e-004    	&  ---&  ---\\
  $$ & $\tau={1}/{20},h={1}/{10}$ &     3.385359e-005      		&      2.0136&  	      4.0273\\ \hline
\end{tabular}
 \end{footnotesize}
 \end{center}
 \end{table}

\begin{example}\label{Ex:3}
(Observe the evolution of coupled system (\ref{eq.1})). Here, we consider system (\ref{eq.1}) with the
 following initial conditions
\begin{eqnarray*}
\begin{array}{lll}
\displaystyle u\left(x,0\right)=\sech\left(x+5\right)
{e}^{8\mathrm{i}x},\;\;
v\left(x,0\right)=\sech\left(x-5\right)
{e}^{-8\mathrm{i}x},\;\;x\in[-10,10].
\end{array}
\end{eqnarray*}\vspace{-0.4cm}
\end{example}

\begin{figure}[!htbp]
\centering
\mbox{
\subfigure
{\includegraphics[width=7.5cm, height=6.5cm]{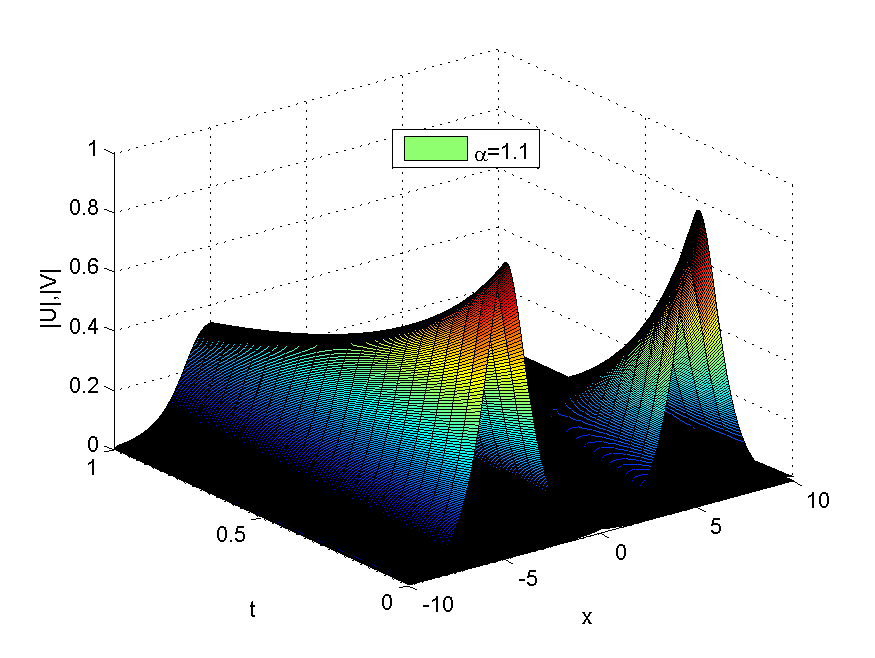}}\hspace{0.1cm}
\quad
\subfigure
{\includegraphics[width=7.5cm, height=6.5cm]{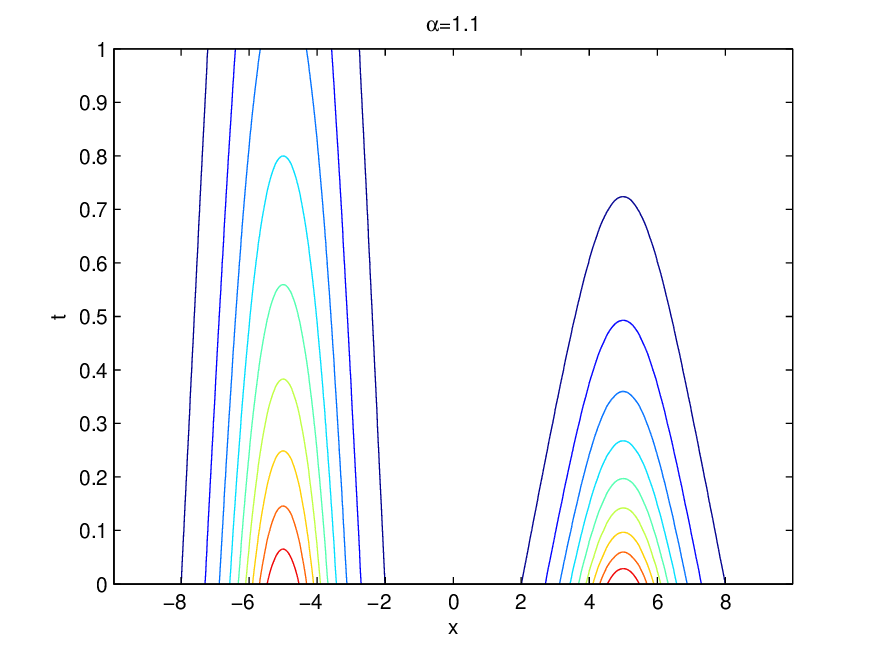}}}
\mbox{
\subfigure
{\includegraphics[width=7.5cm, height=6.5cm]{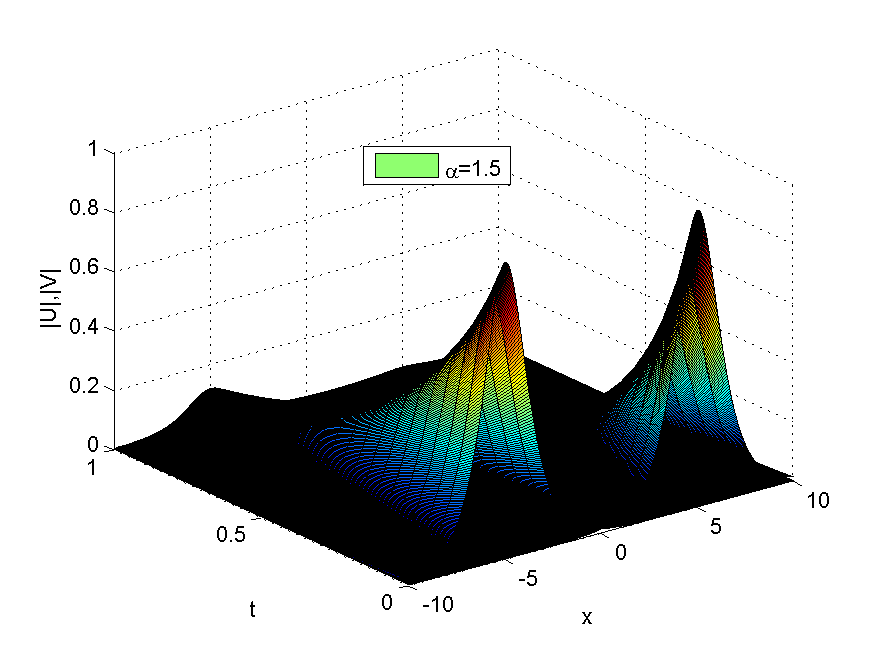}}\hspace{0.1cm}
\quad
\subfigure
{\includegraphics[width=7.5cm, height=6.5cm]{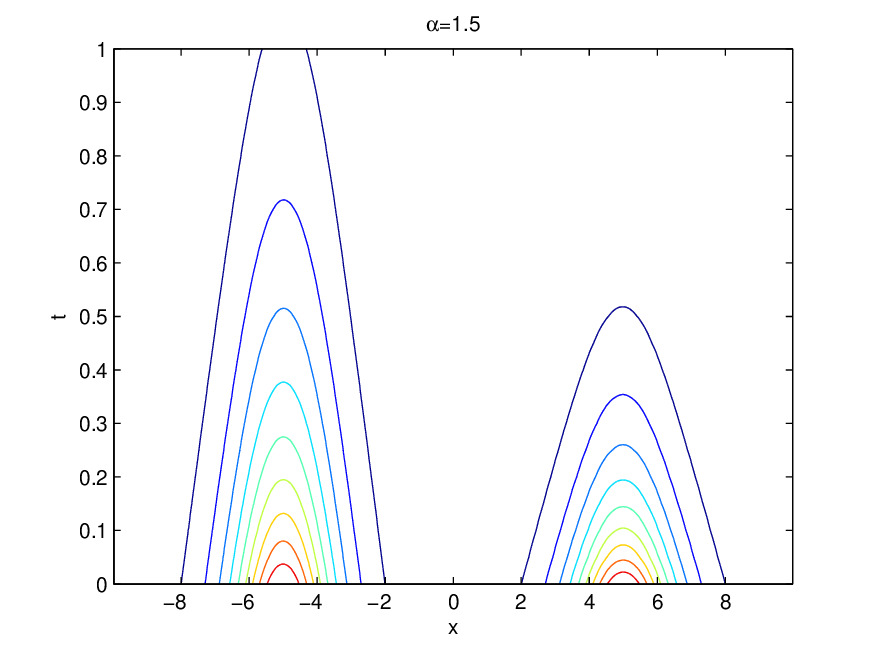}}}
\mbox{
\subfigure
{\includegraphics[width=7.5cm, height=6.5cm]{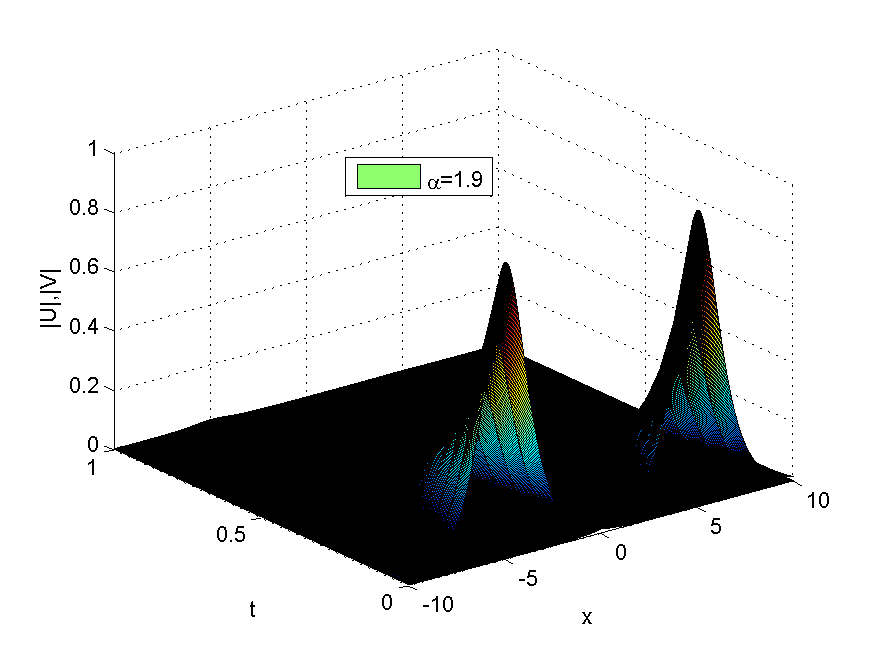}}\hspace{0.1cm}
\quad
\subfigure
{\includegraphics[width=7.5cm, height=6.5cm]{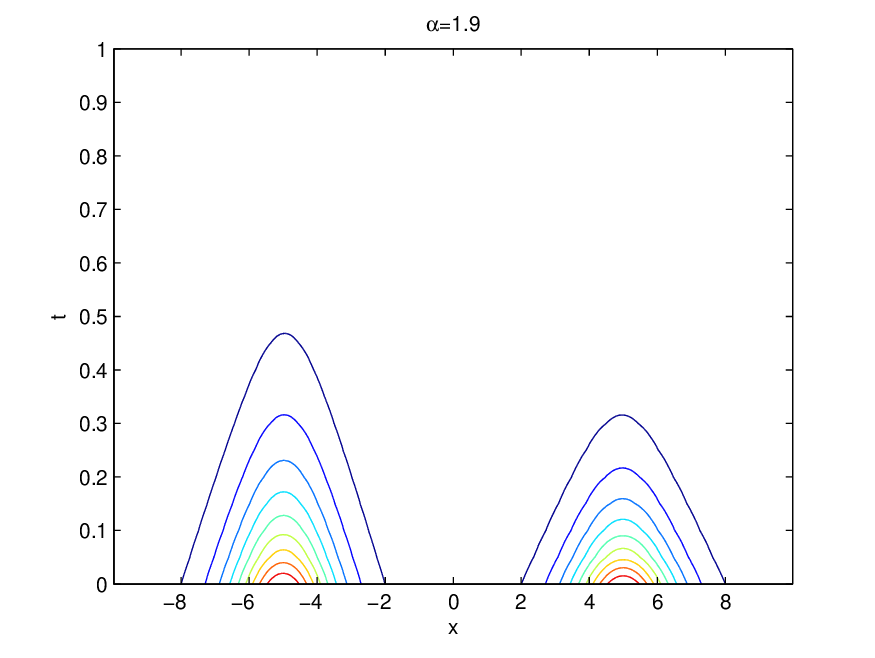}}}
\caption{The profile of the evolution of $|U|$ and $|V|$ (left), and contour plot of
$|U|$ and $|V|$ (right), for different values of $\alpha$ .}\label{Fig:7.1}
\end{figure}

\begin{figure}[!htbp]
\centering
\mbox{
\subfigure
{\includegraphics[width=7.5cm, height=8cm]{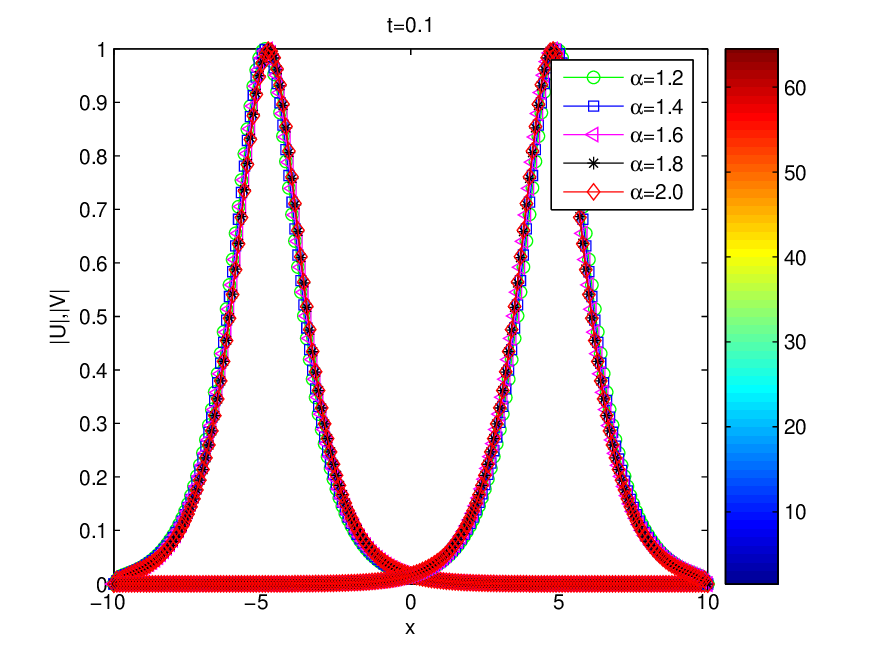}}\hspace{0.1cm}
\quad
\subfigure
{\includegraphics[width=7.5cm, height=8cm]{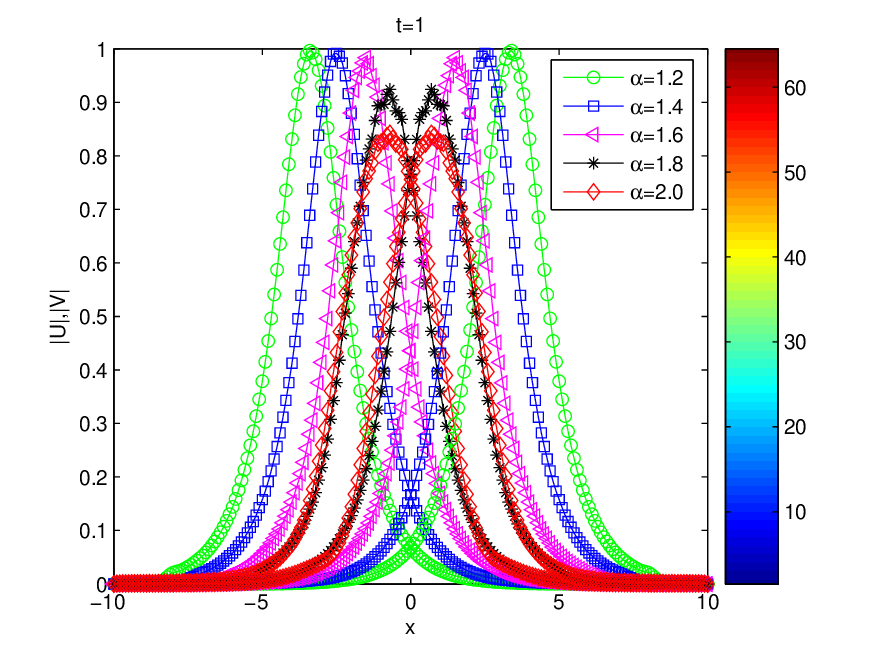}}}
\mbox{
\subfigure
{\includegraphics[width=7.5cm, height=8cm]{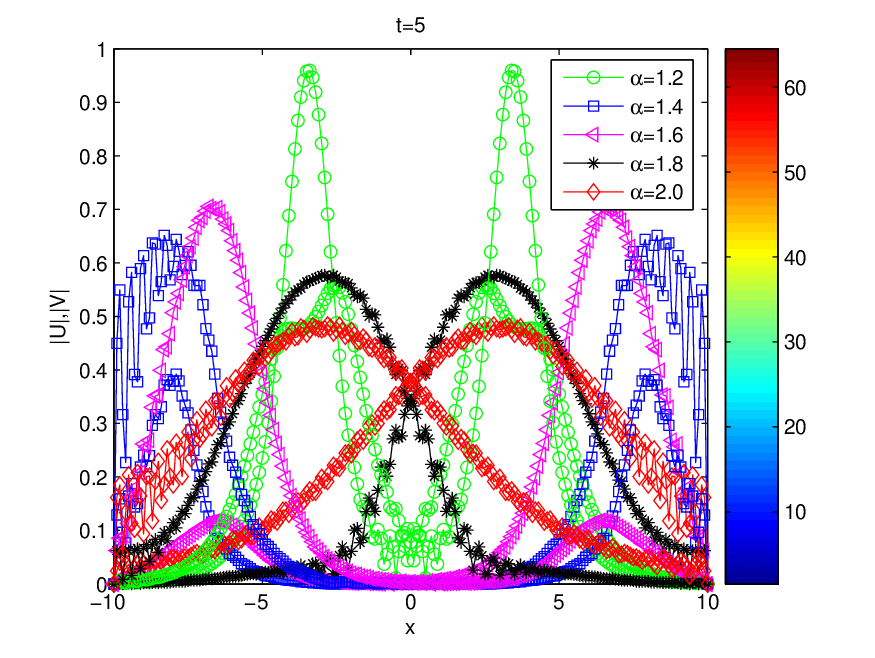}}\hspace{0.1cm}
\quad
\subfigure
{\includegraphics[width=7.5cm, height=8cm]{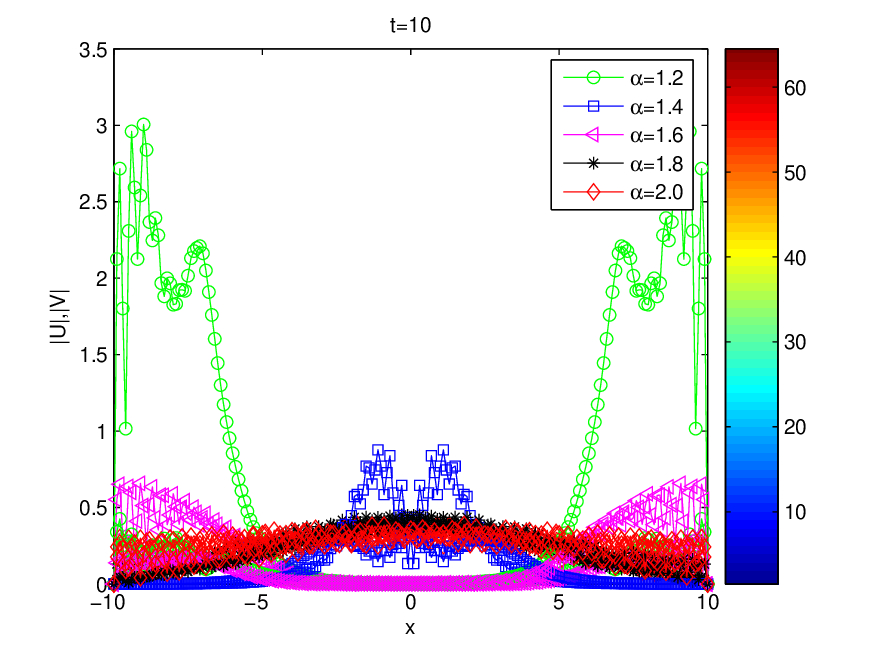}}}
\caption{Shift of profiles of $|U|$ and $|V|$ at different times $t=0.1,1,5,10$ for
 different values of $\alpha$ .}\label{Fig:7.2}
\end{figure}

Firstly, let's examine the influence of the fractional Laplacian operator on the dissipation mechanism.
For this reason, we choose the parameters $\left(\beta_1,\beta_2\right)=\left(10^{-1},10^{-1}\right)$, $\left(\eta_1,\eta_2\right)=\left(10^{-2},10^{-2}\right)$,
 $\left(\mu_1,\mu_2\right)=\left(1,1\right)$, $\left(\zeta_1,\zeta_2\right)=\left(10^{-1},10^{-1}\right)$ and
 $\left(\gamma_1,\gamma_2\right)=\left(0.25,-2\right)$. The numerical solutions with
different values of $\alpha$ are shown in Fig. \ref{Fig:7.1}. From this figure, we can easily find
that the fractional $\alpha$ has a great influence on the shape of solitons,
which is completely different from the classical case, thus further explaining
the nonlocal characteristics of fractional Laplacian.
Furthermore, in order to better show the change of numerical solution with time, using the
parameters $\left(\beta_1,\beta_2\right)=\left(10^{-3},10^{-3}\right)$, $\left(\eta_1,\eta_2\right)=\left(1,1\right)$,
 $\left(\mu_1,\mu_2\right)=\left(10^{-4},10^{-4}\right)$, $\left(\zeta_1,\zeta_2\right)=\left(10^{-1},10^{-1}\right)$ and
 $\left(\gamma_1,\gamma_2\right)=\left(10^{-2},10^{-2}\right)$, the plot of solutions at $t=0.1,1,5,10$ is presented in Fig. \ref{Fig:7.2}.
Based on the figure, we can see how the two waves travel and collide with each other.

\begin{figure}[!htbp]
\centering
\mbox{
\subfigure
{\includegraphics[width=7.5cm, height=6.5cm]{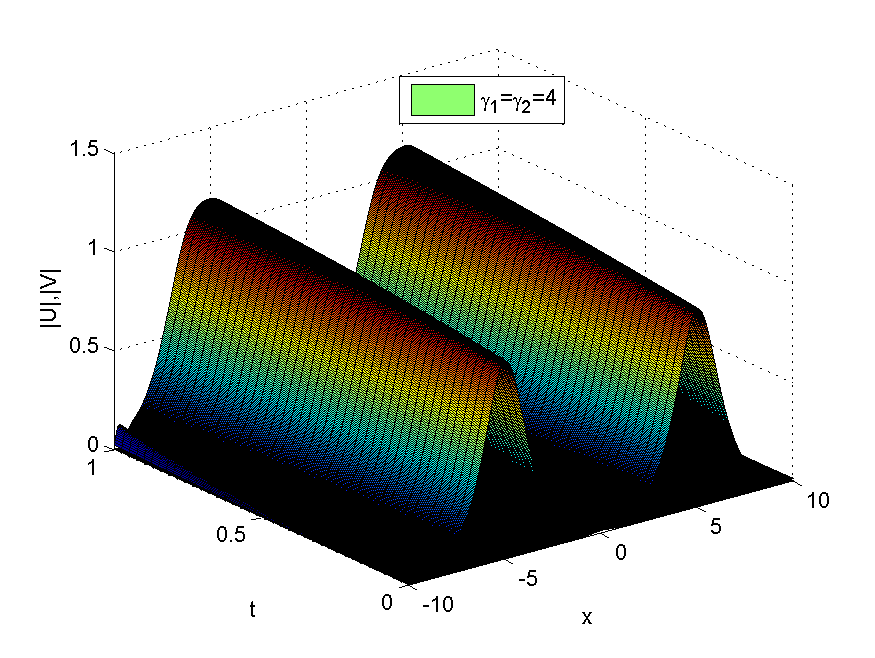}}\hspace{0.1cm}
\quad
\subfigure
{\includegraphics[width=7.5cm, height=6.5cm]{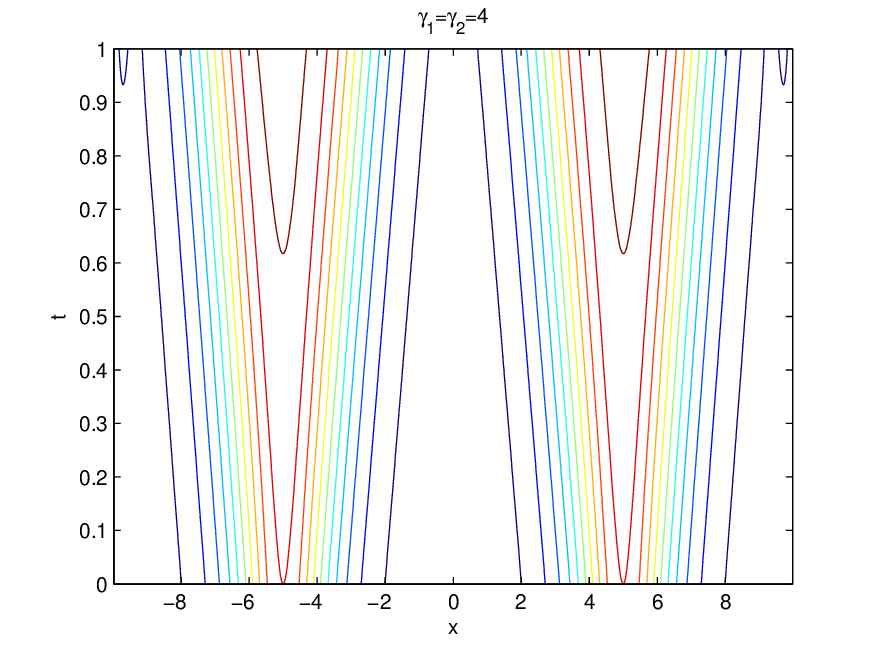}}}
\mbox{
\subfigure
{\includegraphics[width=7.5cm, height=6.5cm]{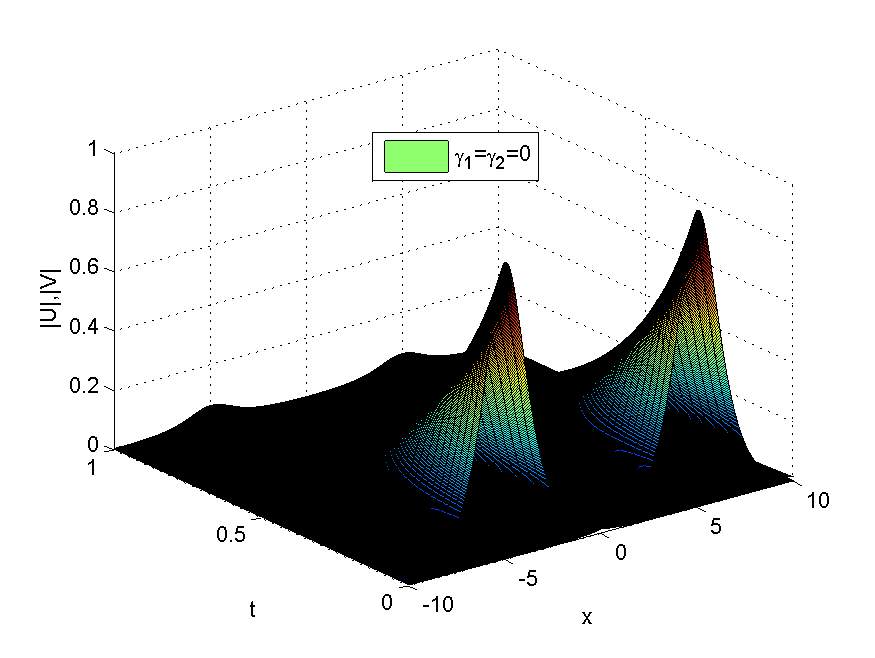}}\hspace{0.1cm}
\quad
\subfigure
{\includegraphics[width=7.5cm, height=6.5cm]{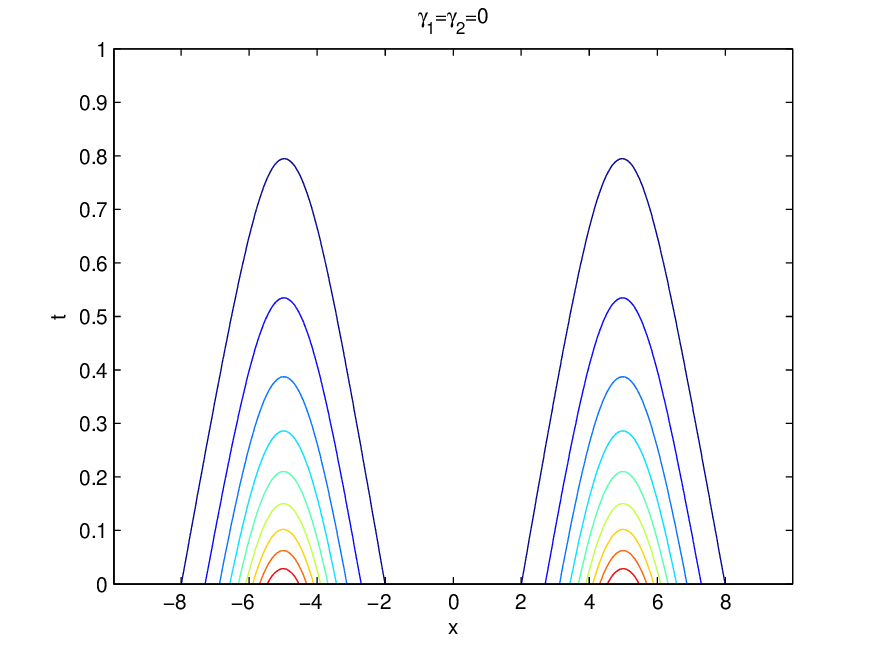}}}
\mbox{
\subfigure
{\includegraphics[width=7.5cm, height=6.5cm]{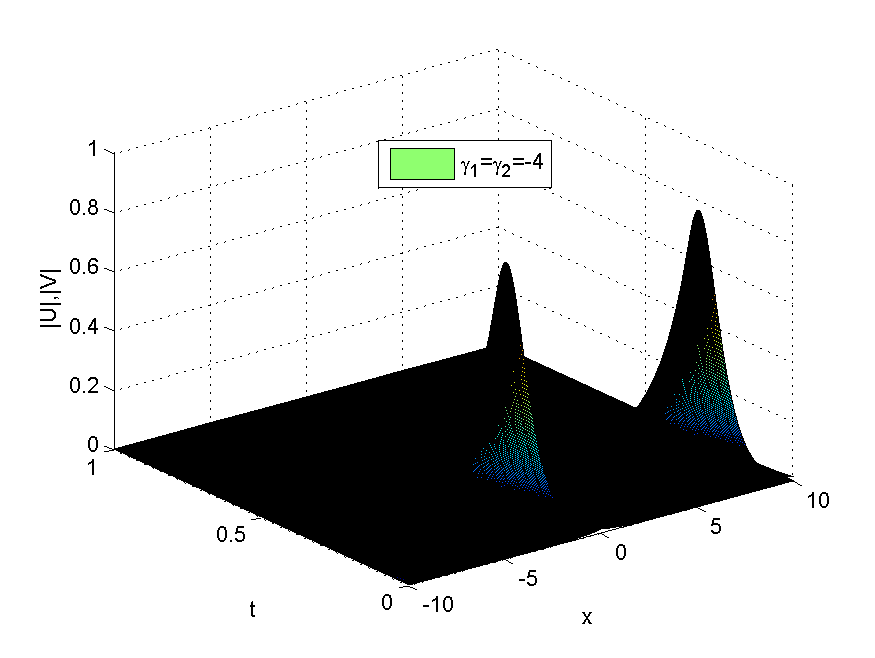}}\hspace{0.1cm}
\quad
\subfigure
{\includegraphics[width=7.5cm, height=6.5cm]{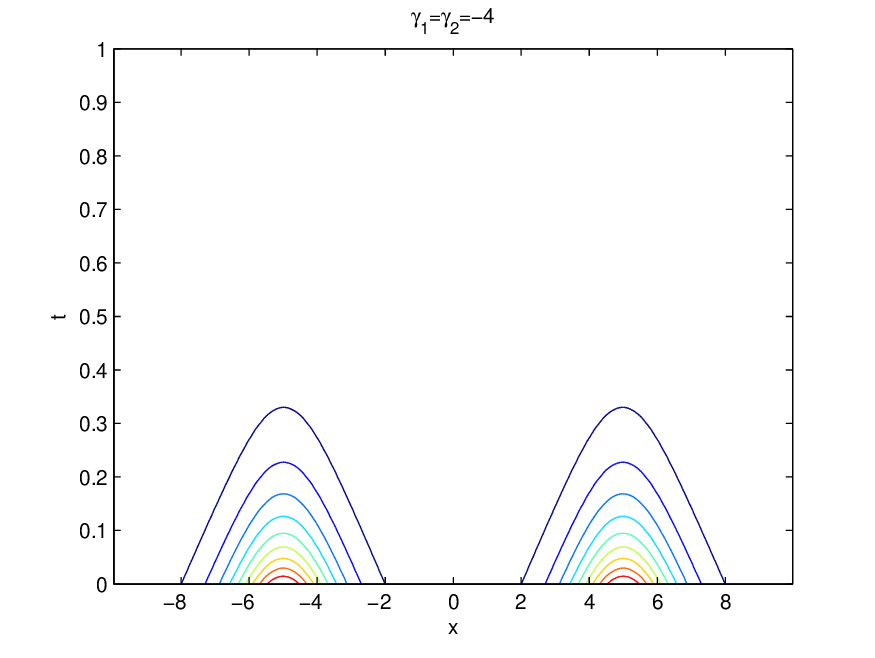}}}
\caption{Plots of position density of $|U|$ and $|V|$ (left), and contour plot of
$|U|$ and $|V|$ (right), for different $\gamma_1$ and $\gamma_2$.}\label{Fig:7.3}
\end{figure}

\begin{figure}[!htbp]
\centering
\mbox{
\subfigure
{\includegraphics[width=7.5cm, height=8cm]{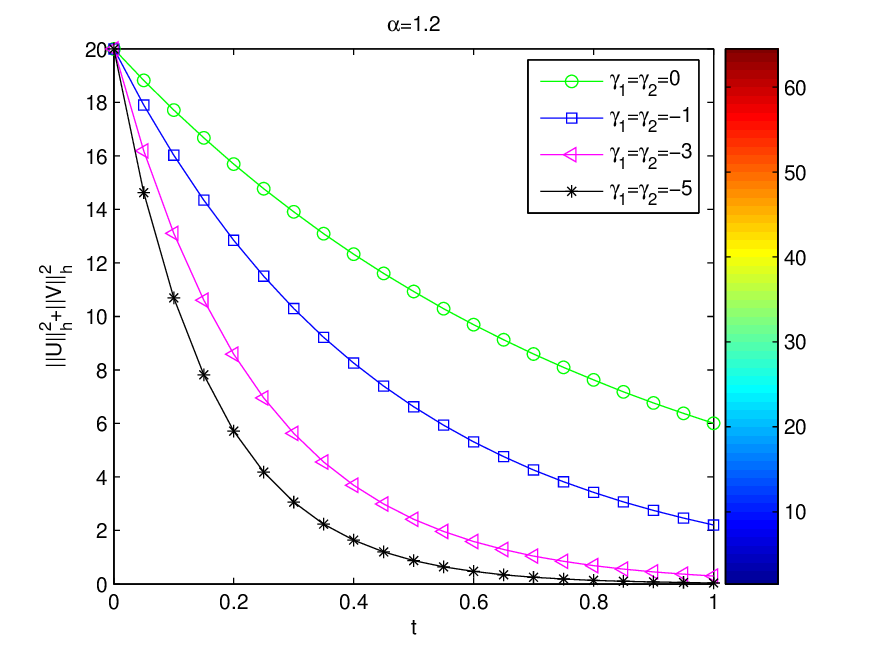}}\hspace{0.1cm}
\quad
\subfigure
{\includegraphics[width=7.5cm, height=8cm]{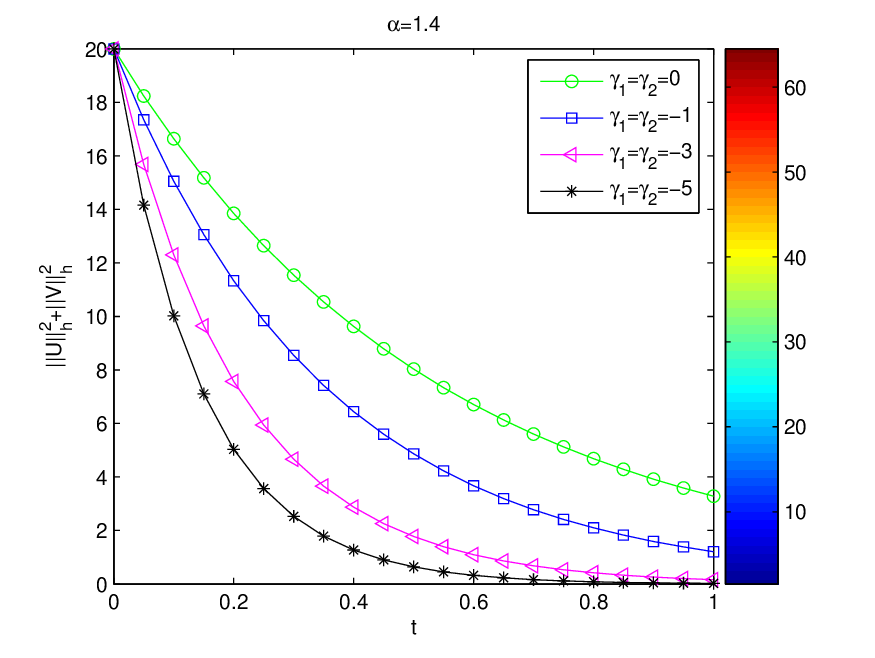}}}
\mbox{
\subfigure
{\includegraphics[width=7.5cm, height=8cm]{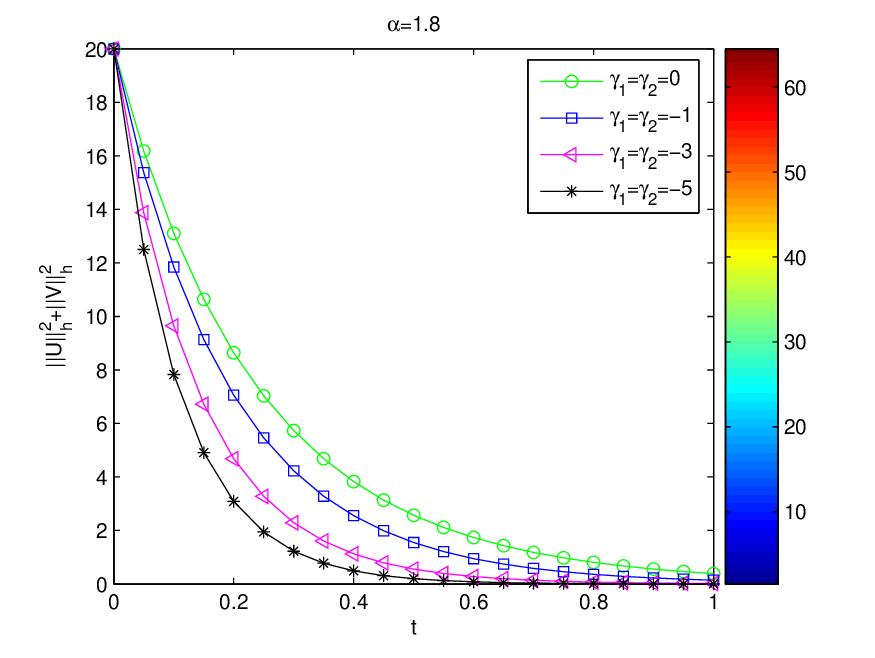}}\hspace{0.1cm}
\quad
\subfigure
{\includegraphics[width=7.5cm, height=8cm]{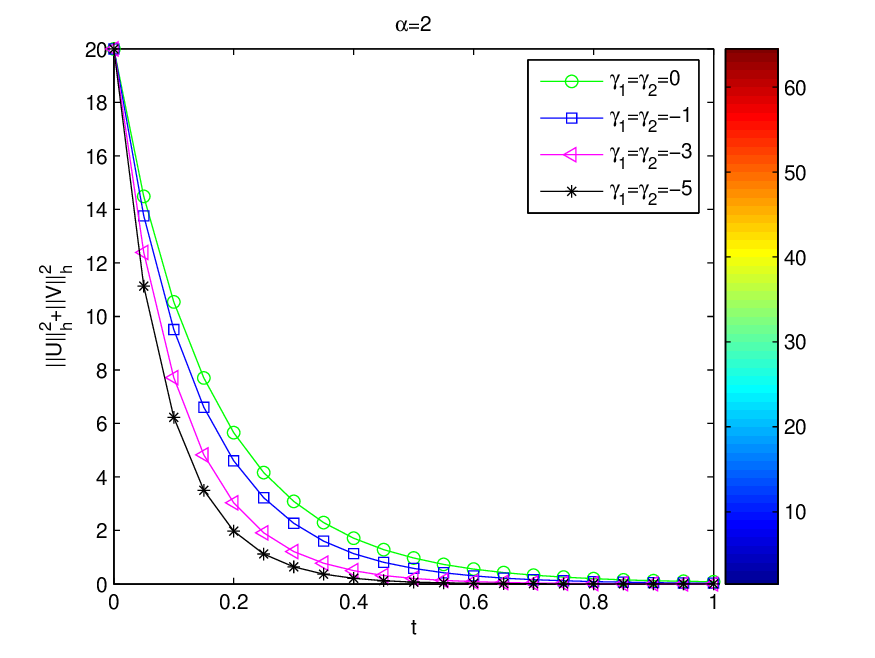}}}
\caption{The discrete norm $\left\|U\right\|_h^2+\left\|V\right\|_h^2$
with different $\gamma_1$ and $\gamma_2$ for different values of $\alpha$ .}\label{Fig:7.4}
\end{figure}

Then, we consider the influence of parameters $\gamma_1$ and $\gamma_2$ on the
 shape of wave evolution. In Fig. \ref{Fig:7.3},
we depict the evolution of $|U|$ and $|V|$ at different $\gamma_1$ and $\gamma_1$ values when parameters
$\left(\beta_1,\beta_2\right)=\left(10^{-1},10^{-1}\right)$, $\left(\eta_1,\eta_2\right)=\left(10^{-2},10^{-2}\right)$,
 $\left(\mu_1,\mu_2\right)=\left(1,1\right)$ and $\left(\zeta_1,\zeta_2\right)=\left(10^{-1},10^{-1}\right)$.
 By comparing the
changing trends of these graphs, we can claim that parameters $\gamma_1$ and $\gamma_1$
 has a significant impact on the evolution of waves.
Moreover, we also depict the evolution of $\left\|U\right\|_h^2+\left\|V\right\|_h^2$ with different $\alpha$
under the parameters
$\left(\beta_1,\beta_2\right)=\left(10^{-1},10^{-1}\right)$, $\left(\eta_1,\eta_2\right)=\left(10^{-4},10^{-4}\right)$,
 $\left(\mu_1,\mu_2\right)=\left(10^{-3},10^{-3}\right)$ and $\left(\zeta_1,\zeta_2\right)
 =\left(10^{-2},10^{-2}\right)$ in Fig. \ref{Fig:7.4}.
From this figure, we can see that as $\gamma_1$ and $\gamma_2$ become smaller,
their corresponding wave evolution decays faster.

\section{Conclusions}
In the paper, we have studied the well-posedness of the weak solution for the CNLSFGL equations
 that involves the fractional Laplacian. Then,
a fourth-order numerical differential formula is developed for the fractional Laplacian.
Using the constructed formula, we propose an effective difference scheme for the
CNLSFGL equations that is unconditionally convergent with  order of
$\mathcal{O}\left(\tau^2+h^4\right)$.
Moreover, we have
investigated the boundedness, existence and uniqueness of the numerical solution.
Finally,
an efficient iterative algorithm is introduced, and some numerical results were presented to confirm the
theoretical results.

In addition, there are three aspects worthy of further clarification: Firstly, the fourth-order
 numerical differential formula that we constructed to approximate the fractional Laplacian
 can be applied to other spatial fractional differential equations. Secondly, some methods and
 techniques in this paper can be used for reference in other spatial fractional differential
  equations. Finally, the coefficient matrix of the algebraic equations corresponding
  to the difference scheme established in this paper has a Toeplitz structure independent
  of time, and some matrix decomposition and iterative techniques can be considered
  to reduce the computation cost and storage requirements, which is our next plan.

${\mathbf{Appendix}}$

{\bf Proof of Lemma \ref{Le.3.3}.}  Here, we will discuss two scenarios.
(a)~For the case where $\alpha$ is a positive integer $n$, that is, when
$\alpha=n\in \mathds{Z}^{+}$, we can see that it is Newton's binomial expansion and is finite,
so the conclusion is obviously valid.

(b)~For the case where $\alpha$ is a non-positive integer, i.e. $\alpha\notin {\mathds{Z}}^{+}$,
we first calculate the $n$-th derivative of $G_2(z)$ as follows
\begin{equation*}
\begin{aligned}
\frac{\mathrm{d}^nG_2(z)}{\mathrm{d}z^n}=&b_0^{\alpha}\sum_{k=0}^{n}\left(n\atop k\right)
\frac{\mathrm{d}^{n-k}\left(1-z\right)^{\alpha}}{\mathrm{d}z^{n-k}}
\frac{\mathrm{d}^{k}\left(1-\frac{b_2}{b_0}z\right)^{\alpha}}{\mathrm{d}z^{k}}
\\
=&\left(-1\right)^nb_0^{\alpha}\left[
\prod_{i=0}^{n-1}\left(\alpha-i\right)\left(1-z\right)^{\alpha-n}
\left(1-\frac{b_2}{b_0}z\right)^{\alpha}+
\prod_{j=0}^{n-1}\left(\alpha-j\right)\left(1-z\right)^{\alpha}\left(\frac{b_2}{b_0}\right)^{n}
\left(1-\frac{b_2}{b_0}z\right)^{\alpha-n}\right.\\
&+\left.
\sum_{k=1}^{n-1}\left(n\atop k\right)\left(\frac{b_2}{b_0}\right)^{k}
\prod_{i=0}^{n-k-1}\left(\alpha-i\right)\prod_{j=0}^{k-1}
\left(\alpha-j\right)\left(1-z\right)^{\alpha-(n-k)}
\left(1-\frac{b_2}{b_0}z\right)^{\alpha-k}
\right],
\;n=1,2,\cdots.
\end{aligned}
\end{equation*}
From this, we can see that the Taylor expansion of function $G_2(z)$ at $z=0$ is
\begin{equation*}
\begin{aligned}G_2(z)=&\;\left(b_0+b_1z+b_2z^2\right)^{\alpha}\\=
&\;b_0^\alpha-\alpha b_0^\alpha\left(1+\frac{b_2}{b_0}\right)z+
\frac{1}{2}\alpha b_0^\alpha\left[\left(\alpha-1\right)+2\alpha\frac{b_2}{b_0}
+\left(\alpha-1\right)\left(\frac{b_2}{b_0}\right)^2\right]z^2\\
&
+\cdots+
\frac{\frac{\mathrm{d}^nG_2(0)}{\mathrm{d}z^n}}{n!}z^n
+\cdots=:\;\sum_{n=0}^{\infty} \kappa_{2,n}^{(\alpha)}z^n,
\end{aligned}
\end{equation*}
where
\begin{equation*}
\begin{aligned}
\kappa_{2,n}^{(\alpha)}=&\left(-1\right)^nb_0^\alpha
\sum_{k=0}^{n}\frac{1}{n!}\left(n\atop k\right)\left(\frac{b_2}{b_0}\right)^{k}
\prod_{i=0}^{n-k-1}\left(\alpha-i\right)\prod_{j=0}^{k-1}
\left(\alpha-j\right)
=\left(-1\right)^nb_0^\alpha
\sum_{k=0}^{n}\left(\frac{b_2}{b_0}\right)^{k}
\left(\alpha\atop k\right)
\left(\alpha\atop n-k\right).
\end{aligned}
\end{equation*}
Using the formula \cite{Ding&Li}
\begin{equation*}
\begin{aligned}
\kappa_{2,n}^{(\alpha)}\sim\frac{\sin\left(\pi\alpha\right)\Gamma\left(\alpha+1\right)}{\pi}
n^{-\alpha-1}\;\;\;\;as\;\;\;\;n\rightarrow\infty,
\end{aligned}
\end{equation*}
we easily know that
\begin{equation*}
\begin{aligned}
\lim_{n\rightarrow\infty}\left|\frac{\kappa_{2,n+1}^{(\alpha)}z^{n+1}}{\kappa_{2,n}^{(\alpha)}z^n}\right|=&
\lim_{n\rightarrow\infty}\left|\frac{\kappa_{2,n+1}^{(\alpha)}}{\kappa_{2,n}^{(\alpha)}}\right|
\left|z\right|=\lim_{n\rightarrow\infty}\left(1+\frac{1}{n}\right)^{-\alpha-1}
|z|
=|z|,
\end{aligned}
\end{equation*}
which indicates that the convergence interval of series (\ref{eq.3.6}) is $(-1, 1)$.

Next, we consider the convergence of series (\ref{eq.3.6}) in $(-1,1)$. The usual way is to consider the
limit of its Cauchy remainder
\begin{equation*}
\begin{aligned}
T_n(z)=&\frac{\frac{\mathrm{d}^{n+1}G_2(\theta z)}{\mathrm{d}z^{n+1}}}{n!}\left(1-\theta\right)^nz^{n+1}\\
=&\left[\left(-1\right)^{n+1}b_0^\alpha\sum_{k=0}^{n+1}\left(n+1\atop k\right)\left(\frac{b_2}{b_0}\right)^{k}
\prod_{i=0}^{n-k}\left(\alpha-i\right)\prod_{j=0}^{k-1}
\left(\alpha-j\right)\left(1-\theta z\right)^{\alpha-(n+1-k)}
\left(1-\frac{b_2}{b_0}\theta z\right)^{\alpha-k}\right]\\&\times\frac{(1-\theta)^n}{n!}z^{n+1}\\
=&\left[\left(-1\right)^{n+1}(n+1)b_0^\alpha\sum_{k=0}^{n+1}\left(\frac{b_2}{b_0}\right)^{k}
\left(\alpha\atop k\right)
\left(\alpha\atop n+1-k\right)\left(\frac{1-\theta z}{1-\frac{b_2}{b_0}\theta z}\right)^{k}
z^{n+1}\right]\\&\times
\left(\frac{1-\theta }{1-\theta z}\right)^n
\left(1-\theta z\right)^{\alpha-1}
\left(1-\frac{b_2}{b_0}\theta z\right)^{\alpha},\;\;\theta\in[0,1],
\end{aligned}
\end{equation*}
but we find that this formula is
extremely complicated and it is very difficult to calculate the limit.
In order to avoid directly studying the remainder, we use other
methods to reconsider the convergence problem.

For this reason,
we assume that the series (\ref{eq.3.6}) converges to the function $S(z)$ in the interval $(-1,1)$, which is
\begin{equation*}
\begin{aligned}
S(z)=\;&b_0^\alpha-\alpha b_0^\alpha\left(1+\frac{b_2}{b_0}\right)z+
\frac{1}{2}\alpha b_0^\alpha\left[\left(\alpha-1\right)+2\alpha\frac{b_2}{b_0}
+\left(\alpha-1\right)\left(\frac{b_2}{b_0}\right)^2\right]z^2\\
\;&+\cdots+
\frac{\frac{\mathrm{d}^nG_2(0)}{\mathrm{d}z^n}}{n!}z^n+\cdots,\;\;z\in(-1,1),
\end{aligned}
\end{equation*}
further calculations can lead to
\begin{equation*}
\begin{aligned}
\alpha\left(b_1+2b_2z\right)S(z)=
\alpha b_1b_0^\alpha+\alpha b_0^\alpha
\left(2b_2+\frac{\alpha b_1^2}{b_0}\right)z+\cdots
+\alpha\left(b_1\kappa_{2,n}^{(\alpha)}+2b_2\kappa_{2,n-1}^{(\alpha)}\right)z^n+\cdots,
\end{aligned}
\end{equation*}
and
\begin{equation*}
\begin{aligned}
\left(b_0+b_1z+b_2z^2\right)S'(z)=\;&\alpha b_1b_0^\alpha+\alpha b_0^\alpha
\left(2b_2+\frac{\alpha b_1^2}{b_0}\right)z+\cdots\\&
+\left[\left(n+1\right)b_0\kappa_{2,n+1}^{(\alpha)}+
n b_1\kappa_{2,n}^{(\alpha)}+\left(n-1\right)b_2\kappa_{2,n-1}^{(\alpha)}\right]z^n+\cdots.
\end{aligned}
\end{equation*}

With the help of the recurrence relation
\begin{equation*}
\begin{aligned}
\kappa_{2,n+1}^{(\alpha)}=\frac{1}{b_0\left(n+1\right)}
\left[b_1\left(\alpha-n\right)\kappa_{2,n}^{(\alpha)}+
b_2\left(2\alpha-n+1\right)\kappa_{2,n-1}^{(\alpha)}\right],
\end{aligned}
\end{equation*}
 we can get
 \begin{equation*}
\begin{aligned}
\left(n+1\right)b_0\kappa_{2,n+1}^{(\alpha)}+
nb_1\kappa_{2,n}^{(\alpha)}+\left(n-1\right)b_2\kappa_{2,n-1}^{(\alpha)}
=\alpha\left(b_1\kappa_{2,n}^{(\alpha)}+2b_2\kappa_{2,n-1}^{(\alpha)}\right),
\end{aligned}
\end{equation*}
from which we can easily find that
\begin{equation*}
\begin{aligned}
\alpha\left(b_1+2b_2z\right)S(z)=
\left(b_0+b_1z+b_2z^2\right)S'(z),
\end{aligned}
\end{equation*}
which means that
\begin{equation}\label{eq.5.1}
\begin{aligned}
\frac{S'(z)}{S(z)}=\frac{\alpha\left(b_1+2b_2z\right)}
{\left(b_0+b_1z+b_2z^2\right)}.
\end{aligned}
\end{equation}

Integrating both sides of the equation (\ref{eq.5.1}) with respect to $z$ on $(0, z)$ can further lead to
\begin{equation*}
\begin{aligned}
S(z)=\left(b_0+b_1z+b_2z^2\right)^\alpha,\;\;z\in(-1,1),
\end{aligned}
\end{equation*}
this further indicates that
\begin{equation*}
\begin{aligned}
G_2(z)=\left(b_0+b_1z+b_2z^2\right)^\alpha=\sum_{n=0}^{\infty} \kappa_{2,n}^{(\alpha)}z^n
\end{aligned}
\end{equation*}
is uniformly valid for $\alpha\geq1$ and $x\in(-1,1)$.

Finally,  we will discuss the convergence of series (\ref{eq.3.6}) at $x=\pm1$.
Let $w_n(z)=\kappa_{2,n}^{(\alpha)}z^n$, and
when $x=\pm1$ is taken, it is recorded as $w_n(z)=w_n$.
For any series $\sum_{n=0}^\infty w_n$, there is
\begin{equation*}
\begin{aligned}
\rho=\lim_{n\rightarrow\infty}n\left[\left|\frac{w_n}{w_{n+1}}\right|-1\right]
=\lim_{n\rightarrow\infty}n\left[\left(1+\frac{1}{n}\right)^{\alpha+1}-1\right]
=\alpha+1>1.
\end{aligned}
\end{equation*}
It follows from the Raabe discriminant method, we know that
$\sum_{n=0}^\infty w_n$ is absolutely convergent.
Therefore, the convergence interval of series
$\sum_{n=0}^{\infty} \kappa_{2,n}^{(\alpha)}z^n$ is $[-1,1]$.

Therefore, based on all the above facts, we can claim that
(\ref{eq.3.6})
is uniformly true for $\alpha\geq1$ and $z\in[-1,1]$. This ends the proof.\\

%
\textbf{\large Declarations}\\

\textbf{Data Availability}
  The authors confirm that the data supporting the findings of this study are available within the article.

\textbf{Conflict of interest} The authors declare no competing interests.

%
%
%
%


\begin{thebibliography}{999}

\bibitem{Podlubny}
{\sc I. Podlubny}, {\em Fractional Differential Equations}, Academic Press, San Diego, CA,
1999.

\bibitem{Metzler}
{\sc R. Metzler, J. Klafter}, {\em The random walk's guide to anomalous diffusion: a
fractional dynamics approach}, Phys. Rep. 339 (1) (2000) 1-77.

\bibitem{Hilfer}
{\sc R. Hilfer}, {\em Applications of Fractional Calculus in Physics}, World Scientific
Publishing Co., Singapore, 2000.

\bibitem{Ortigueira}
{\sc M. Ortigueira}, {\em Fractional Calculus for Scientists and Engineers}, Springer
Netherlands, 2011.

\bibitem{Samko}
{\sc S.G. Samko, A.A. Kilbas, O.I. Marichev}, {\em Fractional integrals and derivatives:
Theory and applications}, Vol. 1. Gordon and Breach Science, 1993.

\bibitem{Landkof}
{\sc N.S. Landkof}, {\em Foundations of modern potential theory},  Springer Berlin Heidelberg,1972.

\bibitem{Aranson}
{\sc I.S. Aranson, L. Kramer}, {\em The world of the complex Ginzburg-Landau equation},
Rev. Mod. Phys. 74 (2002) 99-143.

\bibitem{Du}
{\sc Q. Du, M.D. Gunzburger, J.S. Peterson}, {\em Analysis and approximation of the Ginzburg-Landau model
 of superconductivity}, SIAM Rev. 34 (1992) 54-81.


\bibitem{Ginzburg}
{\sc V.L. Ginzburg, L.D. Landau}, {\em On the theory of superconductivity}, Zh. Eksp. Teor. Fiz. 20 (1950) 1064-1082.

\bibitem{Pan}
{\sc K. Pan, X. Jin, D. He}, {\em Pointwise error estimates of a linearized difference scheme
for strongly coupled fractional Ginzburg-Landau equations}. Math.
Method. Appl. Sci. 43 (2019) 512-35

\bibitem{Tarasov1}
{\sc V. Tarasov, G. Zaslavsky}, {\em Fractional Ginzburg-Landau equation for fractal media},
Physica A. 354 (2005) 249-261.

\bibitem{Tarasov2}
{\sc V. Tarasov, G. Zaslavsky}, {\em Fractional dynamics of coupled oscillators with
long-range interaction}, Chaos. 16 (2006) 023110.





\bibitem{Mvogo}
{\sc A. Mvogo, A. Tambue, G.H. Ben-Bolie, T.C. Kofan\'{e}}, {\em Localized numerical impulse solutions in diffuse
neural networks modeled by the complex
fractional Ginzburg-Landau equation}, Commun. Nonlinear Sci.
Numer. Simul. 39 (2016), 396-410.


\bibitem{Tarasov}
{\sc V.E. Tarasov, G.M. Zaslavsky}, {\em Fractional Ginzburg-Landau equation for fractal media},
 Physica A 354 (2005), 249-61.

\bibitem{Milovanov}
{\sc A. Milovanov, J.J. Rasmussen}, {\em Fractional generalization of the Ginzburg-Landau
 equation: an unconventional approach to critical phenomena in
complex media}, Phys. Lett. A. 337 (2005), 75-80.

\bibitem{Duan}
{\sc J. Duan, E. Titi, P. Holmes}, {\em Regularity approximation and asymptotic
dynamics for a generalized Ginzburg-Landau equation}, Nonlinearity. 6 (1993), 915-933.

\bibitem{Doering}
{\sc C.R. Doering, J.D. Gibbon, C.D. Levermore}, {\em Weak and strong solutions
of the complex Ginzburg-Landau equation}, Physica D. 71 (1994), 285-318.

\bibitem{Gao}
{\sc H. Gao, G.Lin, J. Duan}, {\em Asymptotics for the generalized two-dimensional
Ginzburg-Landau equation}, J. Math. Anal. Appl. 247 (2000), 198-216.

\bibitem{Guo}
{\sc B. Guo, G. Wang, D. Li}, {\em The attractor of the stochastic generalized
Ginzburg-Landau equation}, Sci. China Ser. A. 51 (2008), 955-964.

\bibitem{Huo}
{\sc Z. Huo, Y. Jia}, {\em Global well-posedness for the generalized 2D
Ginzburg-Landau equation}, J. Differ. Equations. 247 (2009), 260-276.

\bibitem{Pu}
{\sc X. Pu, B. Guo}, {\em Well-posedness and dynamics for the fractional
 Ginzburg-Landau equation}, Appl. Anal.
92 (2013), 318-334.

\bibitem{Gu}
{\sc X.M. Gu, L. Shi, T.H. Liu}, {\em Well-posedness of the fractional Ginzburg-Landau equation},
Appl. Anal. 98 (2019), 2545-2558.

\bibitem{Millot}
{\sc V. Millot, Y. Sire}, {\em On a fractional Ginzburg-Landau equation
and 1/2-harmonic maps into spheres}, Arch. Ration.
Mech. Anal. 215 (2015), 125-210.

\bibitem{Lu}
{\sc H. Lu, S. L\"{u}, Z. Feng}, {\em Asymptotic dynamics of 2D fractional complex
Ginzburg-Landau equation}, Int. J. Bifurc.
Chaos. 23 (2013), 1350202

\bibitem{Wang}
{\sc P. Wang P, C. Huang}, {\em An implicit midpoint difference scheme for
 the fractional Ginzburg-Landau equation}, J. Comput. Phys. 312 (2026), 31-49.

\bibitem{Hao}
{\sc Z. Hao, Z.Z. Sun}, {\em  A linearized high-order difference scheme for the fractional
Ginzburg-Landau equation}, Numer. Methods Partial Differential
Equations 33(2016), 105-124.

\bibitem{He}
{\sc D. He, K Pan}, {\em An unconditionally stable linearized difference scheme for the
fractional Ginzburg-Landau equation}, Numer. Algorithms
79 (2018), 899-925.

\bibitem{Wang1}
{\sc N. Wang, C. Huang}, {\em An efficient split-step quasi-compact finite difference method
for the nonlinear fractional Ginzburg-Landau equations},
Comput. Math. Appl. 75 ( 2018), 2223-2242

\bibitem{Ding}
{\sc H. Ding, C. Li}, {\em High-order numerical algorithm and error analysis for the two-dimensional
nonlinear spatial fractional complex Ginzburg-Landau
equation}, Commun. Nonlinear Sci. Numer. Simul. 120 (2023), 107160.

\bibitem{Mohebbi}
{\sc A. Mohebbi A}, {\em Fast and high-order numerical algorithms for the solution
 of multidimensional nonlinear fractional Ginzburg-Landau equation}, Eur.
Phys. J. Plus. 133 (2018), 67.

\bibitem{Zhang1}
{\sc Q. Zhang, X. Lin, K. Pan, Y. Ren}, {\em Linearized ADI schemes for two-dimensional
 space-fractional nonlinear Ginzburg-Landau equation}, Comput. Math.
Appl. 80 (2020), 1201-1220.

\bibitem{Zhang2}
{\sc Q. Zhang, J.S. Hesthaven, Z.Z. Sun, Y. Ren}, {\em Pointwise error estimate in difference
setting for the two-dimensional nonlinear fractional complex
Ginzburg-Landau equation}, Adv. Comput. Math. 47 (2021), 35.


\bibitem{Li}
{\sc  M. Li, C. Huang, N. Wang}, {\em Galerkin finite element method for the nonlinear
 fractional Ginzburg-Landau equation}, Appl. Numer. Math.
118 (2017), 131-149.

\bibitem{Zhang3}
{\sc  Z. Zhang, M. Li, Z. Wang}, {\em A linearized Crank-Nicolson Galerkin FEMs for the
nonlinear fractional Ginzburg-Landau equation}, Appl. Anal.
98 (2018), 2648-2667.

\bibitem{Fei}
{\sc M. Fei, C. Huang, N. Wang, G. Zhang}, {\em Galerkin-Legendre spectral method for
 the nonlinear Ginzburg-Landau equation with the Riesz fractional
derivative}, Math. Methods Appl. Sci. 44(2019), 2711-2730.

\bibitem{Shokri}
{\sc  A. Shokri, M. Dehghan}, {\em A meshless method using radial basis
 functions for the numerical solution of two-dimensional complex Ginzburg-Landau
equation}, Comput. Model. Eng. Sci. 84 (2012), 333.

\bibitem{Abbaszadeh}
{\sc M. Abbaszadeh, M. Dehghan}, {\em The fourth-order time-discrete scheme
 and split-step direct meshless finite volume method for solving
 cubicquintic complex Ginzburg-Landau equations on complicated
geometries}, Eng. Comput. 38 (2020), 1543-1557.

\bibitem{Zeng}
{\sc  W. Zeng, A. Xiao, X. Li}, {\em Error estimate of Fourier pseudo-spectral method
 for multidimensional nonlinear complex fractional Ginzburg-Landau
equations} Appl. Math. Lett. 93 (2019), 40-50.

\bibitem{Lu1}
{\sc  H. Lu, S. L\"{u}, M. Zhang}, {\em Fourier spectral approximations to the dynamics
of 3D fractional complex Ginzburg-Landau equation}, Discrete Contin.
Dyn. Syst. A 37 (2017), 2539-2564.

\bibitem{Li&Huang}
{\sc M. Li, C. Huang}, {\em An efficient difference scheme for the coupled nonlinear
fractional Ginzburg-Landau equations with the fractional Laplacian},
Numer. Methods Partial Differential Equations 35 (2018), 394-421.

\bibitem{Guo1}
{\sc B. Guo, Y. Han, J. Xin}, {\em Existence of the global smooth solution to the
period boundary value problem of fractional nonlinear Schr\"{o}dinger
equation}, Appl. Math. Comput. 204 (2008), 468--477.

\bibitem{Evans}
{\sc
L. C. Evans}, {\em Partial differential equation}, 2nd edition, 1999.

\bibitem{Macias&Diaz}
{\sc J. Mac\'{\i}as-D\'{\i}az}, {\em A structure-preserving method for a class of nonlinear dissipative
wave equations with Riesz space-fractional derivatives},
J. Comput. Phys. 351 (2017), 40--58.

\bibitem{Ding1}
{\sc H. Ding, H. Qu, Q. Yi}, {\em Construction and Analysis of Structure-Preserving
Numerical Algorithm for Two-Dimensional Damped Nonlinear Space Fractional
 Schr\"{o}dinger equation}, J. Sci. Comput. 99 (2024), 60.

\bibitem{Ding2}
{\sc H. Ding, Q. Yi}, {\em The construction of higher-order numerical approximation
 formula for Riesz derivative and its application to nonlinear
 fractional differential equations (I)}, Commun. Nonlinear Sci. 110 (2022), 106394.

\bibitem{Yang}
{\sc Q. Yang , F. Liu , I. Turner}, {\em Numerical methods for fractional partial
differential equations with Riesz space fractional derivatives},
Appl. Math. Model. 34 (2010), 200--218.

\bibitem{Ding&Li}
{\sc H. Ding, C. Li}, {\em High-order numerical algorithms for Riesz derivatives
via constructing new generating functions}, J. Sci. Comput. 71 (2017), 759--784.

\bibitem{Tadjeran&Meerschaert}
{\sc C. Tadjeran, M. M. Meerschaert}, {\em A second-order accurate numerical
method for the two-dimensional fractional diffusion equation},
J. Comput. Phys. 220 (2007), 813--823.

\bibitem{Jian}
{\sc H.Y. Jian, T.Z. Huang, X.M. Gu, X.L. Zhao, Y.L. Zhao}, {\em Fast implicit integration factor method for
nonlinear space Riesz fractional reaction-diffusion equations},
 J. Comput. Appl. Math. 378 (2020), 112935.

 \bibitem{Moler}
{\sc C. Moler, C.V. Loan}, {\em Nineteen dubious ways to compute the exponential of a matrix,
twenty-five years later}, SIAM Rev.
45 (2003), 3--49.

 \bibitem{Bhayo}
{\sc B.A. Bhayo, S. J\'{o}zsef}, {\em On Jordan's and Kober's inequality}, Acta et Commentationes
Universitatis Tartuensis de Mathematica 20 (2016), 111--116.

 \bibitem{Chan}
{\sc R. Chan and X. Jin}, {\em An Introduction to Iterative Toeplitz Solvers}, SIAM, Philadelphia, 2007.

\bibitem{Akrivis}
{\sc G. D. Akrivis}, {\em Finite difference discretization of the cubic
Schr\"{o}dinger equation}, IMA J. Numer. Anal. 13 (1993), 115--124.

\bibitem{Sun}
{\sc Z. Sun, D. Zhao}, {\em On the $L^{\infty}$ convergence of a difference scheme for coupled nonlinear
Schr\"{o}dinger equations}, Comput. Math. Appl.
59 (2010), 3286--3300.









%
%
%
%
%
%
%
%
%
%
%
%
%
%
%
%
%
%
%
%
%
%
%
%
%
%
%
%
%
%
%
%
%
%
%
%
%
%
%
%
%
%
%
%
%
%
%
%
%
%
%













































%
%
%
%
%
%
%
%
%
%
%
%
%
%
%
%
%
%
%
%
%
%
%
%
%


\end{thebibliography}
\end{document}